\documentclass[12pt,letterpaper]{amsart}
\usepackage{amsmath,amsthm,amsfonts,amssymb,amscd,amstext}
\usepackage[dvips]{graphicx}
\usepackage{psfrag,euscript,a4wide}
\usepackage{pslatex}

\newcommand{\quand}{\quad\text{and}\quad}

\theoremstyle{plain}
\newtheorem{maintheorem}{Theorem}
\newtheorem{maincorollary}{Corollary}

\newtheorem{theorem}{Theorem}[section]
\newtheorem{proposition}[theorem]{Proposition}
\newtheorem{corollary}[theorem]{Corollary}
\newtheorem{lemma}[theorem]{Lemma}

\theoremstyle{definition}
\newtheorem{remark}[theorem]{Remark}
\newtheorem{definition}[theorem]{Definition}

\renewcommand{\angle}{\sphericalangle}

\newcommand{\re}{{\mathbb R}}
\newcommand{\nat}{{\mathbb N}}

\newcommand{\var}{\varphi}

\newcommand{\de}{\delta}

\newcommand{\cal}{\mathcal}
\newcommand{\real}{{\mathbb R}}
\newcommand{\vr}{\varphi}
\newcommand{\vep}{\varepsilon}
\newcommand{\diam}{\operatorname{diam}}
\newcommand{\esup}{\operatorname{ess\,sup}}
\newcommand{\osc}{\operatorname{osc}}
\renewcommand{\epsilon}{\varepsilon}
\newcommand{\dist}{\operatorname{dist}}

\newcommand{\interior}{\operatorname{int}}

\newcommand{\deter}{\operatorname{det}}
\newcommand{\Leb}{\operatorname{Leb}}
\newcommand{\si}{\operatorname{Sing}}

\newcommand{\supp}{\operatorname{supp}} 
 
\newcommand{\reff}[1]{(\ref{#1})} 
 
\newcommand{\inte}{\operatorname{int}}

\newcommand{\suporte}{\operatorname{supp}}

\newcommand{\cT}{\EuScript{T}}
\newcommand{\cF}{\EuScript{F}}
\newcommand{\cE}{\EuScript{E}}
\newcommand{\F}{\EuScript{F}}
\newcommand{\cP}{\EuScript{P}}
\newcommand{\cO}{\EuScript{O}}
\newcommand{\V}{\EuScript{V}}

\title[Singular-hyperbolic attractors are chaotic]
{Singular-hyperbolic attractors are chaotic}
\author{V. Araujo, M. J. Pacifico, E. R. Pujals, M. Viana}
\begin{thanks} {V.A. was partially supported by CMUP-FCT
    (Portugal), CNPq (Brazil) and grants BPD/16082/2004 and
    POCI/MAT/61237/2004 (FCT-Portugal) while enjoying a
    post-doctorate leave from CMUP at PUC-Rio and IMPA.
    M.J.P., E.R.P. and M.V. were partially supported by
    PRONEX, CNPq and FAPERJ-Brazil. }
\end{thanks}

\begin{document}

\address{V\'\i tor Ara\'ujo, 
  Instituto de Matem\'a\-tica,
Universidade Federal do Rio de Janeiro,
C. P. 68.530, 21.945-970, Rio de Janeiro, RJ-Brazil
\emph{and} 
Centro de Matem\'atica da
  Universidade do Porto, Rua do Campo Alegre 687, 4169-007
  Porto, Portugal
}
\email{vitor.araujo@im.ufrj.br \emph{and} vdaraujo@fc.up.pt}

\address{Maria Jos\'e Pacifico, Instituto de Matem\'atica,
Universidade Federal do Rio de Janeiro,
C. P. 68.530, 21.945-970 Rio de Janeiro, Brazil}
\email{pacifico@im.ufrj.br \emph{and} pacifico@impa.br}

\address{Enrique R. Pujals,  IMPA, Estrada D. Castorina 110,
22460-320 Rio de Janeiro, Brazil}
\email{enrique@impa.br}

\address{Marcelo Viana, IMPA, Estrada D. Castorina 110,
22460-320 Rio de Janeiro, Brazil}
\email{viana@impa.br}

\begin{abstract}
  We prove that a singular-hyperbolic 
  attractor of a $3$-dimen\-sional flow is chaotic, in two
  strong different senses.  Firstly, the flow is expansive:
  if two points remain close for all times, possibly with
  time reparametrization, then their orbits coincide.
  Secondly, there exists a physical (or Sinai-Ruelle-Bowen)
  measure supported on the attractor whose ergodic basin
  covers a full Lebesgue (volume) measure subset of the
  topological basin of attraction. Moreover this measure has
  absolutely continuous conditional measures along the
  center-unstable direction, is a $u$-Gibbs
  state and an equilibrium state for the logarithm of the
  Jacobian of the time one map of the flow along the
  strong-unstable direction.  
  
  This extends to the class of singular-hyperbolic
  attractors the main elements of the ergodic theory of
  uniformly hyperbolic (or Axiom A) attractors for flows.
  
  In particular these results can be applied (i) to the flow
  defined by the Lorenz equations, (ii) to the geometric
  Lorenz flows, (iii) to the attractors appearing in the
  unfolding of certain resonant double homoclinic loops,
  (iv) in the unfolding of certain singular cycles and (v)
  in some geometrical models which are singular-hyperbolic
  but of a different topological type from the geometric
  Lorenz models. In all these cases the results show that
  these attractors are expansive and have physical measures
  which are $u$-Gibbs states.
\end{abstract}


\renewcommand{\subjclassname}{\textup{2000} Mathematics Subject Classification}


\maketitle



\section{Introduction}

The theory of uniformly hyperbolic dynamics was initiated in the
1960's by Smale~\cite{Sm67} and, through the work of his students
and collaborators, as well as mathematicians in the Russian
school, immediately led to extraordinary development of the whole
field of Dynamical Systems. However, despite its great successes,
this theory left out important classes of dynamical systems, which
do not conform with the basic assumptions of uniform
hyperbolicity. The most influential examples of such systems are,
arguably, the H\'enon map~\cite{He76}, for the discrete time case,
and the Lorenz flow~\cite{Lo63}, for the continuous time case.

The Lorenz equations highlighted, in a striking way, the fact that
for continuous time systems, robust dynamics may occur outside the
realm of uniform hyperbolicity and, indeed, in the presence of
equilibria that are accumulated by recurrent periodic orbits. This
prompted the quest for an extension of the notion of uniform
hyperbolicity encompassing all continuous time systems with robust
dynamical behavior. A fundamental step was carried out by Morales,
Pacifico, Pujals~\cite{MPP98,MPP04}, who proved that a robust
invariant attractor of a $3$-dimensional flow that contains some
equilibrium must be \emph{singular hyperbolic}, that is, it must
admit an invariant splitting $E^s\oplus E^{cu}$ of the tangent
bundle into a $1$-dimensional uniformly contracting sub-bundle and
a $2$-dimensional volume-expanding sub-bundle.

In fact, Morales, Pacifico, Pujals proved that any robust
invariant \emph{set} of a $3$-dimensional flow containing some
equilibrium is a singular hyperbolic attractor or repeller. In the
absence of equilibria, robustness implies uniform hyperbolicity.
The first examples of singular hyperbolic sets included the Lorenz
attractor \cite{Lo63,Tu99} and its geometric models
\cite{Gu76,ABS77,GW79,Wil79}, and the singular-horseshoe
\cite{LP86}, besides the uniformly hyperbolic sets themselves.
Many other examples have been found recently, including attractors
arising from certain resonant double homoclinic loops \cite{MPS05}
or from certain singular cycles \cite{MPu97}, and certain models
across the boundary of uniform hyperbolicity \cite{MPP00}.

The next natural step is to try and understand what are the
dynamical consequences of singular hyperbolicity. Indeed, it is
now classical that uniform hyperbolicity has very precise
implications on the dynamics (symbolic dynamics, entropy), the
geometry (invariant foliations, fractal dimensions), the
statistics (physical measures, equilibrium states) of the
invariant set. It is important to know to what extent this remains
valid in the singular hyperbolic domain. There is substantial
advance in this direction at the topological level
\cite{MPP99,CMP00,MP01,MP03,MPP04,Mo04,BM04}, but the ergodic
theory of singular hyperbolic systems remains mostly open
(for a recent advance in the particular case of the Lorenz
attractor see~\cite{LMP05}). The
present paper is a contribution to such a theory.

Firstly, we prove that the flow on a singular hyperbolic set is
\emph{expansive}. Roughly speaking, this means that any two orbits
that remain close at all times must actually coincide. However,
the precise formulation of this property is far from obvious in
this setting of continuous time systems is far from obvious. The
definition we use here was introduced by Komuro~\cite{Km84}:
other, more naive versions, turn out to be inadequate in this
context.

Another main result, extending \cite{clm2005}, is that
typical orbits in the basin of the attractor have
well-defined statistical behavior: for Lebesgue almost every
point the forward Birkhoff time average converges, and is
given by a certain physical probability measure. We also
show that this measure admits absolutely continuous
conditional measures along the center-unstable directions on
the attractor. As a consequence, it is a $u$-Gibbs
state
and an equilibrium state for the flow.

The main technical tool for the proof of these results is a
construction of convenient cross-sections and invariant
contracting foliations for a corresponding Poincar\'e map,
reminiscent of \cite{BR75}, that allow us to reduce the flow
dynamics to certain $1$-dimensional expanding transformations.
This construction will, no doubt, be useful in further analysis of
the dynamics of singular hyperbolic flows.
   
Let us give the precise statements of these results.


\subsection{Singular-hyperbolicity} 
\label{sec:sing-hyperb}

Throughout, $M$ is a compact boundaryless 3-dimensional manifold
and ${\cal X}^1(M)$ is the set of $C^1$ vector fields on $M$,
endowed with the $C^1$ topology.
From now on we fix some smooth Riemannian structure on $M$
and an induced normalized volume form $m$ that we call
Lebesgue measure. 
We write also $\dist$ for the induced distance on $M$.
Given $X \in {\cal X}^1(M)$, we denote by $X_t$, $t \in \re$ the flow
induced by $X$, and if $x \in M$ and $[a,b]\subset \re$ then
$X_{[a,b]}(x)= \{X_t(x), a\leq t \leq b\}$.

Let $\Lambda$ be a compact invariant set  of $X\in
{\cal X}^1(M)$.  We say that $\Lambda$ is \emph{isolated} if
there exists an open set $U\supset \Lambda$ such that
$$
\Lambda =\bigcap_{t\in \re}X_t(U)
$$
If $U$ above can be chosen such that $X_t(U)\subset U$
for $t>0$, we say that $\Lambda$ is an \emph{attracting
  set}. The \emph{topological basin} of an attracting set $\Lambda$
  is the set
  \[
  W^s(\Lambda)=
\{ x\in M : 
\lim_{t\to+\infty}\dist\big( X_t(x) , \Lambda\big) =0 \}.
  \]
  We say that an attracting set $\Lambda$ is
\emph{transitive} if it coincides with the $\omega$-limit
set of a regular $X$-orbit.
\begin{definition}\label{def:attractor}
  An \emph{attractor} is a transitive attracting set, and a
  \emph{repeller} is an attractor for the reversed vector
  field $-X$.
\end{definition}

An attractor, or repeller, is \emph{proper} if it is not the
whole manifold.  An invariant set of $X$ is
\emph{non-trivial} if it is neither a periodic orbit nor a
singularity.

\begin{definition}
\label{d.dominado}
Let $\Lambda$ be a compact invariant  set of $X \in 
{\cal X}^r(M)$ , $c>0$, and $0 < \lambda < 1$.
We say that $\Lambda$ has a $(c,\lambda)$-dominated splitting if
the  bundle over $\Lambda$ can be written as a continuous 
$DX_t$-invariant sum of sub-bundles
$$
T_\Lambda M=E^1\oplus E^2,
$$
such that for every $t > 0$ and every $x \in \Lambda$, we have
\begin{equation}\label{eq.domination}
\|DX_t \mid E^1_x\| \cdot 
\|DX_{-t} \mid E^2_{X_t(x)}\| < c \, \lambda^t.
\end{equation}
\end{definition}

The domination condition \eqref{eq.domination} implies that the
direction of the flow is contained in one of the sub-bundles.

We stress that we only deal with flows in dimension $3$.  In
all that follows, the first sub-bundle $E^1$ will be
one-dimensional, and the flow direction will be contained in
the second sub-bundle $E^2$, that we call \emph{central
  direction} and denote by $E^{cu}$.

We say that a $X$-invariant subset $\Lambda$ of $M$ is
\emph{partially hyperbolic} if it has a
$(c,\lambda)$-dominated splitting, for some $c>0$ and
$\lambda\in(0,1)$, such that the sub-bundle $E^1=E^s$ is
uniformly contracting: for every $t > 0$ and every $x \in
\Lambda$ we have
$$
\|DX_t \mid E^s_x\| < c \, \lambda^t.
$$

For $x\in \Lambda$ and $t\in\real$ we let $J_t^c(x)$ be the absolute 
value of the determinant of the linear map 
$$
DX_t \mid E^{cu}_x:E^{cu}_x\to E^{cu}_{X_t(x)}.
$$
We say that the sub-bundle $E^{cu}_\Lambda$ of the
partially hyperbolic invariant set $\Lambda$ is \emph{volume
  expanding} if $J_t^c(x)\geq c\, e^{-\lambda t}$ for every
$x\in \Lambda$ and $t\geq 0$.  In this case we say that
$E^{cu}_\Lambda$ is {\em $(c,\lambda)$-volume expanding} to
indicate the dependence on $c,\lambda$.

\begin{definition}
\label{d.singularset}
Let $\Lambda$ be a compact invariant set of $X \in {\cal
  X}^r(M)$ with singularities.  We say that $\Lambda$ is a
\emph{singular-hyperbolic set} for $X$ if all the
singularities of $\Lambda$ are hyperbolic, and $\Lambda$ is
partially hyperbolic with volume expanding central
direction.
\end{definition}


\subsection{Expansiveness}
\label{sec:expansiveness}

The flow is \emph{sensitive to initial data} if there is $\delta>0$
such that, for any $x \in M$ and any neighborhood $N$ of $x$,
there is $y \in N$ and $t \geq 0$ such that
$\dist(X_t(x),X_t(y)) > \delta$. 

We shall work with a much stronger property, called
\emph{expansiveness}.  Denote by $S(\re)$ the set of
surjective increasing continuous functions $h:\re\to \re$.
We say that the flow is \emph{expansive} if for every
$\vep>0$ there is $\delta >0$ such that, for any $h \in
S(\re)$, if 
$$
\dist(X_t(x),X_{h(t)}(y)) \leq \delta
\quad\text{for all\ \ } t \in \re, 
$$
then $X_{h(t_0)}(y)
\in X_{[t_0-\vep,t_0+\vep]}(x)$, for some $t_0 \in \re$.  We
say that an invariant compact set $\Lambda$ is expansive if
the restriction of $X_t$ to $\Lambda$ is an expansive flow.

This notion was proposed by Komuro in~\cite{Km84}, and he called it
$K^*$-expansiveness. He proved that a geometric Lorenz
attractor is expansive in this sense. 
Our first main result generalizes this to any singular-hyperbolic
attractor.

\begin{maintheorem}
\label{sci}
Let $\Lambda$ be a singular-hyperbolic attractor of $X\in
{\cal X}^1(M)$.  Then $\Lambda$ is expansive.
\end{maintheorem}

An immediate consequence of this theorem is the following

\begin{maincorollary}
  A singular-hyperbolic attractor of a $3$-flow is sensitive
  to initial data.
\end{maincorollary} 

A stronger notion of expansiveness has been proposed by
Bowen-Walters~\cite{BoWa72}. In it one considers continuous
maps $h:\re\to \re$ with $h(0)=0$, instead.
This turns out to be unsuitable when dealing with singular
sets, because it implies that all singularities are
isolated~\cite[Lemma~1]{BoWa72}.  An intermediate definition
was also proposed by Keynes-Sears~\cite{KS79}: the set of
maps is the same as in \cite{Km84}, but they
require $t_0=0$.  Komuro~\cite{Km84} shows that a
geometric Lorenz attractor does not satisfy this condition.


\subsection{Physical measure}
\label{sec:sinai-ruelle-bowen}
An invariant probability $\mu$ is a \emph{physical} measure
for the flow $X_t$, $t\in\re$ if the set $B(\mu)$ of points
$z\in M$ satisfying
$$
\lim_{T\to+\infty} \frac{1}{T} \int_0^T
\varphi\big( X_t(z) \big) \, dt
= \int\varphi\,d\mu
\quad\text{for all continuous } \varphi: M \to \re
$$
has positive Lebesgue measure: $m\big( B(\mu) \big)>0$.
In that case, $B(\mu)$ is called the \emph{basin} of $\mu$.

Physical measures for singular-hyperbolic attractors were
constructed by Colmen\'arez~\cite{clm2002}.  We need to
assume that $(X_t)_{t\in\real}$ is a flow of class $C^2$
since for the construction of physical measures a bounded
distortion property for one-dimensional maps is needed.
These maps are naturally obtained as quotient maps over the
set of stable leaves, which form a $C^{1+\alpha}$ foliation
of a finite number of cross-sections associated to the flow
if the flow is $C^2$, see
Section~\ref{sec:global-poincare-map}.

\begin{maintheorem}
\label{srb}
Let $\Lambda$ be a singular-hyperbolic
 \emph{attractor}.  Then $\Lambda$ supports a unique physical
  probability measure $\mu$ which is ergodic, hyperbolic
   and its ergodic basin covers a full Lebesgue
 measure subset of the topological basin of attraction, i.e.
 $B(\mu)=W^s(\Lambda),\, m\bmod0$.
\end{maintheorem}

This statement extends the main result in
Colmen\'arez~\cite{clm2002}, where hyperbolicity of the physical
measure was not proved and the author assumed that periodic
orbits in $\Lambda$ exist and are dense. However in another
recent work, Arroyo and Pujals~\cite{AP} show that every
singular-hyperbolic attractor has a dense set of periodic
orbits, so the denseness assumption is no restriction. Here
we give an independent proof of the existence of SRB
measures which does not use denseness of periodic orbits and
that enables us to obtain the hyperbolicity of the SRB
measure.

Here hyperbolicity means \emph{non-uniform hyperbolicity}:
the tangent bundle over $\Lambda$ splits into a sum $T_z M =
E^s_z\oplus E^X_z\oplus F_z$ of three one-dimensional
invariant subspaces defined for $\mu$-a.e.  $z\in \Lambda$
and depending measurably on the base point $z$, where $\mu$
is the physical measure in the statement of
Theorem~\ref{srb}, $E^X_z$ is the flow direction (with zero
Lyapunov exponent) and $F_z$ is the direction with positive
Lyapunov exponent, that is, for every non-zero vector $v\in
F_z$ we have
\[
\lim_{t\to+\infty}\frac1t\log\|DX_t(z)\cdot v\|>0.
\]
We note that the invariance of the splitting implies that
$E^{cu}_z=E^X_z\oplus F_z$ whenever $F_z$ is defined.
For a proof of non-uniform hyperbolicity without using the
existence of invariant measures, but assuming density of
periodic orbits, see Colmen\'arez \cite{clm2005}.
 
Theorem~\ref{srb} is another statement of sensitiveness,
this time applying to the whole open set $B(\Lambda)$.
Indeed, since non-zero Lyapunov exponents express that the
orbits of infinitesimally close-by points tend to move apart
from each other, this theorem means that most orbits in the
basin of attraction separate under forward iteration.  See
Kifer~\cite{Ki88}, and Metzger~\cite{mtz001}, and references
therein, for previous results about invariant measures and
stochastic stability of the geometric Lorenz models.


\subsection{The physical measure is a $u$-Gibbs state}
\label{sec:phys-meas-are}

In the uniformly hyperbolic setting it is well known that
physical measures for hyperbolic attractors admit a
disintegration into conditional measures along the unstable
manifolds of almost every point which are absolutely
continuous with respect to the induced Lebesgue measure on
these sub-manifolds, see \cite{Bo75,BR75,PS82,Vi97b}.

Here the existence of unstable manifolds is guaranteed by
the hyperbolicity of the physical measure: the
strong-unstable manifolds $W^{uu}(z)$ are the ``integral
manifolds'' in the direction of the one-dimensional
sub-bundle $F$, tangent to $F_z$ at almost every
$z\in\Lambda$. The sets $W^{uu}(z)$ are embedded
sub-manifolds in a neighborhood of $z$ which, in general,
depend only measurably (including its size) on the base
point $z\in\Lambda$.
The \emph{strong-unstable manifold} is defined by
\[
W^{uu}(z)=\{ y\in M:
\lim_{ t\to -\infty} \dist(X_t(y),X_t(z))=0 \}
\]
and exists for almost every $z\in\Lambda$ with respect to
the physical and hyperbolic measure obtained in
Theorem~\ref{srb}. We remark that since $\Lambda$ is an
attracting set, then $W^{uu}(z)\subset\Lambda$ whenever
defined.

The tools developed to prove Theorem~\ref{srb} enable us to
prove that the physical measure obtained there has
absolutely continuous disintegration along the
center-unstable direction. To state this result precisely we
need the following notations.

The uniform contraction along the $E^s$ direction ensures
the existence of \emph{strong-stable one-dimensional
  manifolds $W^{ss}(x)$} through every point $x\in\Lambda$,
tangent to $E^s(x)$ at $x$. Using the action of the flow we
define the \emph{stable manifold of $x\in\Lambda$} by
\[
W^{s}(x)=\bigcup_{t\in\real} X_t\big(W^{ss}(x)\big).
\]
Analogously for $\mu$-a.e. $z$ we can define the
\emph{unstable-manifold of $z$} by
\[
W^{u}(z)=\bigcup_{t\in\real} X_t\big(W^{uu}(z)\big).
\]
We note that $E^{cu}_z$ is tangent to $W^{u}(z)$ at $z$ for
$\mu$-a.e. $z$.  Given $x\in\Lambda$ let $S$ be a smooth
surface in $M$ which is everywhere transverse to the vector
field $X$ and $x\in S$, which we call a \emph{cross-section
  of the flow at $x$}. Let $\xi_0$ be the connected
component of $W^s(x)\cap S$ containing $x$. Then $\xi_0$ is
a smooth curve in $S$ and we take a parametrization
$\psi:[-\epsilon,\epsilon]\times[-\epsilon,\epsilon]\to S$
of a compact neighborhood $S_0$ of $x$ in $S$, for some
$\epsilon>0$, such that
\begin{itemize}
\item $\psi(0,0)=x$ and $\psi\big(
  (-\epsilon,\epsilon)\times\{0\} \big) \subset \xi_0$;
\item $\xi_1=\psi\big( \{0\}\times (-\epsilon,\epsilon)
  \big)$ is transverse to $\xi_0$ at $x$:
  $\xi_0\pitchfork\xi_1=\{x\}$.
\end{itemize}

We consider the family $\Pi(S_0)$ of connected components
$\zeta$ of $W^u(z)\cap S_0$ containing $z\in S_0$ which
\emph{cross} $S_0$. We say that a \emph{curve $\zeta$
  crosses $S_0$} if it can be written as the graph of a map
$\xi_1\to\xi_0$.

Given $\delta>0$ we let $\Pi_\delta(x)=\{
X_{(\delta,\delta)} (\zeta) : \zeta\in \Pi(S_0)\}$ be a
family of surfaces inside unstable leaves in a neighborhood
of $x$ crossing $S_0$. The volume form $m$ induces a volume
form $m_\gamma$ on each $\gamma\in\Pi_\delta(x)$
naturally. Moreover, since $\gamma\in\Pi_\delta(x)$ is a
continuous family of curves ($S_0$ is compact and each curve
is tangent to a continuous sub-bundle $E^{cu}$), it forms a
measurable partition of $\hat\Pi_\delta(x)=\cup\{\gamma:
\gamma\in\Pi_\delta(x)\}$. We say that $\Pi_\delta(x)$
is a \emph{$\delta$-adapted foliated neighborhood of $x$}.

Hence $\mu\mid
\hat\Pi_\delta(x)$ can be disintegrated along the
partition $\Pi_\delta(x)$ into a family of measures
$\{\mu_\gamma\}_{\gamma\in\Pi_\delta(x)}$ such that
\[
\mu\mid
\hat\Pi_\delta(x) = \int \mu_\gamma \, d\hat\mu(\gamma),
\]
where $\hat\mu$ is a measure on $\Pi_\delta(x)$ defined by
\[
\hat\mu(A)=\mu\left(\cup_{\gamma\in A} \gamma\right)
\quad\mbox{for all Borel sets   }
A\subset \Pi_\delta(x).
\]
We say that \emph{$\mu$ has an absolutely continuous
  disintegration along the center-unstable direction} if
\emph{for every given $x\in\Lambda$, each $\delta$-adapted
  foliated neighborhood $\Pi_\delta(x)$ of $x$ induces a
  disintegration $\{\mu_\gamma\}_{\gamma\in\Pi_\delta(x)}$
  of $\mu\mid\hat\Pi_\delta(x)$, for all small enough
  $\delta>0$, such that $\mu_\gamma\ll m_\gamma$ for
  $\hat\mu$-a.e. $\gamma\in\Pi_\delta(x)$ }. (See
Section~\ref{sec:absolute-contin-foli} for more details.)

\begin{maintheorem}
\label{thm:srbmesmo}
Let $\Lambda$ be a singular-hyperbolic attractor for a $C^2$
three-dimen\-sional flow.  Then the physical measure $\mu$
supported in $\Lambda$ has a disintegration into absolutely
continuous conditional measures $\mu_\gamma$ along
center-unstable surfaces $\gamma\in\Pi_\delta(x)$ such that
$\frac{d\mu_\gamma}{dm_\gamma}$ is uniformly bounded from
above, for all $\delta$-adapted foliated neighborhoods
$\Pi_\delta(x)$ and every $\delta>0$.  Moreover
$\suporte(\mu)=\Lambda\,$.
\end{maintheorem}

\begin{remark}
  \label{rmk:suppLambda}
  The proof that $\supp(\mu)=\Lambda$ presented here 
  depends on the abosultely continuous disintegration
  property of $\mu$. 
\end{remark}

\begin{remark}
  \label{rmk:boundedensity}
  It follows from our arguments that the densities of the
  conditional measures $\mu_\gamma$ are bounded from below
  away from zero on $\Lambda\setminus B$, where $B$ is any
  neighborhood of the singularities $\si(X\mid\Lambda)$. In
  particular the densities tend to zero as we get closer to
  the singularities of $\Lambda$.
\end{remark}


The absolute continuity
property along the center-unstable sub-bundle given by
Theorem~\ref{thm:srbmesmo} ensures that
\[
h_\mu(X_1)=\int \log\big| \det (DX_1\mid E^{cu})  \big| \, d\mu,
\]
by the characterization of probability measures satisfying
the Entropy Formula~\cite{LY85}.  The above integral is the
sum of the positive Lyapunov exponents along the sub-bundle
$E^{cu}$ by Oseledets Theorem~\cite{Man87,Wa82}.  Since in the
direction $E^{cu}$ there is only one positive Lyapunov
exponent along the one-dimensional direction $F_z$,
$\mu$-a.e. $z$, the ergodicity of $\mu$ then shows that the
following is true.

\begin{maincorollary}
  \label{cor:uGibbs}
  If $\Lambda$ is a singular-hyperbolic attractor for a
  $C^2$ three-dimen\-sional flow $X_t$, then the physical
  measure $\mu$ supported in $\Lambda$ satisfies the Entropy
  Formula
\[
h_\mu(X_1)=\int\log\| DX_1\mid F_z\|\, d\mu(z).
\]
\end{maincorollary}

Again by the characterization of measures satisfying the
Entropy Formula we get that \emph{$\mu$ has absolutely
  continuous disintegration along the strong-unstable
  direction}, along which the Lyapunov exponent is positive,
thus \emph{$\mu$ is a $u$-Gibbs state}~\cite{PS82}.  This
also shows that \emph{$\mu$ is an equilibrium state for the
  potential} $-\log\| DX_1\mid F_z\|$ with respect to the
diffeomorphism $X_1$. We note that the entropy $h_\mu(X_1)$
of $X_1$ is the entropy of the flow $X_t$ with respect to
the measure $\mu$ \cite{Wa82}.

Hence we are able to extend most of the basic results on the
ergodic theory of hyperbolic attractors to the setting of
singular-hyperbolic attractors.


\subsection{Application to the Lorenz and geometric Lorenz
  flows}
\label{sec:Appl-to-Lorenz}

It is well known that geometric Lorenz flows are transitive
and it was proved in~\cite{MPP99} that they are
singular-hyperbolic attractors. Then as a consequence of our
results we get the following corollary.

\begin{maincorollary}
  A geometric Lorenz flow is expansive and has a unique
  physical invariant probability measure whose basin covers
  Lebesgue almost every point of the topological basin of
  attraction. Moreover this measure is a $u$-Gibbs state and
  satisfies the Entropy Formula.
\end{maincorollary} 

Recently Tucker~\cite{Tu2} proved that the flow defined by
the Lorenz equations~\cite{Lo63} exhibits a
singular-hyperbolic attractor. In particular our results
then show the following.

\begin{maincorollary}
  The flow defined by the Lorenz equations is expansive and
  has a unique physical invariant probability measure whose
  basin covers Lebesgue almost every point of the
  topological basin of attraction. Moreover this measure is
  a $u$-Gibbs state and satisfies the Entropy Formula.
\end{maincorollary}

This paper is organized as follows. In
Section~\ref{sec:cross-sect-poinc} we obtain adapted
cross-sections for the flow near $\Lambda$ and deduce some
hyperbolic properties for the Poincar\'e return maps between
these sections to be used in the sequel. Theorem~\ref{sci}
is proved in Section~\ref{sec:proof-expansiveness}.  In
Section~\ref{sec:proof-theorem-B} we outline the proof of
Theorem~\ref{srb}, which is divided into several steps
detailed in Sections~\ref{sec:global-poincare-map} through
\ref{sec:phys-meas-flow}.  In
Section~\ref{sec:global-poincare-map} we reduce the dynamics
of the global Poincar\'e return map between cross-sections
to a one-dimensional piecewise expanding map. In
Sections~\ref{s.fsuspending} and~\ref{sec:phys-meas-flow} we
explain how to construct invariant measures for the
Poincar\'e return map from invariant measures for the
induced one-dimensional map, and also how to obtain
invariant measures for the flow through invariant measures
for the Poincar\'e return map. This concludes the proof of
Theorem~\ref{srb}.

Finally, in Section~\ref{sec:phys-meas-are-1} we again use
the one-dimensional dynamics and the notion of hyperbolic
times for the Poincar\'e return map to prove that the
physical measure is SRB 
and that $\suporte(\mu)=\Lambda$, concluding the proof of
Theorem~\ref{thm:srbmesmo} and of
Corollary~\ref{cor:uGibbs}.
 
\subsection*{Acknowledgments}
We are grateful to the referee for the careful revison of
the paper and the many valuable suggestions which greatly
improved the readibility of the text.

\section{Cross-sections and Poincar\'e maps} 
\label{sec:cross-sect-poinc}

The proof of Theorem~\ref{sci} is based on analyzing
Poincar\'e return maps of the flow to a convenient
cross-section.  In this section we give a few properties of
\emph{Poincar\'e maps}, that is, continuous maps
$R:\Sigma\to\Sigma'$ of the form $R(x)=X_{t(x)}(x)$ between
cross-sections $\Sigma$ and $\Sigma'$.  We always assume
that the Poincar\'e time $t(\cdot)$ is large
(Section~\ref{s.22}). Recall that we assume
singular-hyperbolicity.

Firstly, we observe (Section~\ref{s.21}) that cross-sections
have co-dimension $1$ foliations which are dynamically
defined: the leaves $W^s(x,\Sigma)=W^s_{loc}(x)\cap\Sigma$
correspond to the intersections with the stable manifolds of
the flow.  These leaves are uniformly contracted
(Section~\ref{s.22}) and, assuming the cross-section is
\emph{adapted} (Section~\ref{s.23}) the foliation is
invariant:
$$
R(W^s(x,\Sigma))\subset W^s(R(x),\Sigma')
\quad\text{for all } x \in \Lambda\cap\Sigma.
$$
Moreover, $R$ is uniformly expanding in the transverse
direction (Section~\ref{s.22}). In Section~\ref{s.24} we
analyze the flow close to singularities, again by means of
cross-sections.

\subsection{Stable foliations on cross-sections}
\label{s.21}

We begin by recalling a few classical facts about partially
hyperbolic systems, especially existence of strong-stable
and center-unstable foliations.  The standard reference is
\cite{HPS77}.

Hereafter, $\Lambda$ is a singular-hyperbolic attractor of
$X\in {\cal X}^1(M)$ with invariant splitting $T_\Lambda M =
E^{s}\oplus E^{cu}$ with $\dim E^{cu}=2$.  Let
$\tilde{E}^s\oplus \tilde{E}^{cu}$ be a continuous extension
of this splitting to a small neighborhood $U_0$ of
$\Lambda$.  For convenience, we take $U_0$ to be forward
invariant.  Then $\tilde{E}^s$ may chosen invariant under
the derivative: just consider at each point the direction
formed by those vectors which are strongly contracted by
$DX_t$ for positive $t$.  In general, $\tilde{E}^{cu}$ is
not invariant. However, we can always consider a cone field
around it on $U_0$
$$
C^{cu}_a(x)=\{v=v^s+v^{cu}: v^s\in \tilde{E}^s_x
\text{ and }
v^u\in\tilde{E}^{cu}_x
\text{ with } \|v^s\|\le a\cdot \|v^{cu}\|\}
$$
which is forward invariant for $a>0$:
\begin{equation}
\label{eq.cone3}
DX_t(C^{cu}_a(x)) \subset C^{cu}_a (X_t(x))
\quad\text{for all large $t>0$.}
\end{equation}
Moreover, we may take $a>0$ arbitrarily small, reducing
$U_0$ if necessary.  For notational simplicity, we write
$E^s$ and $E^{cu}$ for $\tilde E^s$ and $\tilde E^{cu}$ in
all that follows.

The next result asserts that there exist locally strong-stable
and center-unstable manifolds, defined at every regular point
$x\in U_0$\,, which are embedded disks tangent to $E^s(x)$ and
$E^{cu}(x)$, respectively. The strong-stable manifolds are
locally invariant. Given any $x\in U_0$\,, define
$$
W^{ss}(x)=
\{y\in M: \dist(X_t(x),X_t(y))\to 0 \text{ as } t \to +\infty\}
$$
$$
W^s(x) = \bigcup_{t\in\real} W^{ss}(X_t(x))
       = \bigcup_{t\in\real} X_t(W^{ss}(x)).
$$
Given $\vep > 0$, denote $I_\vep=(-\vep, \vep)$ and let
${\cal E}^1(I_1,M)$ be the set of $C^1$ embedding maps 
$f:I_1\to M$ endowed with the $C^1$ topology.

\begin{proposition}(stable and center-unstable manifolds)
\label{p.centroinstavel}
There are continuous maps
$\phi^{ss}:U_0\to {\cal E}^1(I_1,M)$ and 
$\phi^{cu}:U_0\to {\cal E}^1(I_1\times I_1, M)$ such that
given any $0 <\vep< 1$ and $x \in U_0$, if we denote 
$W_\vep^{ss}(x)=\phi^{ss}(x)(I_\vep)$
and $W^{cu}_\vep(x)=\phi^{cu}(x)(I_\vep\times I_\vep)$,

\begin{itemize}
\item[(a)] 
$T_x W^{ss}_\vep(x)=E^s(x)$;

\item [(b)] 
$T_xW^{cu}_\vep(x)=E^{cu}(x)$;

\item[(c)]
$W_\vep^{ss}(x)$ is a neighborhood of $x$ inside $W^{ss}(x)$;

\item[(d)]
$y\in W^{ss}(x) \Leftrightarrow$ there is $T\ge 0$ such
that $X_T(y) \in W^{ss}_\vep(X_T(x))$ (local invariance);

\item [(e)]
$d(X_t(x),X_t(y))\le c \cdot \lambda^t \cdot d(x,y)$ for all $t>0$
and all $y\in W_\vep^{ss}(x)$.

\end{itemize}
\end{proposition}

The constants $c>0$ and $\lambda\in(0,1)$ are taken as in
Definition~\ref{d.dominado} and the distance $d(x,y)$ is the
intrinsic distance between two points on the manifold
$W_\vep^{ss}(x)$, given by the length of the shortest smooth
curve contained in $W_\vep^{ss}(x)$ connecting $x$ to $y$.  

Denoting $E^{cs}_x=E^s_x\oplus E^{X}_x$, where $E^X_x$ is
the direction of the flow at $x$, it follows that
\begin{align}\label{eq:Ecs}
  T_x W^{ss}(x)=E^s_x \quad\text{and}\quad T_x
  W^{s}(x)=E^{cs}_x\,.
\end{align}
We fix $\vep$ once and for all. Then we call
$W^{ss}_{\vep}(x)$ the local \emph{strong-stable manifold}
and $W^{cu}_{\vep}(x)$ the local \emph{center-unstable
manifold} of $x$.


Now let $\Sigma$ be a \emph{cross-section} to the flow, that is,
a $C^2$ embedded compact disk transverse to $X$ at every point.
For every $x\in\Sigma$ we define $W^s(x,\Sigma)$ to be the
connected component of $W^s(x)\cap\Sigma$ that contains $x$.
This defines a foliation $\F^{s}_{\Sigma}$ of
$\Sigma$ into co-dimension $1$ sub-manifolds of class $C^1$. 

\begin{figure}[ht]
\begin{center}
  \includegraphics[width=12cm,height=5.5cm]{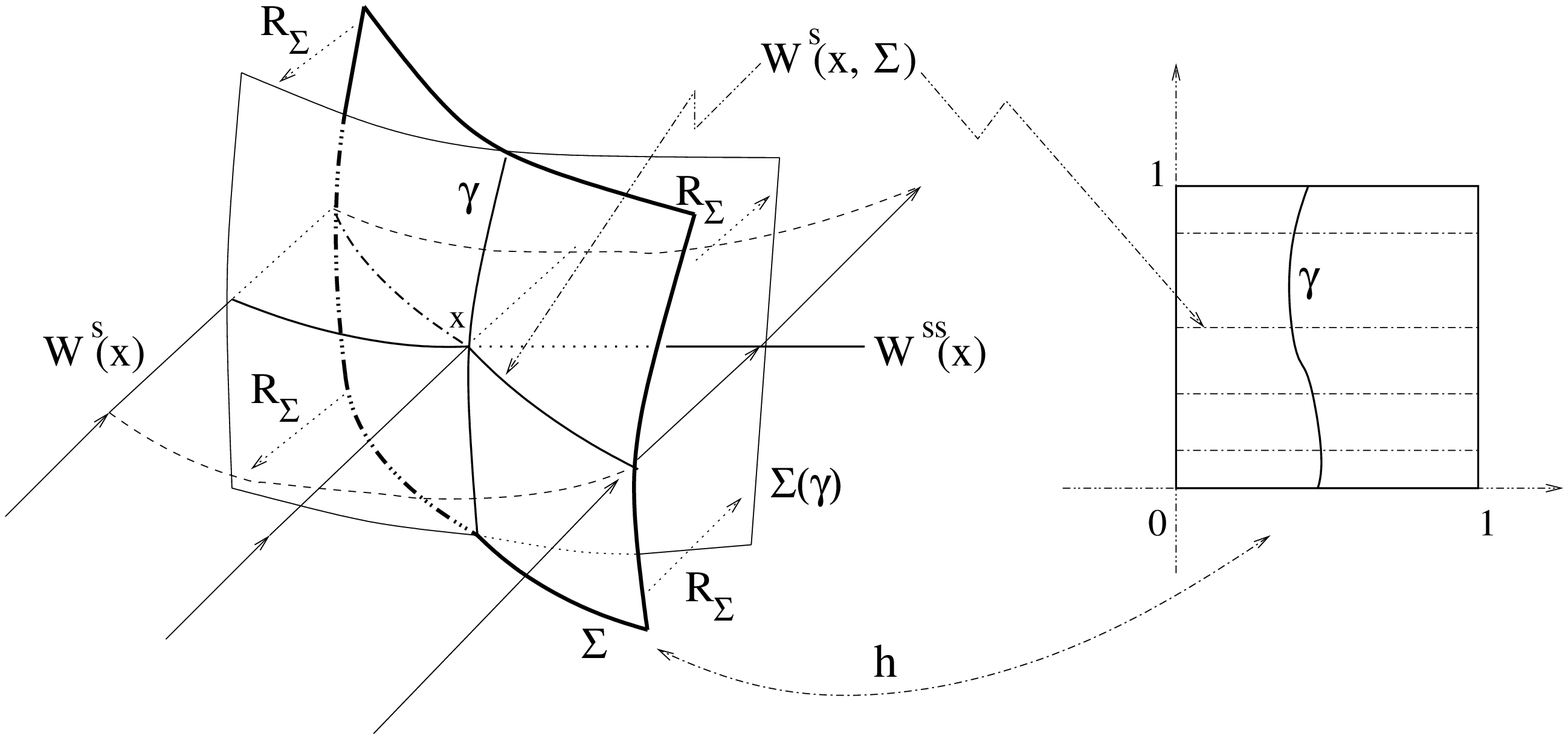}
  \caption{\label{f.squaresection} The sections $\Sigma$,
    $\Sigma(\gamma)$, the manifolds $W^s(x), W^{ss}(x)$,
    $W^s(x,\Sigma)$ and the projection $R_\Sigma$, on the
    right. On the left, the square $[0,1]\times[0,1]$ is
    identified with $\Sigma$ through the map $h$, where
    $\F_\Sigma^s$ becomes the horizontal foliation and the
    curve $\gamma$ is transversal to the horizontal
    direction. Solid lines with arrows indicate the flow
    direction.}
\end{center}
\end{figure}

\begin{remark}\label{r.foliated}
  Given any cross-section $\Sigma$ and a point $x$ in its
  interior, we may always find a smaller cross-section also
  with $x$ in its interior and which is the image of the
  square $[0,1]\times[0,1]$ by a $C^2$ diffeomorphism $h$
  that sends horizontal lines inside leaves of
  $\F^{s}_{\Sigma}$. So, in what follows we always assume
  cross-sections are of the latter kind, see
  Figure~\ref{f.squaresection}. We denote by
  $\interior(\Sigma)$ the image of $(0,1)\times(0,1)$ under
  the above-mentioned diffeomorphism, which we call the
  \emph{interior} of $\Sigma$.

  We also assume that each cross-section $\Sigma$ is
  contained in $U_0$, so that every $x\in\Sigma$ is such
  that $\omega(x)\subset \Lambda$.
\end{remark}

\begin{remark}\label{r.trapaca}
  In general, we can not choose the cross-section such that
  $W^s(x,\Sigma)\subset W^{ss}_\vep(x)$. The reason is that
  we want cross-sections to be $C^2$. Cross-section of class
  $C^1$ are enough for the proof of expansiveness in
  Section~\ref{sec:proof-expansiveness} but $C^2$ is needed
  for the construction of the physical measure in
  Sections~\ref{sec:global-poincare-map} through
  \ref{sec:global-poincare-map} and for the absolute
  continuity results in Section~\ref{sec:phys-meas-are-1}.
  The technical reason for this is explained in
  Section~\ref{sec:absolute-contin-foli}.

  On the one hand $x\mapsto W^{ss}_\vep(x)$ is usually not
  differentiable if we assume that $X$ is only of class
  $C^1$.  On the other hand, assuming that the cross-section
  is small with respect to $\vep$, and choosing any curve
  $\gamma\subset\Sigma$ crossing transversely every leaf of
  $\F_\Sigma^s$\,, we may consider a Poincar\'e map
  $$
  R_\Sigma:\Sigma \to \Sigma(\gamma)=
  \bigcup_{z\in\gamma} W^{ss}_\vep(z)
  $$
  with Poincar\'e time close to zero, see
  Figure~\ref{f.squaresection}.  This is a homeomorphism
  onto its image, close to the identity, such that
  $R_\Sigma(W^s(x,\Sigma))\subset W^{ss}_\vep(R_\Sigma(x))$.
  So, identifying the points of $\Sigma$ with their images
  under this homeomorphism, we may pretend that indeed
  $W^s(x,\Sigma)\subset W^{ss}_\vep(x)$.  We shall often do
  this in the sequel, to avoid cumbersome technicalities.
\end{remark}

\subsection{Hyperbolicity of Poincar\'e maps}
\label{s.22}

Let $\Sigma$ be a small cross-section to $X$ and let
$R:\Sigma\to\Sigma'$ be a Poincar\'e map $R(y)=X_{t(y)}(y)$ to
another cross-section $\Sigma'$ (possibly $\Sigma=\Sigma'$).
Note that $R$ needs not correspond to the first time the orbits
of $\Sigma$ encounter $\Sigma'$\,, nor it is defined
everywhere in $\Sigma$.

The splitting $E^s\oplus E^{cu}$ over $U_0$ induces a continuous 
splitting $E_\Sigma^s\oplus E_\Sigma^{cu}$ of the tangent bundle
$T\Sigma$ to $\Sigma$ (and analogously for $\Sigma'$),
defined by (recall~\eqref{eq:Ecs} for the use of $E^{cs}$)
\begin{equation}\label{eq.splitting}
E_\Sigma^s(y)=E^{cs}_y\cap T_y{\Sigma}
\quad\mbox{and}\quad
E_\Sigma^{cu}(y)=E^{cu}_y\cap T_y{\Sigma}.
\end{equation}
We are going to prove that if the Poincar\'e time $t(x)$ is
sufficiently large then \eqref{eq.splitting} defines a hyperbolic
splitting for the transformation $R$ on the cross-sections, 
at least restricted to $\Lambda$:

\begin{proposition}\label{p.secaohiperbolica}
Let $R:\Sigma\to\Sigma'$ be a Poincar\'e map as before 
with Poincar\'e time $t(\cdot)$.
Then $DR_x(E_\Sigma^s(x)) = E_\Sigma^s(R(x))$ at every
$x\in\Sigma$ and $DR_x(E_\Sigma^{cu}(x)) = E_\Sigma^{cu}(R(x))$
at every $x\in\Lambda\cap\Sigma$.

Moreover for every given $0<\lambda<1$ there exists
$t_1=t_1(\Sigma,\Sigma',\lambda)>0$ such that if
$t(\cdot)>t_1$ at every point, then
$$
\|DR \mid E^s_\Sigma(x)\| < \lambda
\quad\text{and}\quad
\|DR \mid E^{cu}_\Sigma(x)\| > 1/\lambda
\quad\text{at every $x\in\Sigma$.}
$$
\end{proposition}

\begin{remark}
\label{r.angulos1}
In what follows we use $K$ as a generic notation for large
constants depending only on a lower bound for the angles
between the cross-sections and the flow direction, and on
upper and lower bounds for the norm of the vector field on
the cross-sections. The conditions on $t_1$ in the proof of
the proposition depend only on these bounds as well.  In all
our applications, all these angles and norms will be
uniformly bounded from zero and infinity, and so both $K$
and $t_1$ may be chosen uniformly.
\end{remark}

\begin{proof}
The differential of the Poincar\'e map at any point $x\in\Sigma$ 
is given by 
$$
DR(x)= P_{R(x)} \circ DX_{t(x)}\mid T_x\Sigma,
$$
where $P_{R(x)}$ is the projection onto $T_{R(x)}\Sigma'$ along
the direction of $X(R(x))$\,.
Note that $E^s_\Sigma(x)$ is tangent to $\Sigma\cap
W^s(x)\supset W^s(x,\Sigma)$. 
Since the stable manifold $W^s(x)$ is invariant, we have
 invariance of the stable bundle:
$DR(x)\big(E^s_\Sigma(x)\big) = E^s_{\Sigma'}\big(R(x)\big)$.
Moreover for all $x\in\Lambda$ we have
$$
DX_{t(x)} \big(E^{cu}_\Sigma(x)\big)
\subset DX_{t(x)} \big(E_x^{cu}\big)
= E^{cu}_{R(x)}\,.
$$
Since $P_{R(x)}$ is the projection along the vector field, it sends
$E^{cu}_{R(x)}$ to $E^{cu}_{\Sigma'}(R(x))$. This proves that the
center-unstable bundle is invariant restricted to $\Lambda$, i.e.
$DR(x)\big(E^{cu}_\Sigma(x)\big) = E^{cu}_{\Sigma'}(R(x))$.

Next we prove the expansion and contraction statements.  We
start by noting that $\|P_{R(x)}\|\le K$.  Then we consider the
basis $\{\frac{X(x)}{\|X(x)\|},\, e^u_{x}\}$ of
$E^{cu}_{x}$, where $e^u_{x}$ is a unit vector in the
direction of $E^{cu}_\Sigma(x)$.  Since the flow direction is
invariant, the matrix of $DX_{t} \mid E^{cu}_{x}$
relative to this basis is upper triangular:
$$
DX_{t(x)}\mid E^{cu}_{x} =
\left[\begin{array}{cc}
    \frac{\|X(R(x))\|}{\|X(x)\|} & \star  \\
     0  &  \Delta
     \end{array} \right].
$$ 
Moreover
$$
\frac1K\cdot \deter\big(DX_{t(x)}\mid E^{cu}_{x} \big) \le 
\frac{\|X(R(x))\|}{\|X(x)\|} \Delta \le
K \cdot \deter\big(DX_{t(x)}\mid E^{cu}_{x}\big).
$$
Then
\begin{equation*}
\begin{aligned}
\|DR(x)\, e^u_x\| 
& = \|P_{R(x)}\big( DX_{t(x)}(x) \cdot e^u_{x}\big) \|
= \|\Delta \cdot e^u_{R(x)}\| = |\Delta|
\\
&
\ge K^{-3} \, |\deter(DX_{t(x)} \mid  E^{cu}_{x}) |
\ge K^{-3} \lambda^{-t(x)} 
\ge K^{-3} \, \lambda^{-t_1}.
\end{aligned}
\end{equation*}
Taking $t_1$ large enough we ensure that the latter expression is
larger than $1/\lambda$.

To prove $\|DR\mid E^s_\Sigma(x)\| < \lambda$, let us
consider unit vectors $e^s_x \in E^s_x$ and $\hat e^s_x \in
E^s_\Sigma(x)$, and write
$$
e^s_x= a_x \cdot \hat e^s_x + b_x \cdot \frac{X(x)}{\|X(x)\|} \,.
$$
Since $\angle(E^s_x,X(x)) \ge \angle(E_x^s,E^{cu}_x)$ and the latter
is uniformly bounded from zero, we have $|a_x|\ge\kappa$ for 
some $\kappa>0$ which depends only on the flow.
Then
\begin{equation}
\label{eq.contraction}
\begin{aligned}
\|DR(x) \, e_x^s \|
& =\|P_{R(x)}\circ \big(DX_{t(x)}(x) \cdot e^s_x\big)\| 
\\
& = \frac{1}{|a_x|}\,
\left\|P_{R(x)}\circ \Big(DX_{t(x)}(x)
\big(e^s_x - b_x\frac{X(x)}{\|X(x)\|}\big)\Big)\right\|
\\
& = \frac{1}{|a_x|}\,
\left\|P_{R(x)}\circ \big(DX_{t(x)}(x)\cdot \hat e^s_x\big)\right\|
\leq \frac{K}{\kappa} \lambda^{t(x)}
\le \frac{K}{\kappa} \lambda^{t_1}.
\end{aligned}
\end{equation}
Once more it suffices to take $t_1$ large to ensure that the 
right hand side is less than $\lambda$.
\end{proof}

Given a cross-section $\Sigma$, a positive number $\rho$,
and a point $x\in \Sigma$,
we define the unstable cone of width $\rho$ at $x$ by
\begin{equation}
\label{cone}
C_\rho^u(x)=\{v=v^s+v^u : v^s\in E^s_\Sigma(x),\, v^u\in E^{cu}_\Sigma(x)
\mbox{ and } \|v^s\| \le \rho \|v^u\| \}
\end{equation}
(we omit the dependence on the cross-section in our notations). 

Let $\rho>0$ be any small constant. In the following
consequence of Proposition~\ref{p.secaohiperbolica} we
assume the neighborhood $U_0$ has been chose sufficiently
small, depending on $\rho$ and on a bound on the angles
between the flow and the cross-sections.

\begin{corollary}
\label{ccone}
For any $R:\Sigma\to\Sigma'$ as in
Proposition~\ref{p.secaohiperbolica}, with $t(\cdot)>t_1$\,,
and any $x \in\Sigma$, we have
$$
DR(x) (C^u_{\rho}(x)) \subset C_{\rho/2}^u(R(x))
\quad\mbox{and}\quad
\| DR_x(v)\| \ge \frac5{6}\lambda^{-1} \cdot \|v\| 
\quad\mbox{for all}\quad v\in C^u_{\rho}(x).
$$
\end{corollary}

\begin{proof}
  Proposition~\ref{p.secaohiperbolica} immediately implies
  that $DR_x (C^u_{\rho}(x))$ is contained in the cone of
  width $\rho/4$ around $DR(x) \big(E^{cu}_\Sigma(x)\big)$
  relative to the splitting
$$
T_{R(x)}\Sigma^\prime=E^s_{\Sigma^\prime}(R(x))\oplus
DR(x) \big(E^{cu}_\Sigma(x)\big).
$$
(We recall that $E^s_\Sigma$ is always mapped to
$E^s_{\Sigma^\prime}$.)  The same is true for
$E^{cu}_\Sigma$ and $E^{cu}_{\Sigma^\prime}$\,, restricted
to $\Lambda$. So the previous observation already gives the
conclusion of the first part of the corollary in the special
case of points in the attractor.  Moreover to prove the
general case we only have to show that $DR(x)
\big(E^{cu}_\Sigma(x)\big)$ belongs to a cone of width less
than $\rho/4$ around $E^{cu}_{\Sigma'}(R(x))$. This is
easily done with the aid of the flow invariant cone field
$C_a^{cu}$ in \eqref{eq.cone3}, as follows.  On the one
hand,
$$
DX_{t(x)} \big(E^{cu}_\Sigma(x)\big)
 \subset DX_{t(x)} \big(E^{cu}_{x}\big)
 \subset DX_{t(x)} \big(C^{cu}_a(x)\big)
 \subset C^{cu}_a(R(x))\,.
 $$
 We note that $DR(x)\big( E^{cu}_\Sigma (x)\big) =
 P_{R(x)} \circ DX_{t(x)} \big(E^{cu}_\Sigma (x)\big)$.
 Since $P_{R(x)}$ maps $E^{cu}_{R(x)}$ to
 $E^{cu}_{\Sigma'}(R(x))$ and the norms of both $P_{R(x)}$
 and its inverse are bounded by some constant $K$ (see
 Remark~\ref{r.angulos1}), we conclude that
 $DR(x)\big(E^{cu}_\Sigma(x)\big)$ is contained in a cone of
 width $b$ around $E^{cu}_{\Sigma'}(R(x))$, where
 $b=b(a,K)$ can be made arbitrarily small by
 reducing $a$.  We keep $K$ bounded, by assuming the
 angles between the cross-sections and the flow are bounded
 from zero and then, reducing $U_0$ if necessary, we can
 make $a$ small so that $b<\rho/4$.  This concludes
 the proof since the expansion estimate is a trivial
 consequence of Proposition~\ref{p.secaohiperbolica}.
\end{proof}

By a \emph{curve} we always mean the image of a compact interval $[a,b]$
by a $C^1$ map. We use $\ell(\gamma)$ to denote its length.
By a \emph{cu-curve} in $\Sigma$ we mean a curve contained in
the cross-section $\Sigma$ and whose tangent direction
$T_z\gamma\subset C^u_\rho(z)$ for all $z\in\gamma$.
The next lemma says that \emph{cu-curves linking the stable leaves
of nearby points must be short}.

\begin{lemma}
\label{l.lengthversusdistance}
Let us assume that $\rho$ has been fixed, sufficiently small.
Then there exists a constant $\kappa$ such that, for any pair
of points $x, y \in \Sigma$, and any cu-curve $\gamma$ joining $x$ to
some point of $W^s(y,\Sigma)$, we have $\ell(\gamma)\le\kappa
\cdot d(x,y)$.
\end{lemma}

Here $d$ is the intrinsic distance in the $C^2$ surface $\Sigma$.

\begin{proof}
We consider coordinates on $\Sigma$ for which $x$ corresponds
to the origin, $E^{cu}_\Sigma(x)$ corresponds to the
vertical axis, and $E^{s}_\Sigma(x)$ corresponds to the
horizontal axis; through these coordinates we identify
$\Sigma$ with a subset of its tangent space at $x$,
endowed with the Euclidean metric.
In general this identification is not an isometry, but
the distortion is uniformly bounded, and that is taken
care of by the constants $C_1$ and $C_2$ in what follows.
\begin{figure}[ht]
\begin{center}
\includegraphics[height=1.5in]{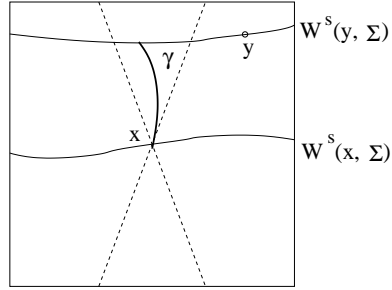}
\caption{\label{f.versus} The stable manifolds on the
  cross-section and the $cu$-curve $\gamma$ connecting
  them.}
\end{center}
\end{figure}
The hypothesis that $\gamma$ is a cu-curve implies that
it is contained in the cone of width $C_1\cdot\rho$ centered
at $x$. On the other hand, stable leaves are close to
being horizontal. It follows (see Figure~\ref{f.versus})
that the length of $\gamma$ is bounded by $C_2\cdot d(x,y)$.
This proves the lemma with $\kappa=C_2$\,.
\end{proof}

In what follows we take $t_1$ in
Proposition~\ref{p.secaohiperbolica} for $\lambda=1/3$. From
Section~\ref{sec:global-poincare-map} onwards we will need
to decrease $\lambda$ once taking a bigger $t_1$.

\subsection{Adapted cross-sections}
\label{s.23}

\begin{figure}[htpb]
  \centering
    \includegraphics[height=4cm]{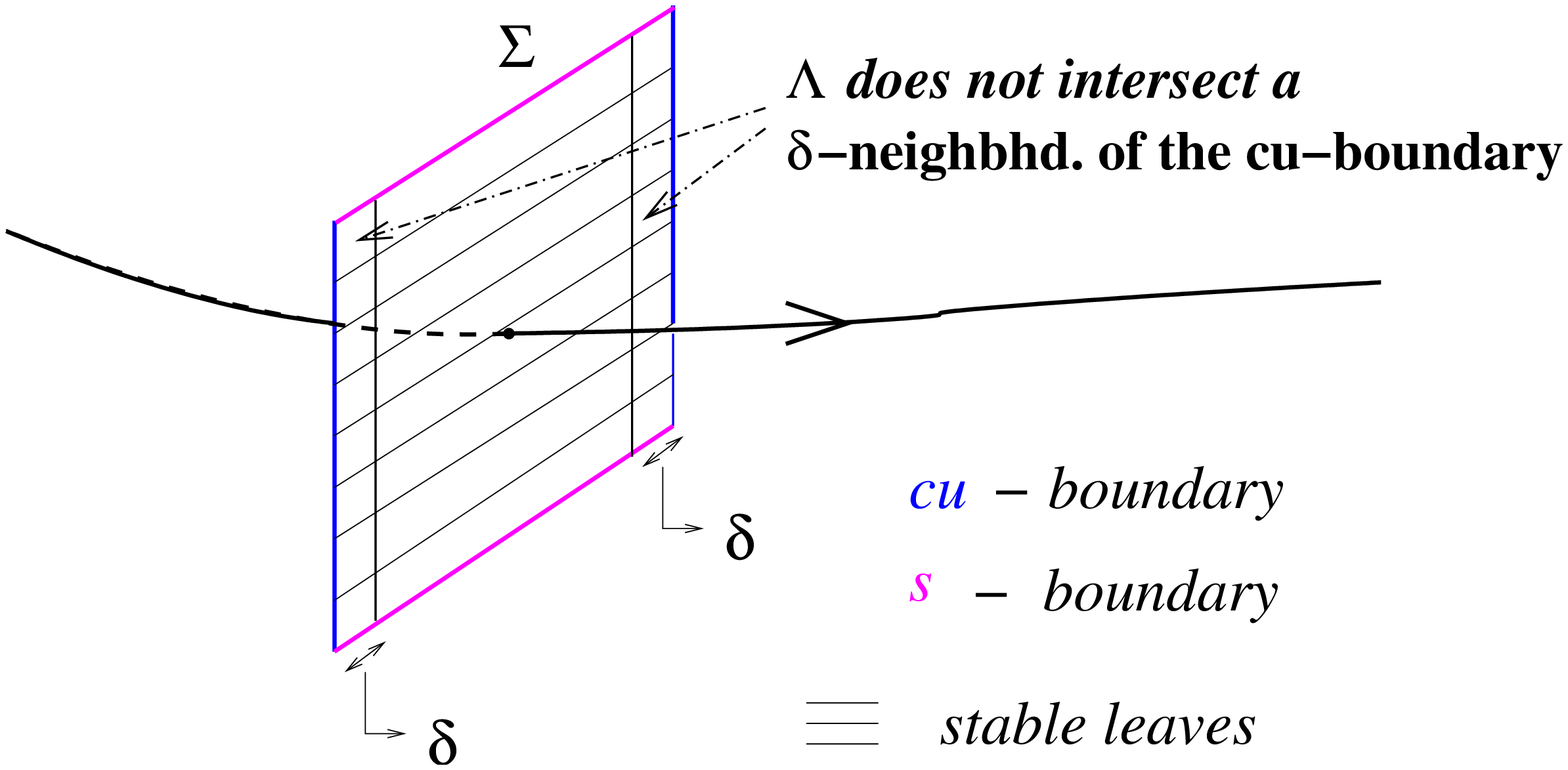}
  \caption{An adapted cross-section for $\Lambda$.}
  \label{fig:an-adapted-cross}
\end{figure}

The next step is to exhibit stable manifolds for Poincar\'e
transformations $R:\Sigma\to\Sigma'$. The natural candidates
are the intersections $W^s(x,\Sigma)=W^s_\vep(x)\cap\Sigma$
we introduced previously. 
These intersections are tangent to the corresponding sub-bundle 
$E^s_\Sigma$ and so, by Proposition~\ref{p.secaohiperbolica},
they are contracted by the transformation.
For our purposes it is also important that the stable foliation
be invariant:
\begin{equation}\label{eq.stableMarkov}
R(W^s(x,\Sigma)) \subset W^s(R(x),\Sigma')
\qquad \text{for every } x\in\Lambda\cap\Sigma.
\end{equation}
In order to have this we restrict somewhat our class of
cross-sections whose center-unstable boundary is disjoint
from $\Lambda$.  Recall (Remark~\ref{r.foliated}) that we
are considering cross-sections $\Sigma$ that are
diffeomorphic to the square $[0,1]\times[0,1]$, with the
horizontal lines $[0,1]\times\{\eta\}$ being mapped to
stable sets $W^s(y,\Sigma)$.  The \emph{stable boundary}
$\partial^{s}\Sigma$ is the image of $[0,1]\times\{0,1\}$.
The \emph{center-unstable boundary} $\partial^{cu}\Sigma$ is
the image of $\{0,1\}\times [0,1]$. The cross-section is
\emph{$\delta$-adapted} if
$$
d(\Lambda \cap \Sigma,\partial^{cu}\Sigma)> \delta,
$$
where $d$ is the intrinsic distance in $\Sigma$, see
Figure~\ref{fig:an-adapted-cross}.  We call \emph{horizontal
  strip} of $\Sigma$ the image $h([0,1]\times I)$ for any
compact subinterval $I$, where $h:[0,1]\times[0,1]\to
\Sigma$ is the coordinate system on $\Sigma$ as in
Remark~\ref{r.foliated}.  Notice that every horizontal strip
is a $\delta$-adapted cross-section.

In order to prove that adapted cross-sections do exist,
we need the following result.

\begin{lemma}\label{l.isolated}
If $\Lambda$ is a singular-hyperbolic attractor, then every point
$x\in\Lambda$ is in the closure of $W^{ss}(x)\setminus\Lambda$.
\end{lemma}

\begin{proof}
The proof is by contradiction.
Let us suppose that there exists $x\in\Lambda$ such that $x$ is in the interior
of $W^{ss}(x)\cap\Lambda$. Let $\alpha(x)\subset\Lambda$ be its
$\alpha$-limit set. Then
\begin{equation}\label{eq.ssinside}
W^{ss}(z)\subset\Lambda \quad\text{for every } z\in \alpha(x),
\end{equation}
since any compact part of the strong-stable manifold of $z$
is accumulated by backward iterates of any small neighborhood
of $x$ inside $W^{ss}(x)$.
It follows that $\alpha(x)$ does not contain any singularity:
indeed, \cite[Theorem~B]{MPP04} proves that the strong-stable 
manifold of each singularity meets $\Lambda$ only at the 
singularity. Therefore by \cite[Proposition~1.8]{MPP04} the
invariant set $\alpha(x)\subset\Lambda$ is hyperbolic.
It also follows from \eqref{eq.ssinside} that the union 
$$
S=\bigcup_{y\in\alpha(x)\cap\Lambda} W^{ss}(y)
$$
of the strong-stable manifolds through the points of
$\alpha(x)$ is contained in $\Lambda$. By continuity of the
strong-stable manifolds and the fact that $\alpha(x)$ is a
closed set, we get that $S$ is also closed. Using 
\cite{MPP04} once more, we see that $S$ does not contain
singularities and, thus, is also a hyperbolic set.

We claim that $W^{u}(S)$, the union of the unstable
manifolds of the points of $S$,  is an open set. To prove this,
we note that $S$ contains the whole stable manifold
$W^s(z)$ of every $z\in S$: this is because $S$ is invariant
and contains the strong-stable manifold of $z$. Now, the
union of the strong-unstable manifolds through the points 
of $W^s(z)$ contains a neighborhood of $z$. This proves 
that $W^u(S)$ is a neighborhood of $S$. Thus the backward
orbit of any point in $W^u(S)$ must enter the interior of
$W^u(S)$. Since the interior is, clearly, an invariant set,
this proves that $W^u(S)$ is open, as claimed. 

Finally,
consider any backward dense orbit in $\Lambda$ (we recall
that for us an attractor is transitive by definition).
On the one hand, its $\alpha$-limit set is the whole $\Lambda$.
On the other hand, this orbit must intersect the open set
$W^u(S)$, and so the $\alpha$-limit set must be contained in $S$.
This implies that $\Lambda\subset S$, which is a contradiction,
because $\Lambda$ contains singularities.
\end{proof}

\begin{corollary}\label{c.interior}
For any $x\in\Lambda$ there exist points $x^+\notin\Lambda$
and $x^-\notin\Lambda$ in distinct connected components of
$W^{ss}(x)\setminus\{x\}$.
\end{corollary}

\begin{proof}
Otherwise there would exist a whole segment of the strong-stable
manifold entirely contained in $\Lambda$. Considering any point
in the interior of this segment, we would get a contradiction to
Lemma~\ref{l.isolated}.
\end{proof}

\begin{lemma}\label{l.existeadaptada}
Let $x \in \Lambda$ be a regular point, that is, such that
$X(x)\neq 0$. Then there exists $\delta>0$ for which there
exists a $\delta$-adapted cross-section $\Sigma$ at $x$.
\end{lemma}

\begin{proof}
Fix $\vep>0$ as in the stable manifold theorem.
Any cross-section $\Sigma_0$ at $x$ sufficiently small with
respect to $\vep>0$ is foliated by the intersections
$W_\vep^s(x)\cap\Sigma_0$\,.
By Corollary~\ref{c.interior}, we may find points
$x^+\notin\Lambda$ and $x^-\notin\Lambda$ in each of the 
connected components of $W^s_\vep(x)\cap\Sigma_0$\,. 
Since $\Lambda$ is closed, there are neighborhoods $V^\pm$ 
of $x^\pm$ disjoint from $\Lambda$.
Let $\gamma\subset\Sigma_0$ be some small curve through $x$,
transverse to $W^s_\vep(x)\cap\Sigma_0$\,.
Then we may find a continuous family of segments inside
$W^s_\vep(y)\cap\Sigma_0$\,, $y\in\gamma$ with endpoints  
contained in $V^\pm$. The union $\Sigma$ of these segments 
is a $\delta$-adapted cross-section, for some $\delta>0$,
see Figure~\ref{fig.adaptedsection}.
\begin{figure}[ht]
\begin{center}
\includegraphics[height=4.5cm]{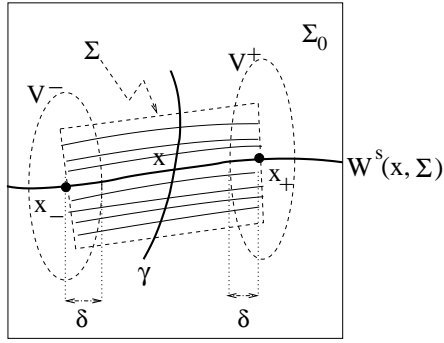}
\caption{\label{fig.adaptedsection} The construction of a
$\de$-adapted cross-section for a regular $x\in\Lambda$.}
\end{center}
\end{figure}
\end{proof}

We are going to show that if the cross-sections are adapted,
then we have the invariance property
\eqref{eq.stableMarkov}.
Given $\Sigma,\Sigma'\in\Xi$ we set 
$\Sigma(\Sigma')=\{ x\in\Sigma: R(x)\in\Sigma'\}$
the domain of the return map from $\Sigma$ to $\Sigma'$.

\begin{lemma}
\label{stablereturnmap}
Given $\delta>0$ and $\delta$-adapted cross-sections
$\Sigma$ and $\Sigma'$, there exists
$t_2=t_2(\Sigma,\Sigma')>0$ such that if
$R:\Sigma(\Sigma')\to\Sigma'$ defined by $R(z)=R_{t(z)}(z)$
is a Poincar\'e map with time $t(\cdot)>t_2$, then
\begin{enumerate}
\item $R\big(W^s(x,\Sigma)\big)\subset W^s(R(x),\Sigma')$
      for every $x\in\Sigma(\Sigma')$, and also
\item $d(R(y),R(z))\le \frac12 \, d(y,z)$ for every $y$, 
      $z\in W^s(x,\Sigma)$ and $x\in\Sigma(\Sigma')$.
\end{enumerate}
\end{lemma}


\begin{proof}
  This is a simple consequence of the relation
  \eqref{eq.contraction} from the proof of
  Proposition~\ref{p.secaohiperbolica}: the tangent
  direction to each $W^s(x,\Sigma)$ is contracted at an
  exponential rate $\lambda$
$$
\|DR(x) \, e_x^s \| \leq \frac{K}{\kappa} \lambda^{t(x)}.
$$ 
Choosing $t_2$ sufficiently large we ensure that 
$$
\frac{1}{\kappa} \lambda^{t_2}\cdot
\sup\{\ell( W^s(x,\Sigma)) : x\in\Sigma\} < \delta.
$$
In view of the definition of $\delta$-adapted cross-section this gives
part (1) of the lemma. Part (2) is entirely analogous: it suffices that
$(K/\kappa)\cdot \lambda^{t_2} < 1/2$.
\end{proof}

\begin{lemma}
\label{l.precisa?}
Let $\Sigma$ be a $\delta$-adapted cross-section. Then,
given any $r>0$ there exists $\rho$ such that
$$
d(y,z) < \rho \quad\Rightarrow\quad \dist(X_s(y),X_s(z))<r
$$
for all $s>0$, every $y$, $z\in W^s(x,\Sigma)$, and every
$x\in\Lambda\cap\Sigma$.
\end{lemma}

\begin{remark}
\label{r.angulos2}
Clearly we may choose $t_2>t_1$\,. Remark~\ref{r.angulos1} applies
to $t_2$ as well. 
\end{remark}

\begin{proof}
Let $y$ and $z$ be as in the statement. As in Remark~\ref{r.trapaca},
we may find $z'=X_{\tau}(z)$
in the intersection of the orbit of $z$ with the strong-stable 
manifold of $y$ satisfying
$$
\frac 1K \le \frac{\dist(y,z^\prime)}{d(y,z)}\le K
\quad\text{and}\quad
|\tau| \le K \cdot d(y,z).
$$
Then, given any $s>0$,
\begin{equation*}\begin{aligned}
\dist(X_s(y),X_s(z))
& \le \dist(X_s(y),X_s(z^\prime)) + \dist(X_s(z^\prime),X_s(z))
\\ & \le C \cdot e^{\gamma s} \cdot \dist(y,z^\prime) + 
\dist(X_{s+\tau}(z),X_s(z))
\\ & \le K C \cdot e^{\gamma s} \cdot d(y,z) + K |\tau|
\le \big(K C + K^2 \big)\cdot d(y,z).
\end{aligned}\end{equation*}
Taking $\rho<r/(K C + K^2)$ we get the statement of the lemma.
\end{proof}

\subsection{Flow boxes around singularities}
\label{s.24}

In this section we collect some known facts about the dynamics
near the singularities of the flow.  
It is known \cite[Theorem~A]{MPP99} that \emph{each singularity
of a singular-hyperbolic attracting set, accumulated by
regular orbits of a $3$-dimensional flow, must be Lorenz-like.}
In particular every singularity $\sigma_k$ of a singular-hyperbolic
\emph{attractor}, as in the setting of Theorem~\ref{sci}, is
Lorenz-like, that is,
the eigenvalues $\lambda_1\,, \lambda_2\,, \lambda_3$ of the
derivative $DX(\sigma_k)$ are all real and satisfy
$$
\lambda_1 > 0 > \lambda_2 > \lambda_3 \quand 
\lambda_1 + \lambda_2 > 0.
$$
In particular, the unstable manifold $W^u(\sigma_k)$ is
one-dimensional, and there is a one-dimensional
strong-stable manifold $W^{ss}(\sigma_k)$ contained in the
two-dimensional stable manifold $W^s(\sigma_k)$.  Most
important for what follows, the attractor intersects the
strong-stable manifold at the singularity only
\cite[Theorem~A]{MPP99}.

\begin{figure}[ht]
\begin{center}
\psfrag{S1}{$\Sigma^{i,+}$}
\psfrag{S2}{$\Sigma^{i,-}$}
\psfrag{S3}{$\Sigma^{o,+}$}
\psfrag{S4}{$\Sigma^{o,-}$}
\psfrag{s}{$\sigma_k$}
\psfrag{L1}{$\ell^+$}
\psfrag{L2}{$\ell^-$}
\psfrag{z}{{\footnotesize $z$}}
\psfrag{R}{{\footnotesize $R(z)$}}
\psfrag{X1}{{\footnotesize $x_1$}}
\psfrag{X2}{{\footnotesize $x_2$}}
\psfrag{X3}{{\footnotesize $x_3$}}
\includegraphics[height=1.5in]{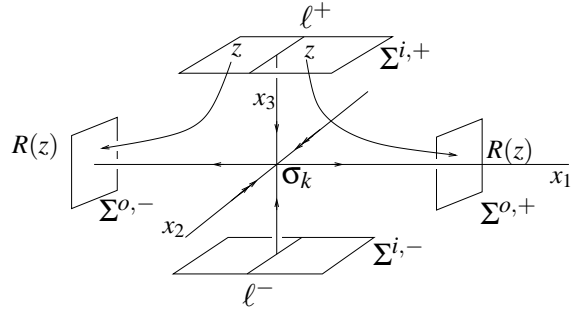}
\caption{\label{f.singularbox}Ingoing and outgoing
  adapted cross-sections near a singularity.}
\end{center}
\end{figure}

Then for some $\de>0$ we may choose $\de$-adapted
cross-sections contained in $U_0$
\begin{itemize}
\item $\Sigma^{o,\pm}$ at points $y^{\pm}$ in different components
of $W^u_{loc}(\sigma_k)\setminus\{\sigma_k\}$ 
\item $\Sigma^{i,\pm}$ at points $x^{\pm}$ in different components
of $W^s_{loc}(\sigma_k)\setminus W^{ss}_{loc}(\sigma_k)$
\end{itemize}
and Poincar\'e maps
$R^\pm:\Sigma^{i,\pm}\setminus\ell^\pm\to
\Sigma^{o,-}\cup\Sigma^{o,+}$,
where $ \ell^\pm=\Sigma^{i,\pm}\cap W^s_{loc}(\sigma_k), $
satisfying (see Figure~\ref{f.singularbox})
\begin{enumerate}
\item every orbit in the attractor passing through a small
  neighborhood of the singularity $\sigma_k$ intersects some
  of the incoming cross-sections $\Sigma^{i,\pm}$;
\item $R^\pm$ maps each connected component of
  $\Sigma^{i,\pm}\setminus\ell^\pm$ diffeomorphically inside
  a different outgoing cross-section $\Sigma^{o,\pm}$,
  preserving the corresponding stable foliations and
  unstable cones.
\end{enumerate}

These cross-sections may be chosen to be planar relative to
some linearizing system of coordinates near $\sigma_k$\,,
e.g. for a small $\delta>0$
\[
\Sigma^{i,\pm}=\{ (x_1,x_2,\pm1): |x_1|\le\delta , |x_2|\le\delta
\}
\quad\mbox{and}\quad 
\Sigma^{o,\pm}=\{ (\pm1,x_2,x_3): |x_2|\le\delta ,
|x_3|\le\delta\},
\] 
where the $x_1$-axis corresponds to the unstable manifold
near $\sigma_k$, the $x_2$-axis to the strong-stable
manifold and the $x_3$-axis to the weak-stable manifold of
the singularity which, in turn, is at the origin, see
Figure~\ref{f.singularbox}.

Reducing the cross-sections if necessary, i.e. taking
$\de>0$ small enough, we ensure that the Poincar\'e times
are larger than $t_2$\,, so that the same conclusions as in
the previous sections apply here. Indeed using linearizing
coordinates it is easy to see that for points
$z=(x_1,x_2,\pm1)\in\Sigma^{i,\pm}$ the time $\tau^\pm$ it
takes the flow starting at $z$ to reach one of
$\Sigma^{o,\pm}$ depends on $x_1$ only and is given by
\[
\tau^\pm(x_1)=-\frac{\log x_1}{\lambda_1}.
\]
We then fix these
cross-sections once and for all and define for small
$\vep>0$ the \emph{flow-box}
\[
U_{\sigma_k}=\bigcup_{x\in\Sigma^{i,\pm}\setminus\ell^{\pm}}
X_{(-\vep,\tau^{\pm}(x)+\vep)}(x)
\cup(-\de,\de)\times(-\de,\de)\times(-1,1)
\]
which is an open neighborhood of $\sigma_k$ with $\sigma_k$
the unique zero of $X\mid U_{\sigma_k}$. We note that the
function $\tau^\pm:\Sigma^{i,\pm}\to\real$ is integrable with
respect to the Lebesgue (area) measure over
$\Sigma^{i,\pm}$: we say that \emph{the exit time function
  in a flow box near each singularity is Lebesgue integrable}.

\medskip

In particular we can determine the expression of the
Poincar\'e maps between ingoing and outgoing cross-sections
easily thought linearized coordinates
\begin{equation}
  \label{eq:nonflatsing}
\Sigma^{i,+}\cap\{x_1>0\}\to \Sigma^{0,+},
\quad
(x_1,x_2,1)\mapsto 
\big(1,x_2\cdot x_1^{-\lambda_3/\lambda_1}, 
x_1^{-\lambda_2/\lambda_1}\big).  
\end{equation}
This shows that the map obtained identifying points with the
same $x_2$ coordinate, i.e. points in the same stable leaf,
is simply $x_1\mapsto x_1^{\beta}$ where
$\beta=-\lambda_2/\lambda_1\in(0,1)$. For the other possible
combinations of ingoing and outgoing cross-sections the
Poincar\'e maps have a similar expression. This will be
useful to construct physical measures for the flow.

\section{Proof of expansiveness}
\label{sec:proof-expansiveness}

Here we prove Theorem~\ref{sci}. The proof is by contradiction:
let us suppose that there exist $\vep>0$, a sequence 
$\delta_n\to 0$, a sequence of functions $h_n\in S(\real)$,
and sequences of points $x_n,\, y_n \in\Lambda$ such that
\begin{equation}
\label{eq.rel1}
d\big(X_t(x_n),X_{h_n(t)}(y_n)\big)\le\delta_n
\quad\text{for all } t\in\real,
\end{equation}
but
\begin{equation}
\label{eq.rel2}
X_{h_n(t)}(y_n) \notin X_{[t-\vep,t+\vep]}(x_n)
\quad\text{for all } t\in\real.
\end{equation}

\subsection{Proof of Theorem \ref{sci}}

The main step in the proof is a reduction to a forward
expansiveness statement about Poincar\'e maps which
we state in Theorem~\ref{t.expansivepoincare} below.

We are going to use the following observation: there exists
some regular (i.e. non-equilibrium) point $z\in\Lambda$
which is accumulated by the sequence of $\omega$-limit sets
$\omega(x_n)$.  To see that this is so, start by observing
that accumulation points do exist, since the ambient space
is compact.  Moreover, if the $\omega$-limit sets accumulate
on a singularity then they also accumulate on at least one
of the corresponding unstable branches which, of course,
consists of regular points.  We fix such a $z$ once and for
all. Replacing our sequences by subsequences, if necessary,
we may suppose that for every $n$ there exists
$z_n\in\omega(x_n)$ such that $z_n\to z$.

Let $\Sigma$ be a $\delta$-adapted cross-section at $z$, for
some small $\delta$. Reducing $\delta$ (but keeping the same
cross-section) we may ensure that $z$ is in the interior of
the subset 
$$
\Sigma_\delta=\{y\in\Sigma: d(y,\partial\Sigma)>\delta\}.
$$ 
By definition the orbit of $x_n$ returns infinitely often to
a neighborhood of $z_n$ which, on its turn, is close to
$z$.  Thus dropping a finite number of terms in our
sequences if necessary, we have that the orbit of $x_n$
intersects $\Sigma$ infinitely many times.  Let $t_n$
be the time corresponding to the first intersection.
Replacing $x_n$, $y_n$, $t$, and $h_n$ by
$x^{(n)}=X_{t_n}(x_n)$, $y^{(n)}=X_{h_n(t_n)}(y_n)$, $t'=t-t_n$,
and $h_n'(t')=h_n(t'+t_n)-h_n(t_n)$, we may suppose that
$x^{(n)}\in\Sigma_\delta$\,, while preserving both relations
\eqref{eq.rel1} and \eqref{eq.rel2}.  Moreover there exists
a sequence $\tau_{n,j}$\,, $j\ge 0$ with $\tau_{n,0}=0$ such
that
\begin{equation}
\label{eq.xnj}
x^{(n)}(j)=X_{\tau_{n,j}}(x^{(n)})\in\Sigma_\delta
\quad\text{and}\quad
\tau_{n,j}-\tau_{n,j-1}>\max\{t_1,t_2\}
\end{equation}
for all $j\ge 1$, where $t_1$ is given by
Proposition~\ref{p.secaohiperbolica} and $t_2$ is given by
Lemma~\ref{stablereturnmap}.

\begin{theorem}
\label{t.expansivepoincare}
Given $\vep_0>0$ there exists $\delta_0>0$ such that if
$x\in\Sigma_\delta$ and $y\in\Lambda$ satisfy
\begin{itemize}
\item[(a)] there exist $\tau_j$ such that
\[
x_j=X_{\tau_j}(x)\in\Sigma_\delta
\quad\text{and}\quad
\tau_{j}-\tau_{j-1}>\max\{t_1,t_2\}
\quad\text{for all $j\ge 1$};
\]
\item[(b)] $\dist \big(X_t(x),X_{h(t)}(y)\big) < \delta_0$, for all $t>0$
and some $h\in S(\real)$;
\end{itemize}
then there exists $s\in\real$ such that
$X_{h(s)}(y)\in W_{\vep_0}^{ss}(X_{[s-\vep_0,s+\vep_0]}(x))$.
\end{theorem}


We postpone the proof of Theorem~\ref{t.expansivepoincare}
until the next section and explain first why it implies Theorem \ref{sci}.
We are going to use the following observation.

\begin{lemma}\label{proximo}
There exist $\rho>0$ small and $c>0$, depending only on the flow,
such that if $z_1, z_2, z_3$ are points in $\Lambda$ satisfying
$z_3\in X^{[-\rho,\rho]}(z_2)$ and $z_2\in W_\rho^{ss}(z_1)$, then
\[
\dist(z_1,z_3) \ge c \cdot \max\{\dist(z_1,z_2),\dist(z_2,z_3)\}.
\]
\end{lemma}

\begin{proof}
  This is a direct consequence of the fact that the angle
  between $E^{ss}$ and the flow direction is bounded from
  zero which, on its turn, follows from the fact that the
  latter is contained in the center-unstable sub-bundle
  $E^{cu}$. Indeed consider for small enough $\rho>0$ the
  $C^1$ surface
  $X^{[-\rho,\rho]}\big(W^{ss}_\rho(z_1)\big)$. The
  Riemannian metric here is uniformly close to the Euclidean
  one and we may choose coordinates on $[-\rho,\rho]^2$
  putting $z_1$ at the origin, sending $W^{ss}_\rho(z_1)$ to
  the segment $[-\rho,\rho]\times\{0\}$ and
  $X^{[-\rho,\rho]}(z_1)$ to $\{0\}\times[-\rho,\rho]$, see
  Figure~\ref{fig:proxima}.
  \begin{figure}[ht]
    \centering
    \psfrag{l}{$-\rho$}\psfrag{r}{$\rho$}\psfrag{a}{$z_1$}
    \psfrag{b}{$z_2$}\psfrag{c}{$z_3$}
    \psfrag{X}{$X^{[-\rho,\rho]}(z_1)$}
    \psfrag{W}{$W^{ss}_\rho(z_1)$}
    \includegraphics[width=4cm]{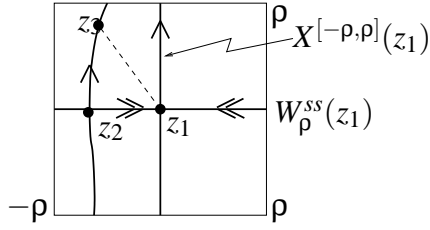}
    \caption{Distances near a point in the stable-manifold.}
    \label{fig:proxima}
  \end{figure}
Then the angle $\alpha$ between
  $X^{[-\rho,\rho]}(z_2)$ and the horizontal is bounded from
  below away from zero and the existence of $c$ follows by
  standard arguments using the Euclidean metric.
\end{proof}

\begin{figure}[ht]
\begin{center}
\includegraphics[height=4.5cm]{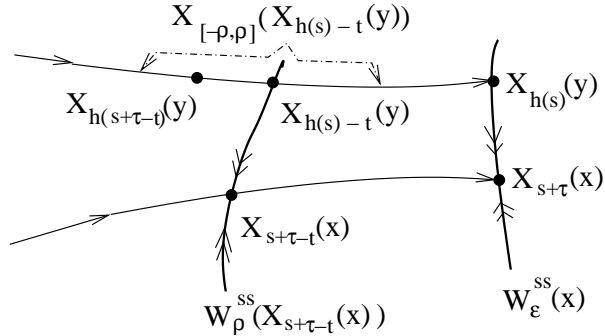}
\caption{\label{fig.contradiction} Sketch of the relative
  positions of the strong-stable manifolds and orbits in the
  argument reducing Theorem~\ref{sci} to
  Theorem~\ref{t.expansivepoincare}.}
\end{center}
\end{figure}

We fix $\vep_0=\vep$ as in \eqref{eq.rel2} and then consider
$\delta_0$ as given by Theorem~\ref{t.expansivepoincare}.
Next, we fix $n$ such that $\delta_n<\delta_0$ and $\delta_n <
c\rho$, and apply Theorem~\ref{t.expansivepoincare} to
$x=x^{(n)}$ and $y=y^{(n)}$ and $h=h_n$\,.  Hypothesis (a) in the
theorem corresponds to \eqref{eq.xnj} and, with these
choices, hypothesis (b) follows from \eqref{eq.rel1}.
Therefore we obtain that $X_{h(s)}(y)\in
W_\vep^{ss}(X_{[s-\vep,s+\vep]}(x))$.  In other words, there
exists $|\tau|\le\vep$ such that $X_{h(s)}(y)\in
W_\vep^{ss}(X_{s+\tau}(x))$.  Hypothesis \eqref{eq.rel2}
implies that $X_{h(s)}(y) \neq X_{s+\tau}(x)$.  Since
strong-stable manifolds are expanded under backward
iteration, there exists $\theta>0$ maximum such that
$$
X_{h(s)-t}(y)\in W^{ss}_\rho(X_{s+\tau-t}(x))
\quad\text{and}\quad
X_{h(s+\tau-t)}(y)\in X_{[-\rho,\rho]}(X_{h(s)-t}(y))
$$
for all $0\le t\le \theta$, see
Figure~\ref{fig.contradiction}. Since $\theta$ is maximum
\begin{align*}
\text{either  }  
\dist\big( X_{h(s)-t}(y) , X_{s+\tau-t}(x) \big) = \rho,\,
\text{  or  }
\dist\big( X_{h(s+\tau-t)}(y) , X_{h(s)-t}(y) \big) = \rho
\text{   for  } t=\theta.
\end{align*}

Using Lemma~\ref{proximo}, we conclude 
that 
$$
\dist(X_{s+\tau-t}(x),X_{h(s+\tau-t)}(y)) \ge c \rho >\delta_n
$$ 
which contradicts \eqref{eq.rel1}.
This contradiction reduces the proof of Theorem~\ref{sci}
to that of Theorem~\ref{t.expansivepoincare}.


\subsection{Infinitely many coupled returns}
\label{s.infinitely}

We start by outlining the proof of
Theorem~\ref{t.expansivepoincare}.  There are three steps.
\begin{itemize}
\item The first one, which we carry out in the present
  section, is to show that to each return $x_j$ of the orbit
  of $x$ to $\Sigma$ there corresponds a nearby return $y_j$
  of the orbit of $y$ to $\Sigma$.  The precise statement is
  in Lemma~\ref{l.returny} below.
\item The second, and most
  crucial step, is to show that there exists a smooth
  Poincar\'e map, with large return time, defined on the
  whole strip of $\Sigma$ in between the stable manifolds of
  $x_j$ and $y_j$\,. This is done in
  Section~\ref{s.poincaretube}.
\item The last step,
  Section~\ref{s.endoftheproof}, is to show that these
  Poincar\'e maps are uniformly hyperbolic, in particular,
  they expand $cu$-curves uniformly (recall the definition
  of $cu$-curve in Section~\ref{s.22}).
\end{itemize}
The theorem is then easily deduced: to prove that
$X_{h(s)}(y)$ is in the orbit of $W_\vep^{ss}(x)$ it
suffices to show that $y_j\in W^s(x_j,\Sigma)$, by
Remark~\ref{r.trapaca}. The latter must be true, for
otherwise, by hyperbolicity of the Poincar\'e maps, the
stable manifolds of $x_j$ and $y_j$ would move apart as
$j\to\infty$, and this would contradict condition (b) of
Theorem~\ref{t.expansivepoincare}. See
Section~\ref{s.endoftheproof} for more details.


\begin{lemma}
\label{l.returny}
There exists $K>0$ such that, in the setting of
Theorem~\ref{t.expansivepoincare}, there 
exists a sequence $(\upsilon_j)_{j\ge 0}$ such that 
\begin{enumerate}
\item $y_j=X_{\upsilon_j}(y)$ is in $\Sigma$ for all $j\ge 0$;
\item $|\upsilon_j-h(\tau_j)|<K\cdot\delta_0$, and
\item $d(x_j,y_j)<K\cdot\delta_0$.
\end{enumerate}
\end{lemma}

\begin{proof}
  By assumption $d(x_j,X_{h(\tau_j)}(y))<K\cdot \delta_0$ for all
  $j\ge 0$.  In particular $y'_j=X_{h(\tau_j)}(y)$ is close
  to $\Sigma$.  Using a flow box in a neighborhood of
  $\Sigma$ we obtain $X_{\epsilon_j}(y_j')\in\Sigma$
  for some $\epsilon_j\in(-K\cdot\delta_0,K\cdot\delta_0)$.  The constant
  $K$ depends only on the vector field $X$ and the
  cross-section $\Sigma$ (more precisely, on the angle
  between $\Sigma$ and the flow direction). Taking
  $\upsilon_j=h(\tau_j)+\epsilon_j$ we get the first two
  claims in the lemma. The third one follows from the
  triangle inequality; it may be necessary to replace $K$ by
  a larger constant, still depending on $X$ and $\Sigma$
  only.
\end{proof}


\subsection{Semi-global Poincar\'e map}
\label{s.poincaretube}

Since we took the cross-section $\Sigma$ to be adapted, we may 
use Lemma~\ref{stablereturnmap} to conclude that there exist
Poincar\'e maps $R_j$ with $R_j(x_j)=x_{j+1}$ and
$R_j(y_j)=y_{j+1}$ and sending $W^{s}_\vep(x_j\,,\Sigma)$ and
$W^{s}_\vep(y_j\,,\Sigma)$ inside $W^{s}_\vep(x_{j+1}\,,\Sigma)$
and $W^{s}_\vep(y_{j+1}\,,\Sigma)$, respectively.
The goal of this section is to prove that \emph{$R_j$ extends to a
smooth Poincar\'e map on the whole strip $\Sigma_j$ of
$\Sigma$ bounded by the stable manifolds of $x_j$ and $y_j$}\,.

\begin{figure}[htb]
\begin{center}
\psfrag{T}{$\mathcal T$}
\includegraphics[height=5cm]{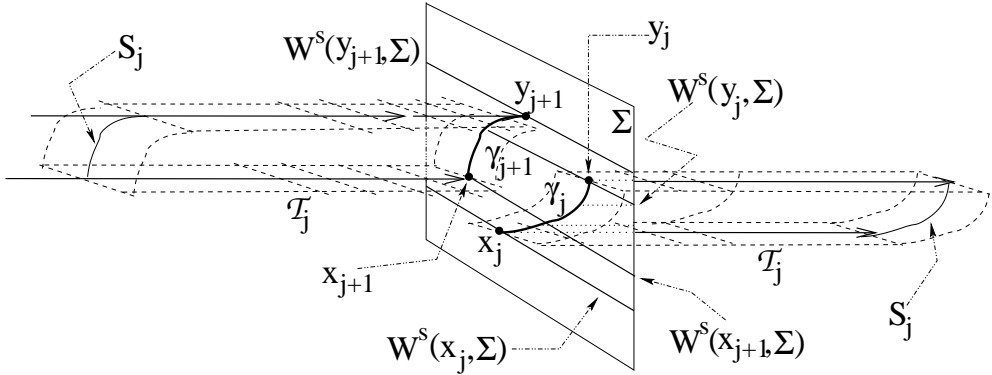}
\caption{\label{fig.tube} A tube-like domain.}
\end{center}
\end{figure}

We first outline the proof. For each $j$ we choose a curve
$\gamma_j$ transverse to the stable foliation of $\Sigma$,
connecting $x_j$ to $y_j$ and such that $\gamma_j$ is
disjoint from the orbit segments $[x_{j}\,,x_{j+1}]$ and
$[y_{j}\,,y_{j+1}]$.  Using Lemma~\ref{stablereturnmap} in
the same way as in the last paragraph, we see that it
suffices to prove that $R_j$ extends smoothly to
$\gamma_j$\,.  For this purpose we consider a tube-like
domain $\cT_j$ consisting of local stable manifolds through
an \emph{immersed surface $S_j$} 
whose boundary is formed by $\gamma_j$ and $\gamma_{j+1}$
and the orbit segments $[x_j\,, x_{j+1}]$ and
$[y_j\,,y_{j+1}]$\,, see Figure~\ref{fig.tube}.  We will
prove that the orbit of any point in $\gamma_j$ must leave
the tube through $\gamma_{j+1}$ in finite time.
We begin by showing that the tube contains no singularities.
This uses hypothesis (b) together with the local dynamics
near Lorenz-like singularities.  Next, using hypothesis (b)
together with a Poincar\'e-Bendixson argument on $S_j$\,, we
conclude that the forward orbit of any point in $\cT_j$ must
leave the tube.  Another argument, using hyperbolicity
properties of the Poincar\'e map, shows that orbits through
$\gamma_j$ must leave $\cT_j$ through $\gamma_{j+1}$\,.  In
the sequel we detail these arguments.


\subsubsection{A tube-like domain without singularities}
\label{sec:tube-like-domain}
\begin{figure}[htb]
\begin{center}
\psfrag{T}{$\mathcal T$}
\psfrag{l}{$\ell$}
\includegraphics[height=5cm]{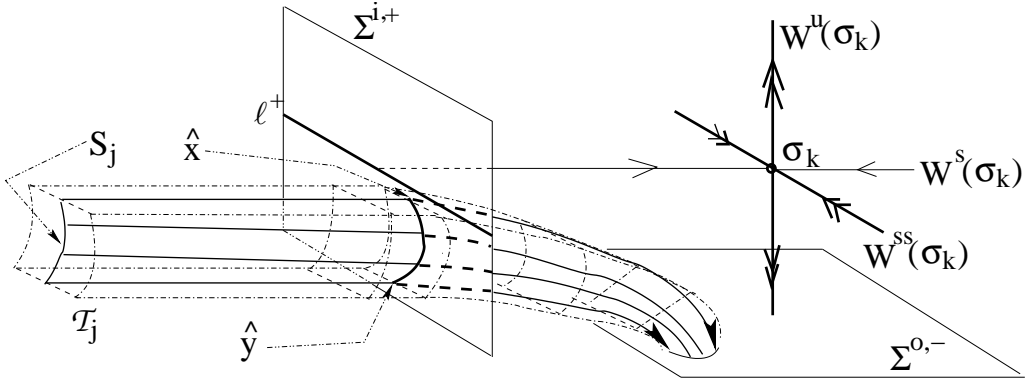}
\caption{\label{fig.tubesingular} Entering the flow box of a singularity.}
\end{center}
\end{figure}

Since we took $\gamma_j$ and $\gamma_{j+1}$ disjoint from 
the orbit segments $[x_j\,,x_{j+1}]$ and $[y_j\,,y_{j+1}]$,
the union of these four curves is an embedded circle.
We recall that the two orbit segments are
close to each other, by hypothesis (b)
$$
d(X_t(x),X_{h(t)}(y))<\delta_0
 \quad\text{for all $t\in [t_j,t_{j+1}]$.}  
 $$
 Assuming that $\delta_0$ is smaller than the radius of
 injectiveness of the exponential map of the ambient
 manifold (i.e. $\exp_x:T_xM\to M$ is locally invertible in
 a $\delta_0$-neighborhood of $x$ in $M$ for any $x\in M$),
 there exists a unique geodesic linking each $X_t(x)$ to
 $X_{h(t)}(y)$, and it varies continuously (even smoothly)
 with $t$.  Using these geodesics we easily see that the
 union of $[y_j\,,y_{j+1}]$ with $\gamma_j$ and
 $\gamma_{j+1}$ is homotopic to a curve inside the orbit of
 $x$, with endpoints $x_j$ and $x_{j+1}$, and so it is also
 homotopic to the segment $[x_j,x_{j+1}]$.  This means that
 the previously mentioned embedded circle is homotopic to
 zero.  It follows that there is a \emph{smooth immersion}
 $\phi:[0,1]\times[0,1] \to M$ such that
\begin{itemize}
\item $\phi(\{0\}\times[0,1])=\gamma_j$ and
      $\phi(\{1\}\times[0,1])=\gamma_{j+1}$;
\item $\phi([0,1]\times\{0\})=[y_j\,,y_{j+1}]$ and
      $\phi([0,1]\times\{1\})=[x_j\,,x_{j+1}]$.
\end{itemize}
Moreover $S_j=\phi([0,1]\times[0,1])$ may be chosen such
that
\begin{itemize}
\item all the points of $S_j$ are at distance less than
$\delta_1$ from the orbit segment $[x_j\,,x_{j+1}]$,
for some uniform constant $\delta_1 >\delta_0$ which can
be taken arbitrarily close to zero, reducing $\delta_0$
if necessary, see Figure~\ref{fig.tube};
\item the intersection of $S_j$ with an incoming
  cross-section of any singularity (Section~\ref{s.24}) is
  transverse to the corresponding stable foliation, see
  Figure~\ref{fig.tubesingular}.
\end{itemize}
Then we define $\cT_j$ to be the union of the local
stable manifolds through the points of that disk.

\begin{proposition}
\label{p.nosingularities}
The domain $\cT_j$ contains no singularities of the flow.
\end{proposition}

\begin{proof}
By construction, every point of $\cT_j$ is at distance
$\le\vep$ from $S_j$ and, consequently, at distance
$\le \vep+\delta_1$ from $[x_j\,,x_{j+1}]$.
So, taking $\vep$ and $\delta_0$ much smaller than the sizes
of the cross-sections associated to the singularities
(Section~\ref{s.24}), we immediately
get the conclusion of the proposition in the case 
when $[x_j\,,x_{j+1}]$ is disjoint from the incoming
cross-sections of all singularities.
In the general case we must analyze the intersections
of the tube with the flow boxes at the singularities.
The key observation is in the following statement whose
proof we postpone.

\begin{lemma}
\label{l.sameside}
Suppose $[x_j\,,x_{j+1}]$ intersects an incoming cross-section
$\Sigma_k^i$ of some singularity $\sigma_k$ at some point 
$\hat x$ with $d(\hat x,\partial\Sigma_k^i)>\delta$.  
Then $[y_j\,,y_{j+1}]$ intersects $\Sigma_k^i$ at some point
$\hat y$ with $d(\hat x,\hat y)<K\cdot\delta_0$ and, moreover
$\hat x$ and $\hat y$ are in the same connected
component of $\Sigma_k^i\setminus W_{loc}^s(\sigma_k)$.
\end{lemma}

Let us recall that by construction the intersection of $S_j$
with the incoming cross-section $\Sigma_k^i$ is transverse
to the corresponding stable foliation,
see Figure~\ref{fig.tubesingular}.
By the previous lemma this intersection is entirely contained
in one of the connected components of
$\Sigma_k^i\setminus W^s_{loc}(\sigma_k)$.
Since $\cT_j$ consists of local stable manifolds through
the points of $S_j$\, its intersection with $\Sigma_k^i$
is contained in the region bounded by the stable manifolds
$W^s(\hat x,\Sigma_k^i)$ and $W^s(\hat y,\Sigma_k^i)$,
and so it is entirely contained in a connected component
of $\Sigma_k^i\setminus W^s_{loc}(\sigma_k)$.
In other words, the crossing of the tube $\cT_j$ through
the flow box is disjoint from $W^s_{loc}(\sigma_k)$, in
particular, it does not contain the singularity.
Repeating this argument for every intersection of the 
tube with a neighborhood of some singularity, we get the
conclusion of the proposition.
\end{proof}

\begin{proof}[Proof of Lemma~\ref{l.sameside}]
The first part is proved in exactly the same way as 
Lemma~\ref{l.returny}. We have 
$$
\hat x = X_{r_0}(x) \quad\text{and}\quad \hat y = X_{s_0}(y)
$$
with $|s_0-h(r_0)|<K\delta_0$\,.
The proof of the second part is by contradiction and relies,
fundamentally, on the local description of the dynamics near
the singularity.
Associated to $\hat x$ and $\hat y$ we have the points
$\tilde x = X_{r_1}(x)$ and $\tilde y = X_{s_1}(y)$,
where the two orbits leave the flow box associated to the
singularity. If $\hat x$ and $\hat y$ are in opposite sides
of the local stable manifold of $\sigma_k$, then $\tilde x$
and $\tilde y$ belong to different outgoing cross-sections
of $\sigma_k$\,. Our goal is to find some $t\in\real$ such
that
$$
\dist\big(X_t(x),X_{h(t)}(y)\big) > \delta_0\,,
$$
thus contradicting hypothesis (b).

We assume by contradiction that $\hat x, \hat y$ are in
different connected components of $\Sigma_k^{i,\pm}\setminus
\ell^{\pm}$.  There are two cases to consider. We suppose
first that $h(r_1)> s_1$ and note that $s_1 \gg s_0 \approx
h(r_0)$, so that $s_1>h(r_0)$.  It follows that there exists
$t\in (r_0,r_1)$ such that $h(t)=s_1$ since $h$ is
non-decreasing and continuous.  Then $X_t(x)$ is on one side
of the flow box of $\sigma_k$\,, whereas $X_{h(t)}(y)$
belongs to the outgoing cross-section at the other side of
the flow box. Thus $\dist\big(X_t(x),X_{h(t)}(y)\big)$ has
the order of magnitude of the diameter of the flow box,
which we may assume to be much larger than $\delta_0$\,.

Now we suppose that $s_1 \ge h(r_1)$ and observe that
$h(r_1)>h(r_0)$, since $h$ is increasing. We recall also
that $X_{h(r_0)}(y)$ is close to $\hat y$, near the incoming
cross-section, so that the whole orbit segment from
$X_{h(r_0)}(y)$ to $X_{s_1}(y)$ is contained in (a small
neighborhood of) the flow box, to one side of the local
stable manifold of $\sigma_j$\,. The previous observation
means that this orbit segment contains $X_{h(r_1)}(y)$.
However $X_{r_1}(x)$ belongs to the outgoing cross-section
at the opposite side of the flow box, and so
$\dist\big(X_{r_1}(x),X_{h(r_1)}(y)\big)$ has the order of
magnitude of the diameter of the flow box, which is much
larger than $\delta_0$\,.
\end{proof}


\subsubsection{Every orbit leaves the tube}
\label{sec:every-orbit-leaves}
Our goal in this section is to show that the forward orbit
of every point $z\in\cT_j$ leaves the tube in finite time.
The proof is based on a Poincar\'e-Bendixson argument
applied to the flow induced by $X^t$ on the immersed disk
$S_j$\,.

We begin by defining this induced flow.
For the time being, we make the following simplifying assumption:
\begin{itemize}
\item[(H)] $S_j=\phi([0,1]\times[0,1])$ is an embedded disk and
the stable manifolds $W^s_\vep(\xi)$ through the points
$\xi\in S_j$ are pairwise disjoint.
\end{itemize}
This condition provides a well-defined continuous projection
$\pi:\cT_j\to S_j$ by assigning to each point $z\in \cT_j$
the unique $\xi\in S_j$ whose local strong-stable manifold
contains $z$.  The (not necessarily complete) flow $Y_t$
induced by $X_t$ on $S_j$ is given by
$Y_t(\xi)=\pi(X_t(\xi))$ for the largest interval of
values of $t$ for which this is defined.  It is clear, just
by continuity, that given any subset $E$ of $S_j$ at a
positive distance from $\partial S_j$\,, there exists
$\vep>0$ such that $Y_t(\xi)$ is defined for all $\xi\in
E$ and $t\in[0,\vep]$. In fact this remains true even if $E$
approaches the curve $\gamma_j$ (since $\Sigma$ is a
cross-section for $X_t$\,, the flow at $\gamma_j$ points
inward $S_j$) or the $X_t$-orbit segments $[x_j\,,x_{j+1}]$
and $[y_j\,,y_{j+1}]$ on the boundary of $S_j$ (because they
are also $Y_t$-orbit segments).  Thus we only have to
worry with the distance to the remaining boundary segment:
\begin{itemize}
\item[(U)] given any subset $E$ of $S_j$ at positive distance
from $\gamma_{j+1}$\,, there exists $\vep>0$ such that $Y_t(\xi)$
is defined for all $\xi\in E$ and $t\in[0,\vep]$.
\end{itemize}
We observe also that for points $\xi$ close to $\gamma_{j+1}$ the
flow $Y_t(\xi)$ must intersect $\gamma_{j+1}$\,, after which
it is no longer defined.

\emph{Now we explain how to remove condition (H).} In this
case, the induced flow is naturally defined on
$[0,1]\times[0,1]$ rather than $S_j$\,, as we now explain.
We recall that $\phi:[0,1]\times [0,1]\to M$ is an
\emph{immersion}.  So given any $w\in[0,1]\times[0,1]$ there
exist neighborhoods $U$ of $w$ and $V$ of $\phi(w)$ in $S_j$
such that $\phi:U\to V$ is a diffeomorphism. Moreover, just
by continuity of the stable foliation, choosing $V$
sufficiently small we may ensure that each strong-stable
manifold $W_\vep^{ss}(\xi)$, $\xi\in V$, intersects $V$ only
at the point $\xi$. This means that we have a well-defined
projection $\pi$ from $\cup_{\xi\in V}W_\vep^{ss}(\xi)$ to
$V$ associating to each point $z$ in the domain the unique
element of $V$ whose stable manifold contains $z$. Then we
may define $Y_t(w)$ for small $t$, by
$$
Y_t(w)= \phi^{-1}(\pi(X_t(\phi(w))).
$$
As before, we extend $Y_t$ to a maximal domain.
This defines a (partial) flow on the square $[0,1]\times[0,1]$,
such that both $[0,1]\times\{i\}$, $i\in\{0,1\}$ 
are trajectories. 

\begin{remark}
  \label{rmk:sing-Yt}
  A singularity $\zeta$ for the flow $Y_t$ corresponds to a
  singularity of $X$ in the local strong-stable manifold of
  $\zeta$ in $M$ by the definition of $Y_t$ through the
  projection $\pi$.
\end{remark}

Notice also that forward trajectories
of points in $\{0\}\times[0,1]$ enter the square. Hence, the
only way trajectories may exit is through $\{1\}\times [0,1]$.
So, we have the following reformulation of property (U):
\begin{itemize}
\item[(U)] given any subset $E$ of $[0,1]\times[0,1]$ at positive
distance from $\{1\}\times[0,1]$, there exists $\vep>0$ such that
$Y_t(w)$ is defined for all $w\in E$ and $t\in[0,\vep]$.
\end{itemize}
Moreover for points $w$ close to $\{1\}\times[0,1]$ the
flow $Y_t(\xi)$ must intersect $\{1\}\times[0,1]$, after which
it is no longer defined.

\begin{proposition}
\label{p.escapadotubo}
Given any point $z\in\cT_j$ there exists $t>0$ such that 
$X_t(z)\notin\cT_j$\,.
\end{proposition}

\begin{proof}
The proof is by contradiction. First, we assume condition (H).
Suppose there exists $z\in\cT_j$ whose forward orbit remains in the
tube for all times. Let $z_0=\pi(z)$. Then $Y_t(z_0)$ is 
defined for all $t>0$, and so it makes sense to speak of the 
$\omega$-limit set $\omega(z_0)$. The orbit $Y_t(z_0)$ can
not accumulate on $\gamma_{j+1}$ for otherwise it would leave $S_j$\,.
Therefore $\omega(z_0)$ is a compact subset of $S_j$ at positive
distance from $\gamma_{j+1}$.
Using property (U) we can find a uniform constant $\vep>0$ such
that $Y_t(w)$ is defined for every $t\in[0,\vep]$ and every
$w\in\omega(z_0)$. Since $\omega(z_0)$ is an invariant
set, we can extend $Y_t$ to a complete flow on it.

In particular we may fix $w_0\in\omega(z_0)$,
$w\in\omega(w_0)$ and apply the arguments in the proof
of the Poincar\'e-Bendixson Theorem.  On the one hand, if we
consider a cross-section $S$ to the flow at $w$, the
forward orbits of $z_0$ and $w_0$ must intersect it on
monotone sequences; on the other hand, every intersection of
the orbit of $w_0$ with $S$ is accumulated by points in the
orbit of $z_0$. This implies that $w$ is in the orbit of
$w_0$ and, in fact, that the later is periodic.  

We consider the disk $D\subset S_j$ bounded by the orbit of
$w_0$.  The flow $Y_t$ is complete restricted to $D$ and so
we may apply Poincar\'e-Bendixson's Theorem (see
\cite{PM82}) once more, and conclude that $Y_t$ has some
singularity $\zeta$ inside $D$.  This implies by
Remark~\ref{rmk:sing-Yt} that $X_t$ has a singularity in the
local stable manifold of $\zeta$\,, which contradicts
Proposition~\ref{p.nosingularities}.  This contradiction
completes the proof of the proposition, under assumption
(H). The general case is treated in the same way, just
dealing with the flow induced on $[0,1]\times[0,1]$ instead
of on $S_j$\,.
\end{proof}


\subsubsection{The Poincar\'e map is well-defined on $\Sigma_j$}
\label{sec:poincare-map-well}
We have shown that for the induced flow $Y_t$ on $S_j$ (or, 
more generally, on $[0,1]\times[0,1]$) every orbit must
eventually cross $\gamma_{j+1}$ (respectively, $\{1\}\times[0,1]$).
Hence there exists a continuous Poincar\'e map
$$
r:\gamma_j\to\gamma_{j+1}, \quad r(\xi)=Y_{\theta(\xi)}(\xi).
$$
By compactness the Poincar\'e time $\theta(\cdot)$ is bounded.
We are going to deduce that every forward $X_t$-orbit eventually
leaves the tube $\cT_j$ through $\Sigma_{j+1}$ (the strip in
$\Sigma$ between the stable manifolds $W^s(x_{j+1},\Sigma)$
and $W^s(y_{j+1},\Sigma)$), which proves 
that $R_j$ is defined on the whole strip of $\Sigma_j$ between the
manifolds $W^s(x_j,\Sigma_j)$ and $W^s(y_j,\Sigma_j)$, as
claimed in Section~\ref{s.infinitely}.

To this end, let $\gamma$ be a \emph{central-unstable curve in}
$\Sigma_\delta$ \emph{connecting the stable manifolds}
$W^s(x_j,\Sigma)$ and $W^s(y_j,\Sigma)$. Observe that
$\gamma$ is inside $\cT_j$.  For each $z\in\gamma$, let
$t(z)$ be the smallest positive time for which $X_{t(z)}$ is
on the boundary of $\cT_j$.

The crucial observation is that, in view of the construction
of $Y_t$\,, \emph{each $X_{t(\xi)}(\xi)$ belongs to the
  stable manifold of $Y_{t(z)}\big(\pi(z)\big)$.}
We observe also that for $\{\xi\}=\gamma\cap W^s(x_j,\Sigma)$ we
have $Y_t(\xi)=X_t(\xi)$ and so $t(\xi)=\theta(\xi)$.

Now we take $z\in\gamma$ close to $\xi$.  Just by continuity
 the
$X_t$-trajectories of $\xi$ and $z$ remain close, and by the
forward contraction along stable manifolds, the
$X_t$-trajectory of $\xi$ remains close to the segment
$[x_j,x_{j+1}]$.  Moreover the orbit of $z$ cannot
leave the tube through the union of the local strong stable
manifolds passing through $[x_j,x_{j+1}]$, for otherwise it
would contradict the definition of $Y_t$. Hence the
trajectory of $z$ must leave the tube through $\Sigma_{j+1}$. In
other words $X_{t(z)}(z)$ is a point of $\Sigma_{j+1}$,
close to $X_{t(\xi)}(\xi)$.

Let $\hat\gamma\subset\gamma$ be the largest connected
subset containing $\xi$ such that
$X_{t(z)}(z)\in\Sigma_{j+1}$ for all $z\in\hat\gamma$.  We
want to prove that $\hat\gamma=\gamma$ since this implies
that $R_j$ extends to the whole $\gamma_j$ and so, using
Lemma~\ref{stablereturnmap}, to the whole strip of $\Sigma_j$\,.

\begin{figure}[htb]
\begin{center}
\includegraphics[height=4cm]{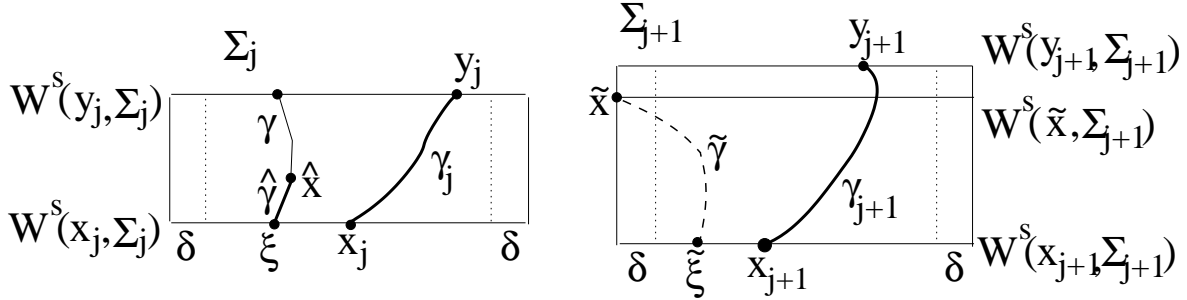}
\caption{\label{fig.exitingcross} Exiting the tube at $\Sigma_{j+1}$.}
\end{center}
\end{figure}

The proof is by contradiction. We assume $\hat\gamma$ is not
the whole $\gamma$\,, and let $\hat x$ be the endpoint
different from $\xi$\,. Then by definition of $\cF_\Sigma^s$
and of $Y_t$ (from Section~\ref{sec:every-orbit-leaves})
$\tilde x = X_{t(\hat x)}(\hat x)$ is on the center-unstable
boundary $\partial^{cu}\Sigma_{j+1}$ of the cross-section
$\Sigma_{j+1}$, between the stable manifolds
$W^s(x_{j+1},\Sigma_{j+1})$ and $W^s(y_{j+1},\Sigma_{j+1})$,
see Figure~\ref{fig.exitingcross}.
By the choice of $\gamma$ and by Corollary~\ref{ccone}, $
\tilde\gamma=\{X_{t(z)}(z): z\in\hat\gamma\} $ is a
$cu$-curve.  On the one hand, by
Lemma~\ref{l.lengthversusdistance}, the distance between
$\tilde x$ and $\tilde\xi=X_{t(\xi)}(\xi)$ dominates the
distance between their stable manifolds and
$\ell(\tilde\gamma)$
\[
\ell(\tilde\gamma)\le \kappa\cdot d(\xi,\tilde x)
\le \kappa\cdot d\big(W^s(x_{j+1},\Sigma), W^s(\tilde
x,\Sigma)\big).
\]
We note that $\ell(\tilde\gamma)$ is larger than $\delta$,
since $\xi$ is in $\Sigma_\delta$ and the section
$\Sigma_{j+1}$ is adapted.  On the other hand, the distance
between the two stable manifolds is smaller than the
distance between the stable manifold of $x_{j+1}$ and the
stable manifold of $y_{j+1}$\,, and this is smaller than
$K\cdot\delta_0$\,. Since $\delta_0$ is much smaller than
$\delta$, this is a contradiction.  This proves the claim
that $X_{t(z)}(z)\in\Sigma_{j+1}$ for all $z\in\gamma$.


\subsubsection{Conclusion of the proof of
  Theorem~\ref{t.expansivepoincare}}
\label{s.endoftheproof}

\begin{figure}[htpb]
  \centering
\psfrag{a}{$y_j$}
\psfrag{b}{$x_j$}
\psfrag{c}{$y_{j+1}$}
\psfrag{d}{$x_{j+1}$}
\psfrag{e}{$\gamma_{j+1}$}
\psfrag{f}{$\gamma_{j}$}
\psfrag{g}{$\Sigma^{j}$}
\psfrag{h}{$\Sigma^{j+1}$}
\includegraphics[height=3.5cm]{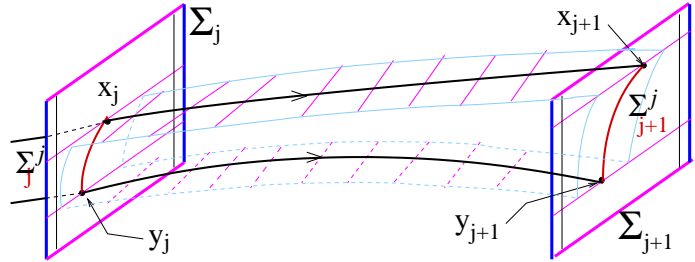}
  \caption{Expansion within the tube.}
  \label{fig:expans-within-tube}
\end{figure}

We have shown that there exists a well defined Poincar\'e
return map $R_j$ on the whole strip between the stable
manifolds of $x_j$ and $y_j$ inside $\Sigma$. By
Proposition~\ref{p.secaohiperbolica} and
Corollary~\ref{ccone} we know that the map $R_j$ is
hyperbolic where defined and, moreover, that the length of
each $cu$-curve is expanded by a factor of $3$ by $R_j$
(since we chose $\lambda=1/3$ in Section~\ref{s.22}).
Hence the distance between the stable manifolds
$R_j\big(W^s(x_j,\Sigma)\big)$ and
$R_j\big(W^s(y_j,\Sigma)\big)$ is increased by a factor
strictly larger than one, see
Figure~\ref{fig:expans-within-tube}. This contradicts item
(2) of Lemma~\ref{l.returny} since this distance will
eventually become larger than $K\cdot\delta_0$.  Thus $y_j$
must be in the stable manifold $W^s(x_j,\Sigma)$.  Since the
strong-stable manifold is locally flow-invariant and
$X_{h(\tau_j)}(y)$ is in the orbit of
$y_j=X_{\upsilon_j}(y)$, then $X_{h(\tau_j)}(y)\in
W^s(x_j)=W^s\big(X_{\tau_j}(x)\big)$, see
Lemma~\ref{l.returny}.

According to Lemma~\ref{l.returny} we have
$|\upsilon_j-h(\tau_j)|<K\cdot\delta_0$ and, by
Remark~\ref{r.trapaca}, there exits a small $\vep_1>0$ such
that
\[
R_\Sigma(y_j)=X_t(y_j)\in
W^{ss}_{\vep}(x_j)\quad\mbox{with}\quad
|\,t|<\vep_1.
\]
Therefore the piece of orbit
$\cO_y=X_{[\upsilon_j-K\cdot\delta_0-\vep_1, \upsilon_j +
  K\cdot\delta_0 +\vep_1]}(y)$ contains $X_{h(\tau_j)}(y)$. We
note that this holds for all sufficiently small values of
$\delta_0>0$ fixed from the beginning.

Now let $\vep_0>0$ be given and let us consider the piece of
orbit $\cO_x=X_{[\tau_j-\vep_0,\tau_j+\vep_0]}(x)$ and the
piece of orbit of $x$ whose strong-stable manifolds
intersect $\cO_y$, i.e.
\[
{\cO}_{xy}=\{ X_s(x) : 
\exists \tau\in[\upsilon_j-K\cdot\delta_0-t, \upsilon_j +
  K\cdot\delta_0 +t]\mbox{  such that  }
X_{\tau}(y)\in W^{ss}_{\vep}\big( X_s(x) \big) \}.
\]
Since $y_j\in W^s(x_j)$ we conclude that $\cO_{xy}$ is a
neighborhood of $x_j=X_{\tau_j}(x)$ which can be made as
small as we want taking $\delta_0$ and $\vep_1$ small
enough. In particular we can ensure that $\cO_{xy}\subset
\cO_x$ and so $X_{h(\tau_j)}(y)\in
W^{ss}_{\vep}\big(X_{[\tau_j-\vep_0,\tau_j+\vep_0]}(x)\big)$.
This finishes the proof of Theorem~\ref{t.expansivepoincare}.

 \section{Construction of physical measures}
 \label{sec:proof-theorem-B}
 
 Here we start the proof of Theorem~\ref{srb}.

\subsection{The starting point}
\label{sec:starting-point}

We show in Section~\ref{sec:global-poincare-map} that
choosing a \emph{global Poincar\'e section} $\Xi$ (with
several connected components) for $X$ on $\Lambda$, we can
reduce the transformation $R$ to the quotient over the
stable leaves. We can do this using
Lemma~\ref{stablereturnmap} with the exception of finitely
many leaves $\Gamma$, corresponding to the points whose
orbit falls into the local stable manifold of some
singularity or are sent into the stable boundary
$\partial^s\Sigma$ of some $\Sigma\in\Xi$ by $R$, where the
return time function $\tau$ is discontinuous.

As will be explained in Section~\ref{sec:reduct-quot-leaf},
the global Poincar\'e map $R:\Xi\to\Xi$ induces in this way
a map $f:\F\setminus\Gamma\to\F$ on the leaf space,
diffeomorphic to a finite union of open intervals $I$.
which is piecewise expanding and admits finitely many
$\upsilon_1,\dots,\upsilon_l$ ergodic absolutely continuous
(with respect to Lebesgue measure on $I$) invariant
probability measures (acim) whose basins cover Lebesgue
almost all points of $I$.

Moreover the Radon-Nikodym derivatives (densities)
$\frac{d\upsilon_k}{d\lambda}$ are 
\emph{bounded from above} and \emph{the support of each
  $\upsilon_k$ contains nonempty open intervals}, so the
basin $B(\upsilon_k)$ contains nonempty open intervals
Lebesgue modulo zero, $k=1,\dots,l$.

\subsection{Description of the construction}
\label{sec:steps-constr}

Afterwards we unwind the reductions made in
Section~\ref{sec:global-poincare-map} and obtain a physical
measure for the original flow at the end. 

We divide the construction of the physical measure for
$\Lambda$ in the following steps.

\begin{enumerate}
\item The compact metric space $\Xi$ is endowed with a
  partition $\cF$ and map $R:\Xi\setminus\Gamma\to\Xi$,
  where $\Gamma$ is a finite set of elements of $\cF$ (see
  Section~\ref{sec:finite-number-strips}).  The map $R$
  preserves the partition $\cF$ and contracts its elements
  by Lemma~\ref{stablereturnmap}. We have a finite family
  $\upsilon_1,\dots,\upsilon_l$ of absolutely continous
  invariant probability measures for the induced quotient
  map $f:\cF\setminus\Gamma\to\cF$.

  We show in Section~\ref{sec:reduct-quot-map} that each
  $\upsilon_i$ defines a $R$-invariant ergodic probability
  measure $\eta_i$. In Section~\ref{sec:phys-meas-glob} we
  show that the basin $B(\eta_i)$ is a union of strips of
  $\Xi$, and $\eta_i$ are therefore physical measures for
  $R$.  Moreover these basins cover $\Xi$:
  \begin{align*}
    \lambda^2\big(\Xi\setminus (B(\eta_1)\cup\dots\cup
    B(\eta_l))\big)=0,
  \end{align*}
  where $\lambda^2$ is the area measure on $\Xi$.
\item We then pass from $R$-invariant physical measures
  $\eta_1,\dots,\eta_l$ to invariant probability measures
  $\nu_1,\dots,\nu_l$ for the suspension semiflow over $R$
  with roof function $\tau$.  In the process we keep the
  ergodicity (Section~\ref{sec:phys-meas-susp}) and the
  basin property (Section~\ref{sec:phys-meas-susp}) of the
  measures: the whole space $\Xi\times[0,+\infty)/\sim$
  where the semiflow is defined equals the union of the
  ergodic basins of the $\nu_i$ Lebesgue modulo zero.
\item Finally in Section~\ref{sec:phys-meas-flow} we convert
  each physical measure $\nu_i$ for the semiflow into a
  physical measure $\mu_i$ for the original flow. We use
  that the semiflow is semiconjugated to $X_t$ on a
  neighborhood of $\Lambda$ by a local diffeomorphism.
  Uniqueness of the physical measure $\mu$ is then deduced
  in Section~\ref{sec:trans-uniq-phys} through the existence
  of a dense regular orbit in $\Lambda$ (recall that our
  definition of attractor \emph{demands} transitivity) and
  by the observation that the basin of $\mu$ contains open
  sets Lebesgue modulo zero. In Section~\ref{sec:hyphysmeas}
  we show that $\mu$ is (non-uniformly) hyperbolic.
\end{enumerate}

The details are exposed in the following sections.


\section{Global Poincar\'e maps and reduction to a
  one-dimensional map}
\label{sec:global-poincare-map}
Here we construct a global Poincar\'e map for the flow near
the singular-hyperbolic attractor $\Lambda$. We then use the
hyperbolicity properties of this map to reduce the dynamics
to a one-dimensional piecewise expanding map through a
quotient map over the stable leaves.


\subsection{Cross-sections and invariant foliations}
\label{sec:cross-sect-invar}

We observe first that by Lemma~\ref{l.existeadaptada} we can
take a $\de$-adapted cross-section at each non-singular
point $x\in\Lambda$. We know also that near each singularity
$\sigma_k$ there is a flow-box $U_{\sigma_k}$ as in
Section~\ref{s.24}, see Figure~\ref{f.singularbox}.

Using a tubular neighborhood construction near any given
adapted cross-section $\Sigma$, we linearize the flow in an
open set
$U_\Sigma=X_{(-\vep,\vep)}(\interior(\Sigma))$ for a
small $\vep>0$, containing the interior of the
cross-section.  This provides an open cover of the compact
set $\Lambda$ by flow-boxes near the singularities and
tubular neighborhoods around regular points.

We let $\{ U_{\Sigma_i}, U_{\sigma_k} : i=1,\dots, l; \,
k=1,\dots,s\}$ be a finite cover of $\Lambda$, where $s\ge1$
is the number of singularities in $\Lambda$, and we set
$t_3>0$ to be an upper bound for the time each point $z\in
U_{\Sigma_i}$ takes to leave the tubular neighborhood by the
action of the flow, for any $i=1,\dots,l$. We assume without
loss of generality that $t_2>t_3$.
 
To define the Poincar\'e map $R$, for any point $z$ in one
of the cross-sections in
 \[
\Xi=
\{\Sigma_j,\Sigma^{i,\pm}_{\sigma_k},\Sigma^{o,\pm}_{\sigma_k}: 
  j=1,\dots,l; k=1,\dots, s\},
\]
we consider $\hat z=X_{t_2}(z)$ and wait for the next time
$t(z)$ the orbit of $\hat z$ hits again one of the
cross-sections. Then we define $R(z)=X_{t_2+t(z)}(z)$ and
say that $\tau(z)=t_2+t(z)$ is the \emph{Poincar\'e time} of
$z$.  If the point $z$ never returns to one of the
cross-sections, then the map $R$ is not defined at $z$ (e.g.
at the lines $\ell^\pm$ in the flow-boxes near a
singularity).  Moreover by Lemma~\ref{stablereturnmap}, if
$R$ is defined for $x\in\Sigma$ on some $\Sigma\in\Xi$, then
$R$ is defined for every point in $W^s(x,\Sigma)$. Hence
\emph{the domain of $R\mid\Sigma$ consists of strips of
  $\Sigma$}. The smoothness of $(t,x)\mapsto X_t(x)$ ensures
that the strips
\begin{equation}
  \label{eq:strip}
  \Sigma(\Sigma') = \{ x\in\Sigma : R(x)\in\Sigma'\}
\end{equation}
have non-empty interior in $\Sigma$ for every
$\Sigma,\Sigma'\in\Xi$.
When $R$ maps to an outgoing strip near a singularity
$\sigma_k$, there might be a boundary of the strip
corresponding to the line $\ell^\pm_k$ of points which fall
in the stable manifold of $\sigma_k$. 


\begin{remark}
  \label{rmk:firstreturn}
  Consider the Poincar\'e map given by the \emph{first
    return map} $R_0:\Xi\to\Xi$ defined simply as
  $R_0(z)=X_{T(z)}(z)$, where
  \[
  T(z)=\inf\{t>0: X_t(z)\in\Xi\}
  \]
  is the time the $X$-orbit of $z\in\Xi$ takes to arrive
  again at $\Xi$.  This map $R_0$ is not defined on those
  points $z$ which do not return and, moreover, $R_0$ might
  not satisfy the lemmas of Section~\ref{s.22}, since we do
  not know whether the flow from $z$ to $R_0(z)$ has enough
  time to gain expansion. However the stable manifolds are
  still well defined.  By the definitions of $R_0$ and of
  $R$ we see that \emph{$R$ is induced by $R_0$}, i.e.
  \emph{if $R$ is defined for $z\in\Xi$, then there exists
    an integer $r(x)$ such that}
  \[
  R(z)=R_0^{r(z)}(z).
  \]
  We note that since the number of cross-sections in $\Xi$
  is finite and the time $t_2$ is a constant,
  then the function $r:\Xi\to\nat$ is bounded: there exists
  $r_0\in\nat$ such that $r(x)\le r_0$ for all $x\in\Xi$.
\end{remark}

\subsubsection{Finite number of strips}
\label{sec:finite-number-strips}

We show that fixing a cross-section $\Sigma\in\Xi$ the
family of all possible strips as in~\eqref{eq:strip} covers
$\Sigma$ except for finitely many stable leaves
$W^s(x_i,\Sigma), i=1,\dots, m=m(\Sigma)$.
Moreover we also show that each strip given
by~\eqref{eq:strip} has finitely many connected components.
Thus the number of strips in each cross-section is finite.

We first recall that each $\Sigma\in\Xi$ is contained in
$U_0$, so $x\in\Sigma$ is such that
$\omega(x)\subset\Lambda$. Note that $R$ is locally smooth
for all points $x\in\interior(\Sigma)$ such that
$R(x)\in\inte(\Xi)$ by the flow box theorem and the
smoothness of the flow, where $\inte(\Xi)$ is the union of
the interiors of each cross-section of $\Xi$.  Let
$\partial^s\Xi$ denote the union of all the leaves forming
the stable boundary of every cross-section in $\Xi$.

\begin{lemma}
  \label{le:descont}
  The set of discontinuities of $R$ in
  $\Xi\setminus\partial^s\Xi$ is contained in the set of
  points $x\in\Xi\setminus\partial^s\Xi$ such that:
\begin{enumerate}
\item either $R(x)$ is defined and belongs to $\partial^s\Xi$;
\item or there is some time $0<t \le t_2$ such that $X_t(x)\in
  W^s_{loc}(\sigma)$ for some singularity $\sigma$ of
  $\Lambda$.
\end{enumerate}
Moreover this set is contained in a finite number of stable
leaves of the cross-sections $\Sigma\in\Xi$.
\end{lemma}


\begin{proof} We divide the proof into several steps.

\begin{description}
\item[Step 1] Cases (1) and (2) in the statement of the
  lemma correspond to all possible discontinuities of $R$ in
  $\Xi\setminus\partial^s\Xi$.
\end{description}

Let $x$ be a point in $\Sigma\setminus\partial^s\Sigma$ for
some $\Sigma\in\Xi$, \emph{not} satisfying any of the
conditions in items (1) and (2).  Then $R(x)$ is defined and
$R(x)$ belongs to the interior of some cross-section
$\Sigma'$. By the smoothness of the flow and by the flow box
theorem we have that $R$ is smooth in a neighborhood of $x$
in $\Sigma$. Hence any discontinuity point for $R$ must be
in one the situations (1) or (2).

  \begin{description}
  \item[Step 2] Points satisfying item (2) are contained in
    finitely many stable leaves in each $\Sigma\in\Xi$.
  \end{description}
  
  Indeed if we set $W=X_{[-t_2,0]}\big( \cup_\sigma
  W^s_{loc}(\sigma) \big)$, where the union above is taken
  over all singularities $\sigma$ of $\Lambda$, then $W$ is
  a compact sub-manifold of $M$ with boundary, tangent to the
  center-stable sub-bundle $E^s\oplus E^X$. This means that
  $W$ is transversal to any cross-section of $\Xi$.
  
  Hence the intersection of $W$ with any $\Sigma\in\Xi$ is a
  one-dimensional sub-manifold of $\Sigma$. Thus the number
  of connected components of the intersection is finite in
  each $\Sigma$. This means that there are finitely many
  points $x_1,\dots, x_k\in\Sigma$ such that
\[
W\cap\Sigma\subset W^s(x_1,\Sigma)\cup\dots \cup W^s(x_k,\Sigma).
\]


  \begin{description}
  \item[Step 3] Points satisfying item (1) are contained in
    a finite number of stable leaves of each $\Sigma\in\Xi$.
  \end{description}

  We argue by contradiction. Assume that the set of points
  $D$ of $\Sigma$ sent by $R$ into stable boundary points of
  some cross-section of $\Xi$ is such that
  \[
  L=\{ W^s(x,\Sigma) : x\in D \}
  \]
  has \emph{infinitely many lines}. Note that $D$ in fact
  equals $L$ by Lemma~\ref{stablereturnmap}. Then there
  exists an accumulation line $W^s(x_0,\Sigma)$. Since the
  number of cross-sections in $\Xi$ is finite we may assume
  that $W^s(x_0,\Sigma)$ is accumulated by \emph{distinct}
  $W^s(x_i,\Sigma)$ with $x_i\in D$ satisfying $R(x_i)\in
  W^s(z,\Sigma')\subset \partial^s\Sigma'$ for a fixed
  $\Sigma'\in\Xi$, $i\ge1$. We may assume that $x_i$ tends
  to $x_0$ when $i\to\infty$, that $x_0$ is in the
  interior of $W^s(x_0,\Sigma)$ and that the $x_i$ are all
  distinct --- in particular the points $x_i$ do not belong
  to any periodic orbit of the flow since we can choose
  the $x_i$ anywhere in the stable set $W^s(x_i,\Sigma)$.
  
  As a preliminary result we show that $R(x_i)=X_{s_i}(x_i)$
  is such that $s_i$ is a bounded sequence in the real line.
  For otherwise $s_i\to\infty$ and this means, by
  definition of $R$, that the orbit of $X_{t_2}(x_i)$ is
  very close to the local stable manifold of some
  singularity $\sigma$ of $\Lambda$ and that $R(x_i)$
  belongs to the outgoing cross-section near this
  singularity: $R(x_i)\in\Sigma_\sigma^{o,\pm}$.  Hence we
  must have that $X_{s_i}(x_i)$ tends to the stable manifold
  of $\sigma$ when $i\to\infty$ and that $R(x_i)$ tends to
  the stable boundary of $\Sigma_\sigma^{o,\pm}$.
  Since no point in any cross-section in $\Xi$ is sent by
  $R$ into this boundary line, we get a contradiction.

  Now the smoothness of the flow and the fact that
  $W^s(z,\Sigma')$ is closed imply that $R(x_0)\in
  W^s(z,\Sigma')$ also since we have the following
  \[
  R(x_0)=\lim_{i\to\infty} R(x_i) = \lim_{i\to\infty}
  X_{s_i}(x_i)=X_{s_0}(x_0)\quad\mbox{and}\quad
  \lim_{i\to\infty} s_i = s_0.
  \]
  Moreover $R(W^s(x_0,\Sigma))\subset W^s(z,\Sigma')$ and
  $R(x_0)$ is in the interior of $R(W^s(x_0,\Sigma))$, then
  $R(x_i)\in R(W^s(x_0,\Sigma))$ for all $i$ big enough. 
  This means that there exists a sequence $y_i\in
  W^s(x_0,\Sigma)$ and a sequence of real numbers $\tau_i$
  such that $X_{\tau_i}(y_i)=R(y_i)=R(x_i)$ for all sufficiently big
  integers $i$. By construction we have that $x_i\neq y_i$
  and both belong to the same orbit. Since $x_i,y_i$ are in
  the same cross-section we get that
  $x_i=X_{\alpha_i}(y_i)$ with $|\alpha_i|\ge t_3$ for
  all big $i$.
  
  However we also have that $\tau_i\to s_0$ because
  $R(y_i)=R(x_i)\to R(x_0)$, $y_i\in W^s(x_0,\Sigma)$ and
  $R\mid W^s(x_0,\Sigma)$ is smooth.  Thus $|s_i-\tau_i|\to
  0$. But $|s_i-\tau_i|=|\alpha_i|\ge t_3>0$. This is a
  contradiction.

This proves that $D$ is contained in finitely many stable
leaves.

Combining the three steps above we conclude the proof of the
lemma.
\end{proof}

Let $\Gamma$ be the finite set of stable leaves of $\Xi$
provided by Lemma~\ref{le:descont} together with
$\partial^s\Xi$. Then the complement $\Xi\setminus\Gamma$ of
this set is formed by finitely many open strips where $R$ is
smooth. Each of these strips is then a connected component
of the sets $\Sigma(\Sigma')$ for $\Sigma,\Sigma'\in\Xi$.

\subsubsection{Integrability of the global Poincar\'e return
  time}
\label{sec:integr-glob-poinc}

We claim that \emph{the Poincar\'e time $\tau$ is integrable
  with respect to the Lebesgue area measure on $\Xi$}.
Indeed given $z\in\Xi$, the point $\hat z=X_{t_2}(z)$ either
is inside a flow-box $U_{\sigma_k}$ of a singularity
$\sigma_k$, or not. In the former case, the time $\hat z$
takes to reach an outgoing cross-section
$\Sigma^{o,\pm}_{\sigma_k}$ is bounded by the exit time
function $\tau^{\pm}_{\sigma_k}$ of the corresponding
flow-box, which is integrable, see Section~\ref{s.24}. In
the latter case, $\hat z$ takes a time of at most $2\cdot
t_3$ to reach another cross-section, by definition of $t_3$.
Thus the Poincar\'e time on $\Xi$ is bounded by $t_2+2\cdot
t_3$ plus a sum of finitely many integrable functions, one
for each flow-box near a singularity, by finiteness of the
number of singularities, of the number of cross-sections in
$\Xi$ and of the number of strips at each cross-section.
This proves the claim.

\begin{remark}
  \label{rmk:imageinside}
Given $z\in\Sigma\in\Xi$ we write
$\tau^k(z)=\tau(R^{k-1}(z))+\dots+\tau(z)$ for $k\ge1$ and
so $\tau=\tau^1$. Since
\[
R^k\big(W^s(z,\Sigma)\big)\subset
X_{\tau^k(z)}\big(W^s(z,\Sigma)\big)\subset
X_{\tau^k(z)}(U),
\]
the length $\ell\Big(R^k\big(W^s(z,\Sigma)\big)\Big)$ is
uniformly contracted and $\tau^k(z)\to+\infty$ when
$k\to+\infty$, we get that
$R^k\big(W^s(z,\Sigma)\big)\subset\Sigma^{\,\prime}$ for
some $\Sigma^{\,\prime}\in\Xi$ and
\[
d\Big( R^k\big(W^s(z,\Sigma)\big) ,
\partial^{cu}\Sigma^{\,\prime} \Big) > \delta/2
\]
for all big enough $k$, by the definition of $U$ and of
$\delta$-adapted cross-section. (The distance $d(A,B)$
between two sets $A,B$ means $\inf\{d(a,b):a\in A, b\in
B\}$.) We may assume that this property holds for all stable
leaves $W^s(z,\Sigma)$, all $z\in\Sigma$ and every
$\Sigma\in\Xi$ for all $k\ge k_0$, for some fixed big
$k_0\in\nat$, by the uniform contraction property of $R$ in
the stable direction.
\end{remark}

\subsection{Absolute continuity of foliations}
\label{sec:absolute-contin-foli}

From now on we assume that the flow $(X_t)_{t\in\real}$ is
$C^2$.  Under this condition it is well known
\cite{PT93,Man87} that the stable leaf $W^s(x,\Sigma)$ for
every $x\in \Sigma\in\Xi$ is a $C^2$ embedded disk and these
leaves define a $C^{1+\alpha}$ foliation
$\F^{\,s}_{\Sigma}$, $\alpha\in(0,1)$, of each
$\Sigma\in\Xi$, as we now explain.

\begin{figure}[htpb]
  \centering
  \psfrag{X}{$X_t(x)$}\psfrag{g}{$\gamma_1$}
  \psfrag{j}{$\gamma_2$}\psfrag{z}{$\xi_0$}
  \psfrag{S}{$S$}\psfrag{E}{$D_1$}\psfrag{D}{$D_0$}
  \psfrag{u}{$W^u$\text{ leaves}}
  \includegraphics[height=4cm,width=7.5cm]{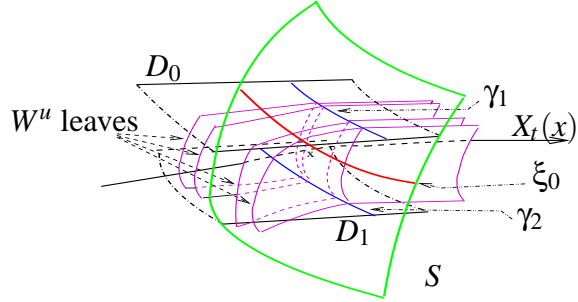}
  \caption{The holonomy map.}
  \label{fig:holonomy-maps}
\end{figure}

Recall the setting presented before the statement of
Theorem~\ref{thm:srbmesmo} to explain the disintegration
along center-unstable manifolds. Let $x\in\Lambda$ and $S$
be a cross-section to the flow at $x$ and $\xi_0$ be the
connected component of $W^s(x)\cap S$ containing $x$.
Assume that $x$ has a unstable leaf $W^u(x)$ and let
$D_1,D_2$ be embedded disk in $M$ transverse to $W^u(x)$ at
$x_1,x_2$, that is $T_{x_i}D_i\oplus T_{x_i}W^u(x)=T_{x_i}
M$, $i=1,2$. Then the strong-unstable leaves through the
points of $D_1$ which cross $D_2$ define a map $h$ between a
subset of $D_1$ to $D_2$: $h(y_1)=y_2=W^{uu}(y_1)\cap D_2$,
called the \emph{holonomy} map of the strong-unstable
foliation between the transverse disks $D_1,D_2$. The
holonomy is injective if $D_1,D_2$ are close enough due to
uniqueness of the strong-unstable leaves through $\mu$-a.e.
point, see Figure~\ref{fig:holonomy-maps}.

We say that $h$ is \emph{absolutely continuous} if there is
a measurable map $J_h:D_1\to[0,+\infty]$, called the
\emph{Jacobian of h}, such that
\begin{align*}
  \Leb_2\big( h(A) \big) = \int_A J_h\, d\Leb_1
  \text{  for all Borel sets  } A\subset D_1,
\end{align*}
and $J_h$ is integrable with respect to $\Leb_1$ on $D_1$,
where $\Leb_i$ denotes the Lebesgue measure induced on $D_i$
by the Riemannian metric, $i=1,2$.

The foliation $\{W^{uu}(x)\}$ is \emph{absolutely
  continuous} (H\"older continuous) if every holonomy map is
absolutely continuous (or $J_h$ is H\"older continuous,
respectively).

Since the pioneering work of Anosov and Sinai
\cite{An67,AS67} it became clear that for $C^{2}$
transformations or flows (in fact it is enough to have
transformations or flows which are $C^1$ with
$\alpha$-H\"older derivative for some $0<\alpha<1$) the
strong-unstable foliation is absolutely continuous and
H\"older continuous. See also \cite{Man87} and~\cite{PS89}
for detailed presentations of these results.
\emph{When the leaves are of codimension one, then the
  Jacobian $J_h$ of the holonomy map $h$ coincides with the
  derivative $h'$ since $h$ is a map between curves in $M$.}
In this case the holonomy map can be seen as a
$C^{1+\alpha}$ transformation between subsets of the real
line. A dual statement is also true for the strong-stable
foliations and corresponding holonomies.

In the case of the stable foliation for a flow, we have that
for any pair of disks $\gamma_1,\gamma_2$ inside $S$
transverse to $W^s(x)\cap S$ at distinct points $y_1,y_2$,
the holonomy $H$ between $\gamma_1$ and $\gamma_2$ along the
leaves $W^{s}(z)\cap S$ crossing $S$ is also H\"older
continuous \emph{if the flow is $C^2$}.

Indeed note that this holonomy map $H$ can be obtained as a
composition of the holonomy map $h$ between two disks
$D_1,D_2$ transverse to the strong-stable leaves which
cross $S$, and the ``projection along the flow'' sending
$w\in X_{(-\delta,\delta)}(S)$ to a point $X_t(w)\in S$
uniquely defined, with $t\in(-\delta,\delta)$. The disks are
defined simply as $D_i=X_{(-\epsilon,\epsilon)}(\gamma_i)$ for
$0<\epsilon<\delta$ and satisfy $D_i\cap S_0=\gamma_1$,
$i=1,2$. Since the the holonomy $h$ is H\"older continuous
and the projection along the flow has the same
differentiability class of the flow (due to the Tubular Flow
Theorem, see e.g.~\cite{PM82}), we see that the holonomy $H$
is also H\"older continuous.

\subsection{Reduction to the quotient leaf space}
\label{sec:reduct-quot-leaf}

We choose once and for all a $C^2$ $cu$-curve
$\gamma_{\,\Sigma}$ transversal to $\F^{\,s}_{\Sigma}$ in
each $\Sigma\in\Xi$.  Then by the discussion in the previous
Section~\ref{sec:absolute-contin-foli} the projection
$p_\Sigma$ along leaves of $\F^{\,s}_{\Sigma}$ onto
$\gamma_{\,\Sigma}$ is a $C^{1+\alpha}$ map. We set
\[
I=\bigcup_{\Sigma,\Sigma'\in\Xi}
\interior\big( \Sigma( \Sigma' ) \big)\cap\gamma_{\,\Sigma}
\]
and observe that by the properties of $\Sigma( \Sigma' )$
obtained in Section~\ref{sec:cross-sect-invar} the set
$I$ is diffeomorphic to a finite union of non-degenerate
open intervals $I_1,\dots,I_m$ by a $C^2$ diffeomorphism and
$p_\Sigma\mid p_\Sigma^{-1}(I)$ becomes a $C^{1+\alpha}$
submersion. Note that since $\Xi$ is finite we can choose
$\gamma_\Sigma$ so that $p_\Sigma$ has bounded derivative:
\begin{align*}
  \text{there exists   } \beta_0>1 \text{  such that  }
  \frac1{\beta_0}\le \big| Dp_\Sigma\mid\gamma \big|
  \le\beta_0
  \text{  \emph{for every $cu$-curve}  } \gamma
  \text{  inside any  } \Sigma\in\Xi.
\end{align*}
In particular, denoting the Lebesgue area
measure over $\Xi$ by $\lambda^2$ and the Lebesgue length
measure on $I$ by $\lambda$, we have
$(p_\Sigma)_*\lambda^2\ll\lambda$.

According to Lemma~\ref{stablereturnmap},
Proposition~\ref{p.secaohiperbolica} and
Corollary~\ref{ccone} the Poincar\'e map $R:\Xi\to\Xi$ takes
stable leaves of $\F^{\,s}_{\Sigma}$ inside stable leaves of
the same foliation and is hyperbolic. In addition a
$cu$-curve $\gamma\subset\Sigma$ is taken by $R$ into a
$cu$-curve $R(\gamma)$ in the image cross-section.  Hence
the map
\[
f:I\to I
\quad\mbox{given by}\quad
I\ni z \mapsto
p_{\Sigma'}\Big( R \big( W^s(z,\Sigma) \cap \Sigma(\Sigma')
\big) \Big)
\]
for $\Sigma,\Sigma'\in\Xi$ is a $C^{1+\alpha}$ map and for
points in the interior of $I_i$, $i=1,\dots,m$
\begin{align}\label{eq:DpSigma}
  \big| Df | = \big| D\big( p_{\Sigma'}\circ R \circ
  \gamma_\Sigma\big) \big| \ge \frac1{\beta_0}\cdot \sigma.
\end{align}
Thus choosing $t_1$ (and consequently $t_2$) big enough so
that $\sigma/\beta_0>3/2>1$ in
Proposition~\ref{p.secaohiperbolica}, we obtain that $f$ is
piecewise expanding.  Moreover $|f'|^{-1}\mid I_j$ is a
$\alpha$-H\"older function since for all $x,y\in I_j$ we
have
\[
\frac1{|f'(x)|}-\frac1{|f'(y)|} \le \frac{|f'(x)-f'(y)|}{|f'(x)
  f'(y)|}
\le \frac{C}{(3/2)^2}\cdot |x-y|^\alpha,
\quad\mbox{for some}\quad 0<\alpha<1.
\]
Thus $f:I\to I$ is a \emph{$C^{1+\alpha}$ piecewise
  expanding map}.

\begin{remark}
  \label{rmk:Poincareleaves}
  By Lemma~\ref{stablereturnmap} the Poincar\'e time $\tau$
  is constant on stable leaves $W^s(x,\Sigma)$ for all
  $x\in\Sigma\in\Xi$. Thus after
  Section~\ref{sec:integr-glob-poinc} there exists a return
  time function $\tau_I$ on $I$ such that $\tau=\tau_I\circ
  p$, where $p:\Xi\to \gamma_\Xi$ is the joining of all
  $p_\Sigma$, $\Sigma\in\Xi$ and
  $\gamma_\Xi=\{\gamma_\Sigma:\Sigma\in\Xi\}$.  The
  integrability of $\tau$ with respect to $\lambda^2$ (see
  Section~\ref{sec:cross-sect-invar}) implies the
  $\lambda$-integrability of $\tau_I$ naturally since
  $(p_\Sigma)_*\lambda^2\ll\lambda$ and $\tau_I\circ
  p=\tau$.
\end{remark}

\subsection{Existence and finiteness of acim's}
\label{sec:existence-finiten-ac}

It is well known \cite{Vi97b,Wo78,HK82} that $C^1$ piecewise
expanding maps $f$ of the interval such that $1/|f'|$ is of
bounded variation, have finitely many absolutely continuous
invariant probability measures whose basins cover Lebesgue
almost all points of $I$.

Using an extension of the notion of bounded variation
(defined below) it was shown in \cite{Ke85} that the results
of existence and finiteness of absolutely continuous ergodic
invariant measures can be extended to $C^1$ piecewise
expanding maps $f$ such that $g=1/|f'|$ is $\alpha$-H\"older
for some $\alpha\in(0,1)$. These functions are of
universally bounded variation, i.e.
\[
\sup_{a=a_0<a_1<\dots<a_n=b} \left(
\sum_{j=1}^n \big|
\varphi(a_i)-\varphi(a_{i-1})
\big|^{1/\alpha}
\right)^\alpha <\infty,
\]
where the supremum is taken over all finite partition of the
interval $I=[a,b]$. 
Moreover from \cite[Theorem 3.2]{Ke85} the densities $\varphi$ of
the absolutely continuous invariant probability measures for
$f$ satisfy the following: there exists constants $A,C>0$
such that
\begin{align*}
  \int\osc(\varphi,\epsilon,x)\,dx \le
  C\cdot\epsilon^\alpha
  \quad\text{for all}\quad 0<\epsilon\le A,
\end{align*}
where $\osc(\varphi,\epsilon,x)=\esup_{y,z\in B(x,\epsilon)}
\big|\varphi(y)-\varphi(z)\big|$ and the essential supremo
is taken with respect to Lebesgue measure. From this we can
find a sequence $\epsilon_n\to0$ such that
$\osc(\varphi,\epsilon_n,\cdot)\xrightarrow[n\to\infty]{}0$
(with respect to Lebesgue measure). This implies that
$\supp(\varphi)$ \emph{contains non-empty open intervals}.

Indeed, for a given small $\delta>0$ let $\alpha>0$ be so
small and $n$ so big that 
\begin{align*}
  W=\{\varphi>\alpha\}\text{  and  }
  V=\{\osc(\varphi,\epsilon_n,\cdot)>\alpha/2\}
  \quad\text{satisfy}\quad
  \lambda(I\setminus W)<\delta \text{   and  }
  \lambda(V)<\delta.
\end{align*}
Then $\lambda(W\cap I\setminus V)>1-2\delta>0$. Let $x$ be a
Lebesgue density point of $W\cap I\setminus V$. Then there
exists a positive Lebesgue measure subset of
$B(x,\epsilon_n)$ where $\varphi>\alpha$. By definition of
$\osc(\varphi,\epsilon_n,x)$ this implies that for Lebesgue almost
every $y\in B(x,\epsilon_n)$ we have $\varphi(y)>\alpha/2>0$,
thus $B(x,\epsilon_n)\subset\supp(\varphi)$.

In addition from \cite[Theorem 3.3]{Ke85} there are finitely
many ergodic absolutely continuous invariant probability
measures $\upsilon_1,\dots,\upsilon_l$ of $f$ and every
absolutely continuous invariant probability measure
$\upsilon$ decomposes into a convex linear combination
$\upsilon=\sum_{i=1}^l a_i\upsilon_i$. From \cite[Theorem
3.2]{Ke85} considering any subinterval $J\subset I$ and the
normalized Lebesgue measure $\lambda_J=(\lambda\mid J)/\lambda(J)$ on
$J$, then every weak$^*$ accumulation point of
$n^{-1}\sum_{j=0}^{n-1} f_*^j(\lambda_J)$ is an absolutely
continuous invariant probability measure $\upsilon$ for $f$
(since the indicator function of $J$ is of generalized
$1/\alpha$-bounded variation). Hence the basin of the
$\upsilon_1,\dots,\upsilon_l$ cover $I$ Lebesgue modulo
zero: $
\lambda\big(I\setminus(B(\upsilon_1)\cup\dots\cup B(\upsilon_l)\big)
=0. $

Note that from \cite[Lemma 1.4]{Ke85} we also know that
\emph{the density $\varphi$ of any absolutely continous
  $f$-invariant probability measure} \emph{is bounded from
  above}.  In what follows we show how to use these
properties to build physical measures for the flow.


\section{Physical measures through suspension}
\label{s.fsuspending}

Here we show, in Section~\ref{sec:reduct-quot-map}, how to
construct an invariant measure for a transformation from an
invariant measure for the quotient map obtained from a
partition of the space. We show also that if the measure is
ergodic on the quotient, then we also obtain ergodicity on
the starting space. In Section~\ref{sec:phys-meas-glob} we
apply these results to the global Poincar\'e map $R$ of a
singular-hyperbolic attractor and its corresponding
one-dimensional quotient map $f$.

In Section~\ref{sec:reduct-poinc-map} we extend the
transformation to a semi-flow through a suspension
construction and show that each invariant and ergodic
measure for the transformation corresponds to a unique
measure for the semi-flow with the same properties. In
Section~\ref{sec:phys-meas-susp} we again apply these
results to the transformation $R$ to obtain physical
measures for the suspension semiflow over $R$ with roof
function $\tau$.


\subsection{Reduction to the quotient map}
\label{sec:reduct-quot-map}
Let $\Xi$ be a compact metric space, $\Gamma\subset \Xi$
and $F:(\Xi\setminus\Gamma)\to\Xi$ be a measurable map.
We assume that there exists a partition $\F$ of $\Xi$ into
measurable subsets, having $\Gamma$ as an element, which is
\begin{itemize}
\item{{\em invariant:\/}} the image of any $\xi\in\F$ distinct
     from $\Gamma$ is contained in some element $\eta$ of $\F$;
\item{{\em contracting:\/}} the diameter of $F^n(\xi)$ goes to
     zero when $n\to\infty$, uniformly over all the $\xi\in\F$
     for which $F^n(\xi)$ is defined.
\end{itemize}
We denote $p:\Xi\to \F$ the canonical projection, i.e. $p$
assigns to each point $x\in\Xi$ the atom $\xi\in\F$ that
contains it.
By definition, $A\subset \F$ is measurable if and only if
$p^{-1}(A)$ is a measurable subset of $\Xi$ and likewise $A$
is open if, and only if, $p_\Sigma^{-1}(A)$ is open in $\Xi$.
The invariance condition means that there is a uniquely
defined map
$$
f:(\F\setminus\{\Gamma\}) \to \F
\quad\text{such that}\quad
f\circ p = p \circ F.
$$
Clearly, $f$ is measurable with respect to the measurable
structure we introduced in $\F$.  We assume from now on that
the leaves are sufficiently regular so that $\Xi/\F$ is a
metric space with the topology induced by $p$.

Let $\mu_f$ be any probability measure on $\F$ invariant
under the transformation $f$.
For any bounded function $\psi:\Xi\to\real$, let
$\psi_{-}:\F\to\real$ and $\psi_{+}: \F\to\real$ be
defined by
$$
\psi_{-}(\xi)=\inf_{x\in\xi}\psi(x)
\qquad\mbox{and}\qquad
\psi_{+}(\xi)=\sup_{x\in\xi}\psi(x).
$$

\begin{lemma}
\label{l.fconvergence}
Given any continuous function $\psi:\Xi\to\real$,
both limits
\begin{equation}
  \label{eq:bothlim}
\lim_n \int (\psi\circ F^n)_{-}\,d\mu_f
\quad\text{and}\quad
\lim_n \int (\psi\circ F^n)_{+}\,d\mu_f  
\end{equation}
exist, and they coincide.
\end{lemma}

\begin{proof}
Let $\psi$ be fixed as in the statement.
Given $\vep>0$, let $\delta>0$ be such that
$|\psi(x_1)-\psi(x_2)|\le\vep$ for all $x_1, x_2$ with
$d(x_1,x_2)\le \delta$.
Since the partition $\F$ is assumed to be contractive,
there exists $n_0\ge 0$ such that $\diam(F^n(\xi))\le\delta$
for every $\xi\in\F$ and any $n\ge n_0$.
Let $n + k \ge n \ge n_0$.   
By definition,
$$
(\psi\circ F^{n+k})_{-}(\xi) - (\psi\circ F^{n})_{-}(f^k(\xi))
= \inf (\psi \mid F^{n+k}(\xi) )
- \inf (\psi \mid F^{n}(f^k(\xi)) ).
$$
Observe that $F^{n+k}(\xi) \subset F^n(f^k(\xi))$.
So the difference on the right hand side is bounded by
$$
\sup \big( \psi \mid F^{n}(f^k(\xi))\big)
- \inf \big(\psi \mid F^{n}(f^k(\xi))\big)
\le \vep.
$$
Therefore
$$
\left|\int (\psi \circ F^{n+k})_{-} \,d\mu_f 
 - \int (\psi \circ F^{n})_{-}\circ f^k \,d\mu_f\right|
\le \vep.
$$
Moreover, one may replace the second integral by
$\int (\psi \circ F^{n})_{-} \, d\mu_f$, because
$\mu_f$ is $f$-invariant.

At this point we have shown that $\{\int (\psi\circ
F^n)_{-}\,d\mu_F\}_{n\ge 1}$ is a Cauchy sequence in
$\real$.  In particular, it converges.  The same argument
proves that $\big\{\int (\psi\circ
F^n)_{+}\,d\mu_F\big\}_{n\ge1}$ is also convergent.
Moreover, keeping the previous notations,
$$
0\le (\psi\circ F^n)_{+}(\xi) - (\psi\circ F^n)_{-}(\xi)
= \sup \big(\psi \mid F^n(\xi)\big)
- \inf \big(\psi \mid F^n(\xi)\big)
\le \vep
$$
for every $n\ge n_0$.  So the two sequences
in~\eqref{eq:bothlim} must have the same limit.
The lemma is proved.
\end{proof}

\begin{corollary}
\label{c.poincareinvariant}
There exists a unique probability measure $\mu_F$ on $\Xi$
such that 
$$
\int \psi\,d\mu_F
=\lim \int (\psi\circ F^n)_{-}\,d\mu_f
=\lim \int (\psi\circ F^n)_{+}\,d\mu_f.
$$
for every continuous function $\psi:\Xi\to\real$.
Besides, $\mu_F$ is invariant under $F$. Moreover the
correspondence $\mu_f\mapsto\mu_F$ is injective.
\end{corollary}

\begin{proof}
Let $\hat\mu(\psi)$ denote the value of the two limits.
Using the expression for $\hat\mu(\psi)$ in terms of
$(\psi\circ F^n)_{-}$ we immediately get that
$$
\hat\mu(\psi_1+\psi_2)
\ge \hat\mu(\psi_1) + \hat\mu(\psi_2).
$$
Analogously, the expression of $\hat\mu(\psi)$ in terms
of $(\psi\circ F^n)_{+}$ gives the opposite inequality.
So, the function $\hat\mu(\cdot)$ is additive. 
Moreover, $\hat\mu(c\psi) = c\hat\mu(\psi)$ for every
$c\in\real$ and every continuous function $\psi$. 
Therefore, $\hat\mu(\cdot)$ is a linear real operator
in the space of continuous functions $\psi:\Xi\to\real$.

Clearly, $\hat\mu(1)=1$ and the operator $\hat\mu$ is
non-negative: $\hat\mu(\psi) \ge 0$ if $\psi\ge 0$.
By the Riesz-Markov theorem, there exists a unique measure
$\mu_F$ on $\Xi$ such that
$\hat\mu(\psi)=\int\psi\,d\mu_F$
for every continuous $\psi$.
To conclude that $\mu_F$ is invariant under $F$ it suffices
to note that
$$
\hat\mu(\psi\circ F)
=\lim_n \int(\psi\circ F^{n+1})_{-}\,d\mu_f
= \hat\mu(\psi)
$$
for every $\psi$. 

To prove that the map $\mu_f\mapsto\mu_F$
is injective, we note that if
$\mu_F=\mu_F'$ are obtained from $\mu_f$ and $\mu_f'$
respectively, then for any continuous function
$\vr:\F\to\real$ we have that $\psi=\vr\circ p:\Xi\to\real$
is continuous. But
\[
\mu_f\big( (\psi\circ F^n)_\pm \big) =
\mu_f\big( (\vr\circ p\circ F^n)_\pm \big)=
\mu_f\big( (\vr\circ f^n\circ p)_\pm \big)=
\mu_f(\vr\circ f^n)=\mu_f(\vr)
\]
for all $n\ge1$ by the $f$-invariance of $\mu_f$. Hence by
definition
\[
\mu_f(\vr)=\mu_F(\psi)=\mu_F'(\psi)=\mu_f'(\vr)
\]
and so $\mu_f=\mu_f'$.
This finishes the proof of the corollary.
\end{proof}


\begin{lemma}
\label{l.SRBmSRBp}
Let $\psi:\Xi\to\real$ be a continuous function
and $\xi\in\F$ be such that 
$$
\lim_n \frac{1}{n} \sum_{j=0}^{n-1} (\psi\circ F^k)_{-}(f^j(\xi))
=\int (\psi\circ F^k)_{-}\,d\mu_f
$$
for every $k\ge 1$.
Then
$\displaystyle{
\lim_n \frac{1}{n} \sum_{j=0}^{n-1} \psi(F^j(x))
=\int \psi \, d\mu_F}$
for every $x\in\xi$.
\end{lemma}

\begin{proof}
Let us fix $\psi$ and $\xi$ as in the statement. Then by
definition of $(\psi\circ F^k)_\pm$ and by the properties of
$\F$ we have
\[
(\psi\circ F^k)_-\big(f^j(\xi)\big)
\le
(\psi\circ F^k)\big(F^j(x)\big)
\le
(\psi\circ F^k)_+\big(f^j(\xi)\big)
\]
for all $x\in\xi$ and $j,k\ge1$. Given $\vep>0$, by
Corollary~\ref{c.poincareinvariant} there exists
$k_0\in\nat$ such that for all $k\ge k_0$
\[
\mu_F(\psi)-\frac\vep{2} \le
\mu_f\big((\psi\circ F^k)_-\big)
\le
\mu_f\big((\psi\circ F^k)_+\big)
\le
\mu_F(\psi) + \frac\vep{2} 
\]
and there is $n_0\in\nat$ such that for all $n\ge
n_0=n_0(k)$
\[
\left|
\frac1n\sum_{j=0}^{n-1} (\psi\circ
F^{k})_-\big(f^j(\xi)\big)
-
\mu_f\big((\psi\circ F^{k})_-\big)
\right|
<\frac\vep{2}.
\]
Hence we have that for all $n\ge n_0(k)$
\begin{align*}
  \mu_F(\psi)-\vep & \le 
\frac1n\sum_{j=0}^{n-1} (\psi\circ F^{k})( F^j(x) )
\\
& =
\frac{n+k}{n}\cdot\frac1{n+k}\sum_{j=0}^{n+k-1}
(\psi\circ F^{j})(x) - \frac1n\sum_{i=0}^{k-1} (\psi\circ
F^{j})(x)
\le
\mu_F(\psi)+\vep.
\end{align*}
Since $n$ can be made arbitrarily big and $\vep>0$ can be
taken as small as we want, we have concluded the proof of
the lemma.
\end{proof}

\begin{corollary}
\label{l.poincareergodic}
If $\mu_f$ is $f$-ergodic, then $\mu_F$ is ergodic for $F$.
\end{corollary}

\begin{proof}
Since $\Xi/\F$ is a metric space with the topology induced
  by $p$ we have that $C^0(\F,\real)$ is dense in
  $L^1(\F,\real)$ for the $L^1$-topology and $p:\Xi\to\cF$
  is continuous. Hence
there exists a subset $\cE$ of $\F$ with $\mu_f(\cE)=1$ such
that the conclusion of Lemma~\ref{l.SRBmSRBp} holds for a
subset $E=p^{-1}(\cE)$ of $\Xi$. To prove the corollary it
is enough to show that $\mu_F(E)=1$.

Let $\var=\chi_E=\chi_{\cE}\circ p$ and take
$\psi_n:\F\to\real$ a sequence of continuous functions such
that $\psi_n\to\chi_{\cE}$ when $n\to+\infty$ in the $L^1$
topology with respect to $\mu_f$. Then $\var_n=\psi_n\circ
p$ is a sequence of continuous functions on $\Xi$ such
that $\psi_n\to\psi$ when $n\to+\infty$ in the $L^1$ norm
with respect to $\mu_F$.

Then it is straightforward to check that
\[
\mu_F(\psi_n)=\lim_{k\to+\infty}\mu_f\Big( (\psi_n\circ
F^k)_- \Big)
=\lim_{k\to+\infty}\mu_f (\var_n\circ f^k) = \mu_f(\var_n)
\]
which converges to $\mu_f(\cE)=1$. Since $\mu_F(\psi_n)$
tends to $\mu_F(E)$ when $n\to+\infty$, we conclude that
$\mu_F(E)=1$, as we wanted.
\end{proof}


\subsection{Physical measure for the global Poincar\'e map}
\label{sec:phys-meas-glob}

Let us now apply these results (with $R$ replacing $F$) to
the case of the global Poincar\'e map for a
singular-hyperbolic attractor.

From the previous results in
Sections~\ref{sec:global-poincare-map}
and~\ref{sec:reduct-quot-map} the finitely many acim's
$\upsilon_1,\dots,\upsilon_l$ for the one-dimensional
quotient map $f$ uniquely induce
$R$-invariant ergodic probability measures
$\eta_1,\dots,\eta_l$ on $\Xi$.

We claim that the basins of each
$\eta_1,\dots,\eta_l$ have positive Lebesgue area
$\lambda^2$ on $\Xi$ and cover $\lambda^2$ almost every
point of $p^{-1}(I)$. Indeed the uniform  contraction 
of the leaves $\F_\Sigma^s\setminus\Gamma$ provided by
Lemma~\ref{stablereturnmap},  implies that the
forward time averages of any pair $x,y$ of points in
$\xi\in\F\setminus p(\Gamma)$ on continuous functions
$\vr:\Xi\to\real$ are equal
\[
\lim_{n\to+\infty}\left[
\frac1n\sum_{j=0}^{n-1}\vr\big( R^j(x) \big)
-
\frac1n\sum_{j=0}^{n-1}\vr\big( R^j(y) \big)
\right]=0.
\]
Hence $B(\eta_i)\supset p^{-1}\big( B(\upsilon_i)
\big), i=1,\dots,l$. This shows that $B(\eta_i)$
contains an entire strip except for a subset of
$\lambda^2$-null measure, because $B(\upsilon_i)$ contains
some open interval $\lambda$ modulo zero.  Since
$p_*(\lambda^2)\ll\lambda$ we get in particular
\[
\lambda^2\big( B(\eta_i)\big) >0 \quad\mbox{and}\quad
\lambda^2\Big( p^{-1}(I)\setminus\bigcup_{i=1}^l
B(\eta_i)\Big)=
p_*(\lambda^2)\Big(I\setminus\bigcup_{i=1}^l
B(\upsilon_i)\Big)=0,
\]
which shows that $\eta_1,\dots,\eta_l$ are
physical measures whose basins cover $p^{-1}(I)$ Lebesgue
almost everywhere.  We observe that $p^{-1}(I)\subset\Xi$ is
forward invariant under $R$, thus it contains
$\Lambda\cap\Xi$.


\subsection{Suspension flow from the Poincar\'e map}
\label{sec:reduct-poinc-map}

Let $\Xi$ be a measurable space, $\Gamma\,$ be some
measurable subset of $\Xi$, and $F:(\Xi\setminus
\Gamma)\to\Xi$ be a measurable map.  Let $\tau: \Xi\to
(0,+\infty]$ be a measurable function such that $\inf\tau
>0$ and $\tau\equiv +\infty$ on $\Gamma$.

Let $\sim$ be the equivalence relation on $\Xi\times
[0,+\infty)$ generated by $(x,\tau(x))\sim (F(x),0)$, that
is, $(x,s) \sim (\tilde{x},\tilde{s})$ if and only if there
exist
$$
(x,s)=(x_0,s_0),
\ (x_1, s_1),
\ \ldots,
\ (x_N,s_N)=(\tilde{x},\tilde{s})
$$
in $\Xi\times (0,+\infty)$ such that, for every $1\le i \le N$
\begin{equation*}
\begin{aligned}
\text{either}\quad & 
x_i=F(x_{i-1}) \quad\text{and}\quad s_i=s_{i-1}-\tau(x_{i-1});
\\
\text{or}\quad &   
x_{i-1}=F(x_i) \quad\text{and}\quad s_{i-1}=s_i-\tau(x_i).
\end{aligned}
\end{equation*}
We denote by $V=\Xi \times [0,+\infty)/\sim$ the
corresponding quotient space and by $\pi:\Xi\to V$ the
canonical projection which induces on $V$ a topology and a
Borel $\sigma$-algebra of measurable subsets of $V$.

\begin{definition}
\label{d.suspension}
The {\em suspension of $F$ with return-time $\tau$\/} is the
semi-flow $(X^t)_{t\ge 0}$ defined on $V$ by
$$
X^t(\pi(x,s))=\pi(x,s+t)
\quad\text{for every $(x,s)\in\Xi\times[0,+\infty)$ and $t>0$.}
$$
\end{definition}

It is easy to see that this is indeed well defined.

\begin{remark}
\label{r.fextendflow}
If $F$ is injective then we can also define
$$
X^{-t}\big(\pi(x,s)\big)
=\pi\big(F^{-n}(x),s+\tau(F^{-n}(x))+\cdots+\tau(F^{-1}(x))-t\big)
$$
for every $x\in F^n(\Xi)$ and $0 < t \le s +
\tau(F^{-n}(x))+\cdots+\tau(F^{-1}(x))$.  The expression on
the right does not depend on the choice of $n\ge 1$.  In
particular, the restriction of the semi-flow $(X_t)_{t\ge
  0}$ to the maximal invariant set
$$
\Lambda = \bigg\{(x,t): x\in\bigcap_{n\ge 0} F^n(\Xi)
              \text{\ \ and\ \ } t\ge 0\bigg\}
$$
extends, in this way, to a flow $(X^t)_{t\in\real}$ on $\Lambda$.
\end{remark}

Let $\mu_F$ be any probability measure on $\Xi$ that is
invariant under $F$.
Then the product $\mu_F\times dt$ of $\mu_F$ by Lebesgue
measure on $[0,+\infty)$ is an infinite measure, invariant
under the trivial flow $(x,s)\mapsto (x,s+t)$ in
$\Xi\times[0,+\infty)$.
In what follows we assume that the return time is integrable
with respect to $\mu_F$, i.e.
\begin{equation}
\label{eq.time}
\mu_F(\tau)=\int\tau\,d\mu_F < \infty.
\end{equation}
In particular $\mu_F(\Gamma)=0$.  Then we introduce the
probability measure $\mu_X$ on $V$ defined by
$$
\int \vr\,d \mu_X = \frac{1}{\mu_F(\tau)}\,
\int\int_0^{\tau(x)} \vr(\pi(x,t))\,dt \,  d\mu_F(x) 
$$
for each bounded measurable $\vr:V\to\real$. 

We observe that
the correspondence $\mu_F\mapsto\mu_X$ defined above is
injective. Indeed for any bounded measurable
$\psi:\Xi\to\real$, defining $\vr$ on
$\{x\}\times[0,\tau(x))$ to equal $\psi(x)$ gives a
bounded measurable map $\vr:V\to\real$ such that
$\mu_X(\vr)=\mu_F(\psi)$. Hence if $\mu_X=\mu_X'$ then
$\mu_F=\mu_F'$. 

\begin{lemma}
\label{p.flowinvariant}
The measure $\mu_X$ is invariant under the
semi-flow $(X^t)_{t\ge 0}$.
\end{lemma}

\begin{proof}
It is enough to show that $\mu_X\big((X^t)^{-1}(B)\big)=\mu_X(B)$
for every measurable set $B\subset V$ and any $0< t < \inf\tau$.
Moreover, we may suppose that $B$ is of the form $B=\pi(A\times J)$
for some $A\subset \Xi$ and $J$ a bounded interval in
$[0,\inf(\tau\mid A))$.
This is because these sets form a basis for the
$\sigma$-algebra of measurable subsets of $V$.

Let $B$ be of this form and $(x,s)$ be any point in
$\Xi$ with $0\le s <\tau(x)$.
Then $X^t(x,s)\in B$ if and only if
$\pi(x,s+t)=\pi(\tilde{x},\tilde{s})$ for some
$(\tilde{x},\tilde{s})\in A \times J$.
In other words, $(x,s)\in (X^t)^{-1}(B)$
if and only if there exists some $n\ge 0$ such that
$$
\tilde{x}=F^n(x)
\quad\text{and}\quad
\tilde{s}=s + t - \tau(x) - \cdots -\tau(F^{n-1}(x)).
$$
Since $s < \tau(x)$, $t < \inf\tau$, and $\tilde{s} \ge 0$,
it is impossible to have $n \ge 2$.
So,
\begin{itemize}
\item either $\tilde{x}=x$ and $\tilde{s} = s+t$
      (corresponding to $n=0$),
\item or $\tilde{x}=F(x)$ and $\tilde{s} = s+t - \tau(x)$
      (corresponding to $n=1$) 
\end{itemize}
The two possibilities are mutually exclusive:
for the first one $(x,s)$ must be such that $ s + t < \tau(x)$,
whereas in the second case $s+t \ge \tau(x)$.
This shows that we can write $(X^t)^{-}(B)$ as a disjoint union
$(X^t)^{-}(B)=B_1 \cup B_2$, with
\begin{equation*}
\begin{aligned}
B_1 & = \pi \big\{(x,s): x\in A \text{ and }
                     s \in (J-t) \cap [0,\tau(x))\big\}
\\
B_2 & = \pi \big\{(x,s): F(x)\in A \text{ and }
                     s \in (J+\tau(x)-t) \cap [0,\tau(x))\big\}.
\end{aligned}
\end{equation*}
Since $t>0$ and $\sup\, J <\tau(x)$, we have
$(J-t) \cap [0,\tau(x)) = (J-t) \cap [0, +\infty)$
for every $x\in A$.
So, by definition, $\mu_X(B_1)$ equals
$$
\int_{A} \ell\Big((J-t) \cap
[0,\tau(x)) \Big) \,  d\mu_F(x)
= \mu_F(A) \cdot \ell\Big((J-t)
\cap [0, +\infty)\Big).
$$
Similarly $\inf J\ge 0$ and $t<\tau(x)$ imply that
$$
(J + \tau(x) - t) \cap [0,\tau(x))
= \tau(x) + (J-t) \cap (-\infty,0).
$$
Hence $\mu_X(B_2)$ is given by 
\begin{equation*}
  \int_{F^{-1}(A)} \ell\Big((J - t) \cap
  (-\infty,0) \Big) \,  d\mu_F(x)
 = \mu_F(F^{-1}(A))\cdot \ell\Big(  (J-t) \cap
  (-\infty,0)\Big).
\end{equation*}
Since $\mu_F$ is invariant under $F$, we may replace
$\mu_F(F^{-1}(A))$ by $\mu_F(A)$ in the last expression.  It
follows that
$$
\mu_X\big((X^t)^{-1}(B)\big)
= \mu_X(B_1) + \mu_X(B_2)
= \mu_F(A) \cdot \ell\big( (J-t) \big).
$$
Clearly, the last term may be written as $\mu_F(A)\cdot\ell(J)$
which, by definition, is the same as $\mu_X(B)$.
This proves that $\mu_X$ is invariant under the semi-flow
and ends the proof.
\end{proof}

Given a bounded measurable function $\vr:V \to\real$, let
$\hat\vr:\Xi\to\real$ be defined by
\begin{equation}
\label{eq.hatvr}
\hat\vr(x)=\int_0^{\tau(x)} \vr(\pi(x,t))\,dt\,.
\end{equation}
Observe that $\hat\vr$ is integrable with respect to $\mu_F$
and by the definition of $\mu_X$
$$
\int \hat\vr \,d\mu_F
= \mu_F(\tau)\cdot \int \vr \,d\mu_X.
$$ 

\begin{lemma}
\label{l.SRBpSRBf}
Let $\vr:V \to\real$ be a bounded function, and $\hat\vr$ be
as above.  We assume that $x\in\Xi$ is such that
$\tau(F^j(x))$ and $\hat\vr(F^j(x))$ are finite for every $j
\ge 0$, and also
\begin{itemize}
\item[(a)] $\displaystyle{\lim_n \frac{1}{n} \sum_{j=0}^{n-1} \tau(F^j(x))
            =\int\tau\,d\mu_F}$, and
\item[(b)] $\displaystyle{\lim_n \frac{1}{n} \sum_{j=0}^{n-1} \hat\vr(F^j(x))
            =\int\hat\vr\,d\mu_F}$.
\end{itemize}
Then $\displaystyle{\lim_{T\to+\infty} \frac{1}{T}
                    \int_{0}^{T} \vr(\pi(x,s+t))\,dt
                   =\int \vr\,d\mu_X}\,$
for every $\,\pi(x,s)\in V$.
\end{lemma}

\begin{proof}
Let $x$ be fixed, satisfying (a) and (b).
Given any $T>0$ we define $n=n(T)$ by
$$
T_{n-1} \le T < T_n
\quad\text{where}\quad T_j=\tau(x)+\cdots\tau(F^j(x))
\text{ for }j\ge 0
$$
Then using $(y,\tau(y))\sim(F(y),0)$ we get
\begin{equation}
\label{eq.threeterms}
\begin{aligned}
  \frac{1}{T} \int_0^T \vr(\pi(x,s+t)) \,dt = \frac{1}{T}
  \left[ \sum_{j=0}^{n-1}\int_0^{\tau(F^j(x))}
    \vr(\pi(F^j(x),t)) \,dt\right.&
  \\
  \left. + \int_0^{T-T_{n-1}} \vr(\pi(F^n(x),t)) \,dt
    -
    \int_0^{s} \vr(\pi(x,t)) \,dt\right]& \,.
\end{aligned}
\end{equation}

Using the definition of $\hat\vr$, we may rewrite the first term
on the right hand side as
\begin{equation}
\label{eq.firstterm}
\frac{n}{T} \cdot \frac{1}{n} \sum_{j=0}^{n-1} \hat\vr(F^j(x)).
\end{equation}
Now we fix $\vep>0$. Assumption (a) and the definition of
$n$ imply that,
$$
n \cdot \Big(\int\tau\,d\mu_F-\vep\Big)
\le T_{n-1} \le T \le T_n
\le (n+1) \cdot \Big(\int\tau\,d\mu_F+\vep\Big),
$$
for every large enough $n$.
Observe also that $n$ goes to infinity as $T\to+\infty$,
since $\tau(F^j(x))<\infty$ for every $j$.
So, for every large $T$, 
$$
\mu_F(\tau) -\vep 
\le \frac{T}{n}
\le \frac{n+1}{n} \mu_F(\tau) + \vep
\le \mu_F(\tau) + 2\vep.
$$
This proves that $T/n$ converges to $\mu_F(\tau)$ when
$T\to+\infty$.
Consequently, assumption (b) implies that \reff{eq.firstterm}
converges to 
$$
\frac{1}{\mu_F(\tau)} \int \hat\vr\,d\mu_F 
= \int \vr\,d\mu_X\,. 
$$

Now we prove that the remaining terms in \reff{eq.threeterms}
converge to zero when $T$ goes to infinity.
Since $\vr$ is bounded
\begin{equation}
\label{eq.secondterm}
\left|\frac{1}{T} \int_0^{T-T_{n-1}} \vr(\pi(F^n(x),t)) \,dt\right|
\le \frac{T-T_{n-1}}{T} \, \sup|\vr|.
\end{equation}
Using the definition of $n$ once more,
$$
T-T_{n-1}
\le T_n -T_{n-1}
\le (n+1)\big(\int\tau\,d\mu_F+\vep\big) - 
n\big(\int\tau\,d\mu_F-\vep\big)
$$
whenever $n$ is large enough. Then
$$
\frac{T-T_{n-1}}{T} \le \frac{\int\tau\,d\mu_F +
  (2n+1)\vep}{n\big(\int\tau\,d\mu_F-\vep\big)} \le
\frac{4\vep}{\int\tau\,d\mu_F-\vep}
$$
for all large enough $T$.
This proves that $(T-T_{n-1})/T$ converges to zero, and then so
does \reff{eq.secondterm}.
Finally, it is clear that
$$
\frac{1}{T} \int_0^{s} \vr(\pi(x,t)) \,dt \to 0
\quad\text{when}\quad T\to+\infty.
$$ 
This completes the proof of the lemma.
\end{proof}

\begin{corollary}
\label{c.flowergodic}
If $\mu_F$ is ergodic then $\mu_X$ is ergodic.
\end{corollary}

\begin{proof}
  Let $\vr: V\to\real$ be any bounded measurable function,
  and $\hat\vr$ be as in \eqref{eq.hatvr}.  As already noted,
  $\hat\vr$ is $\mu_F$-integrable.  It follows that
  $\hat\vr(F^j(x))<\infty$ for every $j\ge 0$, at
  $\mu_F$-almost every point $x\in\Xi$.  Moreover, by the
  Ergodic Theorem, condition (b) in Lemma~\ref{l.SRBpSRBf}
  holds $\mu_F$-almost everywhere.  For the same reasons,
  $\tau(F^j(x))$ is finite for all $j\ge 0$, and condition
  (a) in the lemma is satisfied, for $\mu_F$-almost all
  $x\in\Xi$.

This shows that Lemma~\ref{l.SRBpSRBf} applies to every point
$x$ in a subset $A \subset \Xi$ with $\mu_F(A)=1$.
It follows that
$$
\lim_{T\to+\infty} \frac{1}{T}\int_0^T \vr(X^t(z))\,dt
= \int\vr\,d\mu_X
$$ 
for every point $z$ in $B=\pi(A\times[0,+\infty))$.
Since the latter has $\mu_X(B)=1$, we have shown that the
Birkhoff average of $\vr$ is constant $\mu_X$-almost everywhere.
Then the same is true for any integrable function, as bounded
functions are dense in  $L^1(\mu_X)$.
Thus $\mu_X$ is ergodic and the corollary is proved.
\end{proof}

\subsection{Physical measures for the suspension}
\label{sec:phys-meas-susp}

Using the results from Sections~\ref{sec:phys-meas-glob}
and~\ref{sec:reduct-poinc-map} it is straightforward to
obtain ergodic probability measures $\nu_1,\dots,\nu_l$
invariant under the suspension $(X^t)_{t\ge0}$ of $R$ with
return time $\tau$, corresponding to the $R$-physical
probability measures $\eta_1,\dots,\eta_l$ respectively.

Now we use Lemma~\ref{l.SRBpSRBf} to show that each $\nu_i$
is a physical measure for $(X^t)_{t\ge0}$, $i=1,\dots,l$.
Let $x\in\Sigma\cap B(\nu_i)$ for a fixed $\Sigma\in\Xi$ and
$i\in\{1,\dots, l\}$. According to
Remark~\ref{rmk:Poincareleaves} the return time $\tau_I$ on
$I$ is Lebesgue integrable, thus $\upsilon_i$-integrable
also since $\frac{d\upsilon_i}{d\lambda}$ is bounded. Hence
$\tau$ is $\eta_i$-integrable by the construction of
$\eta_i$ from $\upsilon_i$ (see
Section~\ref{sec:reduct-quot-map}).

Lemma~\ref{l.SRBpSRBf} together with the fact that $\eta_i$
is physical for $R$, ensures that $B(\nu_i)$ contains the
positive $X^t$ orbit of almost every point $(x,0), x\in
B(\nu_i),$ with respect to $\lambda^2$ on $B(\eta_i)$. If we
denote by $\lambda^3=\pi_*(\lambda^2\times dt)$ a natural
volume measure on $V$, then we get $\lambda^3\big( B(\nu_i)
\big)>0$.

This also shows that the basins $B(\nu_1),\dots, B(\nu_l)$
cover $\lambda^3$-almost every point in
$V_0=\pi\big(p^{-1}(I)\times[0,+\infty)\big)$. Notice that
this subset is a neighborhood of the suspension
$\pi\big((\Lambda\cap\Xi\setminus\Gamma)\times[0,+\infty)\big)$
of $\Lambda\cap\Xi\setminus\Gamma$.


\section{Physical measure for the flow}
\label{sec:phys-meas-flow}

Here we extend the previous conclusions of
Section~\ref{s.fsuspending} to the original
flow, completing the proof of Theorem~\ref{srb}.

We relate the suspension $(X^t)_{t\ge0}$ of $R$ with
return time $\tau$ to $(X_t)_{t\ge0}$ in $U$ as follows. We
define
\[
\Phi: \Xi\times[0,+\infty)\to U\quad
\mbox{by}\quad
(x,t)\mapsto X_t(x)
\]
and since
$\Phi\big(x,\tau(x)\big)=\big(R(x),0\big)\in\Xi\times\{0\}$,
this map naturally defines a quotient map
\begin{equation}
  \label{eq:quotient}
\phi:V\to U\quad\mbox{such that}\quad
\phi\circ X^t = X_t\circ\phi,\quad\mbox{for all}\quad
t\ge0,  
\end{equation}
through the identification $\sim$ from
Section~\ref{sec:reduct-poinc-map}. 

Let $\Xi_\tau=\{
(x,t)\in(\Xi\setminus\Gamma)\times[0,+\infty): 0< t
<\tau(x)\}$. Note that $\Xi_\tau$ is a open set in $V$ and
that $\pi\mid\Xi_\tau:\Xi_\tau\to \Xi_\tau$ is a
homeomorphism (the identity). Then the map $\phi\mid
\Xi_\tau$ is a local diffeomorphism into
$V_0=\phi\big(\Xi\times[0,+\infty)\big)\subset U$ by the
natural identification given by $\pi$ and by the Tubular
Flow Theorem, since points in $\Xi_\tau$ are not sent into
singularities of $X$. Notice that $\Xi_\tau$ is a full
Lebesgue ($\lambda^3$) measure subset of $V$. Thus $\phi$ is
a semiconjugation modulo zero. Note also that the number of
pre-images of $\phi$ is globally bounded by $r_0$ from
Remark~\ref{rmk:firstreturn}.

Therefore the measures $\nu_i$ constructed for the semiflow
$X^t$ in Section~\ref{sec:phys-meas-susp} define physical
measures $\mu_i=\phi_*(\nu_i)$, $i=1,\dots,l$, whose basins
cover a full Lebesgue ($m$) measure subset of $V_0$, which
is a neighborhood of $\Lambda$. Indeed the semiconjugacy
\eqref{eq:quotient} ensures that $\phi(B(\nu_i))\subset
B(\mu_i)$ and since $\phi$ is a local diffeomorphisms on a
full Lebesgue measure subset, then
\begin{align*}
  m\Big(\phi\big(B(\nu_1)\cup\dots\cup B(\nu_l)\big)\Big)=0.
\end{align*}
Since $V_0\subset U$ we have
\[
W^s(\Lambda)=\bigcup_{t<0} X_t(V_0).
\]
Moreover $X_t$ is a diffeomorphism for all $t\in\real$, thus
preserves subsets of zero $m$ measure.
Hence $\cup_{t<0} X_t\big( B(\mu_1)\cup\dots\cup B(\mu_l)
\big)$ has full Lebesgue measure in $W^s(\Lambda)$.
In other words, Lebesgue ($m$)
almost every point $x$ in the basin $W^s(\Lambda)$ of
$\Lambda$ is such that $X_t(x)\in B(\mu_i)$ for some $t>0$
and $i=1,\dots,l$.


\subsection{Uniqueness of the physical
  measure}
\label{sec:trans-uniq-phys}

The set $\Lambda$ is an attractor thus according to our
Definition~\ref{def:attractor} there exists $z_0\in\Lambda$
such that $\{X_t(z_0) : t>0\}$ is a dense regular orbit in
$\Lambda$.

We prove uniqueness of the physical measure by
contradiction, assuming that the number $l$ of distinct
physical measures is bigger than one.  Then we can take
distinct physical measures $\eta_1,\eta_2$ for $R$ on
$\Xi$ associated to distinct physical measures $\mu_1,\mu_2$
for $X\mid\Lambda$. Then there are open sets
$U_1,U_2\subset\Xi$ such that
\[
U_1\cap U_2=\emptyset \quad\mbox{and}\quad
\lambda^2\big(B(\eta_i)\setminus U_i\big)=0,\quad i=1,2.
\]
For a very small $\zeta>0$ we consider the open subsets
$V_i=X_{(-\zeta,\zeta)}(U_i), \, i=1,2$ of $U$ such that
$V_1\cap V_2=\emptyset$. According to the construction of
$\mu_i$ we have $\mu_i(B(\mu_i)\setminus V_i)=0, \, i=1,2$.

The transitivity assumption ensures that there are positive
times $T_1<T_2$ (exchanging $V_1$ and $V_2$ if needed) such
that $X_{T_i}(z_0)\in V_i, \, i=1,2$.  Since $V_1,V_2$ are
open sets and $g=X_{T_2-T_1}$ is a diffeomorphism, there
exists a small open set $W_1\subset V_1$ such that $g\mid
W_1 : W_1\to V_2$ is a $C^1$ diffeomorphism into its image
$W_2=g(W_1)\subset V_2$.

Now the $C^1$ smoothness of $g\mid W_1$ ensures that a full
Lebesgue ($m$) measure subset of $W_1$ is sent into a
full Lebesgue measure subset of $W_2$. By the definition of
$g$ and the choice of $V_1,V_2$, there exists a
 point in $B(\mu_1)\cap W_1$ whose positive orbit contains
a point in $B(\mu_2)\cap W_2$, thus $\mu_1=\mu_2$. Hence
\emph{singular-hyperbolic attractors have a unique physical
  probability measure $\mu$}.


\subsection{Hyperbolicity of the physical measure}
\label{sec:hyphysmeas}

For the hyperbolicity of the measure $\mu$ we note that
\begin{itemize}
\item the sub-bundle $E^s$ is one-dimensional and uniformly
  contracting, thus on the $E^s$-direction the Lyapunov
  exponent is negative for every point in $U$;
\item the sub-bundle $E^{cu}$ is two-dimensional, dominates
  $E^s$, contains the flow direction and is volume
  expanding, thus by Oseledets Theorem \cite{Man87,Wa82} the sum of the
  Lyapunov exponents on the direction of $E^{cu}$ is given by
  $\mu_{\,i}(\log|\det DX_1\mid E^{cu}|)>0$. Hence there is a
  positive Lyapunov exponent for $\mu_{\,i}$-almost every point
  on the direction of $E^{cu}$, $i=1,\dots, l$.
\end{itemize}

We will show that the expanding direction in $E^{cu}$ does not
coincide with the flow direction $E^X_z=\{ s\cdot X(z) :
s\in\real\},\, z\in \Lambda$.
Indeed, the invariant direction given by  $E^X_z$ cannot have
positive Lyapunov exponent, since for all $t>0$ and $z\in U$
\begin{equation}
  \label{eq:lyapquo}
\frac1t\log\big\| DX_t(z)\cdot X(z)
\big\|=
\frac1t\log\Big\| X\big( X_t(z) \big) \Big\| \le 
\frac1t \log\|X \|_0,  
\end{equation}
where $\|X\|_0=\sup\{ \| X(z) \| : z\in U\}$ is a constant.
Analogously this direction cannot have positive exponent for
negative values of time, thus the Lyapunov exponent along
the flow direction must be zero at regular points.


This shows that at $\mu$-almost every point $z$ the
Oseledets splitting of the tangent bundle has the form
\[
T_z M= E^s_z\oplus E^X_z\oplus F_z,
\]
where $F_z$ is the one-dimensional measurable sub-bundle of
vectors with positive Lyapunov exponent.  The proof of
Theorem~\ref{srb} is complete.


\section{Absolutely continuous disintegration of the
  physical measure}
\label{sec:phys-meas-are-1}

Here we prove Theorem~\ref{thm:srbmesmo}. We let $\mu$ be a
physical ergodic probability measure for a
singular-hyperbolic attractor $\Lambda$ of a $C^2$-flow in
an open subset $U\subset M^3$, obtained through the sequence
of reductions of the dynamics of the flow $X_t$ to the
suspension flow $X^t$ of the Poincar\'e map $R$ and return
time function $\tau$, with corresponding $X^t$-invariant
measure $\nu$ obtained from the $R$-invariant measure
$\eta$. In addition $\eta$ is obtained through the ergodic
invariant measure $\upsilon$ of the one-dimensional map
$f:I\to I$. This is explained in
Sections~\ref{sec:proof-theorem-B}
through~\ref{sec:phys-meas-flow}.  We know that $\mu$ is
hyperbolic as explained in Section~\ref{sec:phys-meas-flow}.

Let us fix $\de_0>0$ small. Then by Pesin's non-uniformly
hyperbolic theory \cite{Pe76,FHY83,PS89} we know that there
exists a compact subset $K\subset\Lambda$ such that
$\mu(\Lambda\setminus K) <\de_0$ and there exists $\de_1>0$
for which every $z\in K$ admits a strong-unstable manifold
$W^{uu}_{\de_1}(z)$ with inner radius $\de_1$. We refer to
this kind of sets as \emph{Pesin's sets}. The \emph{inner
  radius} of $W^{uu}_{\de_1}(z)$ is defined as the length of
the shortest smooth curve in this manifold from $z$ to its
boundary. Moreover $K\ni z\mapsto W^{uu}_{\de_1}(z)$ is a
continuous map $K\to{\cal E}^{1}(I_1,M)$ (recall the
notations in Section~\ref{s.21}).

The suspension flow $X^t$ defined on $V$ in
Section~\ref{sec:reduct-poinc-map} is conjugated to the
$X_t$-flow on an open subset of $U$ through a finite-to-1
local homeomorphism $\phi$, defined in
Section~\ref{sec:phys-meas-flow}, which takes orbits to
orbits and preserves time as in \eqref{eq:quotient}. Hence
there exists a corresponding set $K'=\phi^{-1}(K)$
satisfying the same properties of $K$ with respect to $X^t$,
where the constants $\de_0,\de_1$ are changed by at most a
constant factor due to $\phi^{-1}$ by the compactness of
$K$.
In what follows we use the measure $\nu=(\phi^{-1})_*\mu$
instead of $\mu$ and write $K$ for $K'$.

We fix a density point $x_0\in K$ of $\nu\mid K$. We may
assume that $x_0\in\Sigma$ for some $\Sigma\in\Xi$.
Otherwise if $x_0\not\in\Xi$, since $x_0=(x,t)$ for some
$x\in\Sigma$, $\Sigma\in\Xi$ and $0<t<T(x)$, then we use
$(x,0)$ instead of $x_0$ in the following arguments, but we
still write $x_0$. Clearly the length of the unstable
manifold through $(x,0)$ is unchanged due to the form of the
suspension flow, at least for small values of $\de_1$.
Since $\nu$ is given as a product measure on the quotient
space $V$ (see Section~\ref{sec:phys-meas-susp}), we may
assume without loss of generality that $x_0$ is a density
point of $\eta$ on $\Sigma\cap K$.

We set $W^u(x,\Sigma)$ to be the connected component of
$W^u(x)\cap\Sigma$ that contains $x$, for $x\in
K\cap\Sigma$, where $W^u(x)$ is defined in
Section~\ref{sec:phys-meas-are}. Recall that
$W^u(x)\subset\Lambda$. Then $W^u(x,\Sigma)$ has inner
radius bigger than some positive value $\de_2>0$ for $x\in
K\cap\Sigma$, which depends only on $\de_1$ and the angle
between $W^{uu}_{\de_1}(x)$ and $T_x\Sigma$. 

Let $\F^s(x_0,\de_2)=\{ W^s(x,\Sigma) : x\in W^u(x_0,\Sigma)
\}$ and $
F^s(x_0,\de_2)=\cup_{\gamma\in\F^s(x_0,\de_2)} \gamma$
be a horizontal
strip in $\Sigma$. Points $z\in F^s(x_0,\de_2)$ can be
specified using coordinates $(x,y)\in
W^u(x_0,\Sigma)\times\real$, where $x$ is given by
$W^u(x_0,\Sigma)\cap W^s(z,\Sigma)$ and $y$ is the length of
the shortest smooth curve connecting $x$ to $z$ in
$W^s(z,\Sigma)$.  Let us consider
\[
\F^u(x_0,\de_2)=\{ W^u(z,\Sigma):
z\in\Sigma\quad\mbox{and}\quad
W^u(z,\Sigma)\quad\mbox{crosses}\quad F^s(x_0,\de_2)\},
\]
where we say that a curve $\gamma$ \emph{crosses}
$F^s(x_0,\de_2)$ if the trace of $\gamma$ can be written as
the graph of a map $W^u(x_0,\Sigma)\to W^s(x_0,\Sigma)$
using the coordinates outlined above. We stress that
$\F^u(x_0,\de_2)$ is not restricted to leaves through points
of $K$.

We may assume that $F^u(x_0,\de_2)=\cup\F^u(x_0,\de_2)$
satisfies $\eta(F^u(x_0,\de_2))>0$ up to taking a smaller
$\de_2>0$, since $x_0$ is a density point of $\eta\mid
K\cap\Sigma$.  Let $\hat\eta$ be the measure on
$\F^u(x_0,\de_2)$ given by
\[
\hat\eta(A)=\eta\Big(\bigcup_{\gamma\in A} \gamma\Big)
\quad\mbox{for every measurable set}\quad A\subset\F^u(x_0,\de_2).
\]

\begin{proposition}
  \label{pr:densidades}
  The measure $\eta\mid F^u(x_0,\de_2)$ admits a
  disintegration into conditional measures $\eta_\gamma$
  along $\hat\eta$-a.e. $\gamma\in \F^u(x_0,\de_2)$
  such that $\eta_\gamma\ll \lambda_\gamma$, where
  $\lambda_\gamma$ is the measure (length) induced on
  $\gamma$ by the natural Riemannian measure $\lambda^2$
  (area) on $\Sigma$. Moreover there exists $D_0>0$ such
  that
\[
\frac1{D_0}\le\frac{d\eta_\gamma}{d\lambda_\gamma}\le D_0,
\quad
\mbox{$\eta_\gamma$-almost everywhere for $\hat\eta$-almost every
$\gamma$.}
\]
\end{proposition}
This is enough to conclude the proof of
Theorem~\ref{thm:srbmesmo} since both $\de_0$ and $\de_2$
can be taken arbitrarily close to zero, so that all unstable
leaves $W^u(x,\Sigma)$ through almost every point with
respect to $\eta$ will support a conditional measure of
$\eta$. 

\begin{figure}[htpb]
  \centering
  \psfrag{S}{$\Sigma$}\psfrag{F}{$F^s$}
  \psfrag{g}{$\gamma\times[0,T]$}
  \psfrag{m}{$\gamma$}
  \includegraphics[height=3cm,width=5cm]{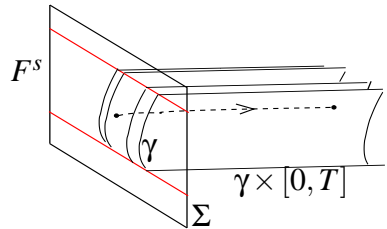}
  \caption{Center-unstable leaves on the suspension flow.}
  \label{fig:center-unstable-leav}
\end{figure}

Indeed, to obtain the disintegration of $\nu$ along the
center-unstable leaves that cross any small ball around a
density point $x_0$ of $K$, we project that neighborhood of
$x_0$ along the flow in negative time on a cross section
$\Sigma$. Then we obtain the family $\{\eta_\gamma\}$, the
disintegration of $\eta$ along the unstable leaves
$\gamma\in\F^u$ on a strip $F^s$ of $\Sigma$, and consider
the family $\{ \eta_\gamma\times dt\}$ of measures on
$\F^u\times[0,T]$ to obtain a disintegration of $\nu$, where
$T>0$ is a fixed time slightly smaller than the return time
of the points in the strip $F^s$, see
Figure~\ref{fig:center-unstable-leav}.

In fact, $\eta_\gamma\times dt \ll \lambda_\gamma\times dt$
and $\lambda_\gamma\times dt$ is the induced (area) measure
on the center-unstable leaves by the volume measure
$\lambda^3$ on $V$, and it can be given by restricting the
volume form $\lambda^3$ to the surface $\gamma\times[0,T]$
which we write $\lambda^3_\gamma$, for $\gamma\in\F^u$.
Thus by Proposition~\ref{pr:densidades} and by the
definition of $\nu$, we have
\[
\nu_\gamma=\eta_\gamma\times dt =
\frac{d\eta_\gamma}{d\lambda_\gamma}\cdot \lambda^3_\gamma,
\quad \gamma\in\F^u
\]
and the densities of the conditional measures
$\eta_\gamma\times dt$ with respect to $\lambda^3_\gamma$
are also uniformly bounded from above and from below away
from zero -- we have left out the constant factor
$1/\mu(\tau)$ to simplify the notation.

Since $\mu=\phi_*\nu$ and $\phi$ is a finite-to-1 local
diffeomorphism when restricted to $\Xi_\tau$, then $\mu$
also has an absolutely continuous disintegration along the
center-unstable leaves.  The densities on unstable leaves
$\gamma$ are related by the expression (where $m_\gamma$
denotes the area measure on the center-unstable leaves
induced by the volume form $m$)
\[
\mu_{\gamma}=\phi_*(\nu_\gamma)=
\phi_*\Big(\frac{d\eta_\gamma}{\lambda_\gamma}\cdot
\lambda^3_\gamma\Big) = 
\left( \frac{1}{\det D(\phi\mid\gamma\times[0,T])}
\cdot\frac{d\eta_\gamma}{\lambda_\gamma}
\right)\circ \phi^{-1}\cdot m_\gamma, \quad
\gamma\in\F^u
\]
which implies that the densities along the center-unstable
leaves are uniformly bounded from above.  

Indeed observe first that the number of pre-images of $x$
under $\phi$ is uniformly bounded by $r_0$ from
Remark~\ref{rmk:firstreturn}, i.e. by the number of
cross-sections of $\Xi$ hit by the orbit of $x$ from time
$0$ to time $t_2$.  Moreover the tangent bundle of
$\gamma\times[0,T]$ is sent by $D\phi$ into the bundle
$E^{cu}$ by construction and recalling that
$\phi(x,t)=X_t(x)$ then, if $e_1$ is a unit tangent vector
at $x\in\gamma$, $\hat e_1$ is the unit tangent vector at
$\phi(x,0)\in W^u(x,\Sigma)$ and $e_2$ is the flow direction
at $(x,t)$ we get
\[
D\phi(x,t)(e_1)=DX_t\big(X_t(x)\big)(\hat e_1)
\quad\mbox{and}\quad
D\phi(x,t)(e_2)=DX_t\big(X_t(x)\big)\big(X(x,0)\big)=
X\big(X_t(x)\big).
\]
Hence $D\big(\phi\mid\gamma\times[0,T]\big)(x,t)=DX_t\mid
E^{cu}_{\phi(x,t)}$ for $(x,t)\in\gamma\times[0,T]$ and so
\[
|\det D\big(\phi\mid\gamma\times[0,T]\big)(x,t)|=J^c_t(x).
\]

Now the volume expanding property of $X_t$ along the
center-unstable sub-bundle, together with the fact that the
return time function $\tau$ is not bounded from above near
the singularities, show that the densities of $\mu_\gamma$
are uniformly bounded from above throughout $\Lambda$ but
not from below. In fact, this shows that these densities
will tend to zero close to the singularities of $X$ in
$\Lambda$.

This finishes the proof of Theorem~\ref{thm:srbmesmo} except
for the proof of Proposition~\ref{pr:densidades} and of
$\suporte(\mu)=\Lambda$, which we present in what follows.

\subsection{Constructing the disintegration}
\label{sec:constr-disint}

Here we prove Proposition~\ref{pr:densidades}. We split the
proof into several lemmas keeping the notations of the
previous sections. 

Let $\lambda^2$, $R:p^{-1}(I)\to\Xi$, $\F^u(x_0,\de_2)$,
$F^u(x_0,\de_2)$ and $\eta$ be as before, where $x_0\in
K\cap\Sigma$ is a density point of $\eta\mid K$ and $K$ is a
compact Pesin set.  We write $\{\eta_\gamma\}$ and
$\{\lambda^2_\gamma\}$ for the disintegrations of $\eta\mid
F^u(x_0,\de_2)$ and $\lambda^2$ along
$\gamma\in\F^u(x_0,\de_2)$.

\begin{lemma}
  \label{le:tudounada}
 Either $\eta_{\,\gamma}\ll\lambda^2_\gamma\,$ for
  $\hat\eta$-a.e. $\gamma\in\F^u(x_0,\de_2)$,
 or $\eta_{\,\gamma}\perp \lambda^2_\gamma\,$ for
  $\hat\eta$-a.e. $\gamma\in\F^u(x_0,\de_2)$.
\end{lemma}

\begin{proof}
  We start by assuming that the first item in the statement
  does not hold and proceed to show that this implies the
  second item. We write $\eta$ for $\eta(
  F^u(x_0,\de_2))^{-1}\cdot\eta\mid F^u(x_0,\de_2)$ to
  simplify the notation in this proof.
  
  Let us suppose that there exists $A\subset F^u(x_0,\de_2)$
  such that $\eta(A)>0$ and $\lambda^2_\gamma(A)=0$ for
  $\hat\eta$-a.e.  $\gamma\in\F^u(x_0,\de_2)$. Let
  $B=\cup_{k\ge0} R^k(A)$. We claim that $\eta(B)=1$.
  
  Indeed, we have $R(B)\subset B$, then $B\subset R^{-1}(B)$
  and $\big(R^{-k}(B)\big)_{k\ge0}$ is a nested increasing
  family of sets.  Since $\eta$ is $R$-ergodic we have for
  any measurable set $C\subset \Xi$
  \begin{equation}
    \label{eq:mediaergodica}
    \lim_{n\to+\infty}
    \frac1n\sum_{j=0}^{n-1} \eta\big(C\cap R^{-j}(B)\big)
    = \eta(C)\cdot\eta(B).
  \end{equation}
  But $\eta\big(\cup_{k\ge0} R^{-k}(B)\big)=1$ because
  this union is $R$-invariant and
  $\eta(B)=\eta\big(R^{-k}(B)\big)>0$ by assumption, for
  any $k\ge0$. Because the sequence is increasing and nested
  we have $\eta\big( R^{-k}(B)\big)\nearrow1$. Hence
  from~\eqref{eq:mediaergodica} we get that
  $\eta(C)=\eta(C)\cdot\eta(B)$ for all sets $C\subset X.$
  Thus $\eta(B)=1$ as claimed.

Therefore $1=\eta(B)=\int\eta_{\,\gamma}(B)\,
d\hat\eta(\gamma)$ and so $\eta_{\,\gamma}(B)=1$ for
$\hat\eta$-a.e. $\gamma\in\F^u(x_0,\de_2)$ since every
measure involved is a probability measure.

We now claim that $\lambda^2_\gamma(B)=0$ for $\hat\mu$-a.e.
$\gamma\in\F^u(x_0,\de_2)$. For if
$R(A)\cap\gamma\neq\emptyset$ for some
$\gamma\in\F^u(x_0,\de_2)$, then $A\cap R^{-1}(\gamma)\cap
F^u(x_0,\de_2)\neq\emptyset$ and so it is enough to consider
only $A\cap F^u_1$, where
$F^u_1=R^{-1}(F^u(x_0,\de_2))\cap F^u(x_0,\de_2)$.  But
$\lambda^2_\gamma(A\cap F^u_1)\le\lambda^2_\gamma(A)=0$ thus
\[
0= \lambda^2_\gamma\big(R_0(A\cap F^u_1)\big) \ge
\lambda^2_\gamma\big(R_0(A)\cap F^u(x_0,\de_2)\big) =
\lambda^2_\gamma(R_0(A))
\]
for $\hat\eta$-a.e. $\gamma$ since $R_0$ is piecewise
smooth, hence a regular map. Therefore we get
$\lambda^2_\gamma(R^k(A))=0$ for all $k\ge1$ implying that
$\lambda^2_\gamma(B)=0$ for $\hat\eta$-a.e. $\gamma$.

This shows that $\eta_\gamma$ is singular with respect to
$\lambda^2_\gamma$ for $\hat\eta$-a.e. $\gamma$. The proof is
finished.
\end{proof}


\subsubsection{Existence of hyperbolic times
  for $f$ and consequences to $R$}
\label{sec:exist-prop-hyperb}

Now we show that a positive measure subset of
$\F^u(x_0,\de_2)$ has absolutely continuous disintegrations,
which is enough to conclude the proof of
Proposition~\ref{pr:densidades} by Lemma~\ref{le:tudounada},
except for the bounds on the densities.


We need the notion of \emph{hyperbolic time} for the
one-dimensional map $f$~\cite{ABV00}. We know that this map
is piecewise $C^{1+\alpha}$ and the boundaries $\Gamma_0$ of
the intervals $I_1,\dots, I_n$ can be taken as a
\emph{singular set} for $f$ (where the map is not defined
or is not differentiable) which behaves like a
\emph{power of the distance to $\Gamma_0$}, as follows.
Denoting by $d$ the usual distance on the intervals $I$,
there exist $B>0$ and $\beta>0$ such that
\begin{itemize}
\item $\frac1B \cdot d(x,\Gamma_0)^\beta
\le \big| f' \big| \le B\cdot d(x,\Gamma_0)^{-\beta}$;
\item $\big| \log|f'(x)| - \log|f'(y)| \big|\le
B\cdot d(x,y)\cdot d(x,\Gamma_0)^{-\beta}$,
\end{itemize}
for all $x,y\in I$ with $d(x,y)<d(x,\Gamma_0)/2$. This is
true of $f$ since in Section~\ref{sec:reduct-quot-leaf}
it was shown that $f'\mid I_j$ either is bounded from
above and below away from zero, or else is of the form
$x^\beta$ with $\beta\in(0,1)$.

Given $\de>0$ we define $d_\delta(x,\Gamma_0)=d(x,\Gamma_0)$ if
$d(x,\Gamma_0)<\de$ and $1$ otherwise.

\begin{definition}
  \label{def:tempo-hiperbolico}
  Given $b,c,\delta>0$ we say that $n\ge1$ is a
  $(b,c,\delta)$-hyperbolic time for $x\in I$ if
  \begin{equation}
    \label{eq:tempo-hip}
\prod_{j=n-k}^{n-1}\big| f'\big(f^j(x)\big)\big|^{-1} \le e^{-ck}
\quad\mbox{and}\quad
\prod_{j=n-k}^{n-1} d_\delta \big(
f^j(x),\Gamma_0\,\big) \ge e^{-bk}
      \end{equation}
      for all $k=0,\dots,n-1$.
\end{definition}

Since  $f$ has positive Lyapunov exponent
$\upsilon$-almost everywhere, i.e.
\[
\lim_{n\to+\infty}\frac1n\log\big| (f^n)'(x)\big| > 0
\quad\mbox{for   }\upsilon\mbox{-almost all   } x\in I, 
\]
and $\frac{d\upsilon}{d\lambda}$ is bounded from above
(where $\lambda$ is the Lebesgue length measure on $I$),
thus $|\log d(x,\Gamma_0)|$ is $\upsilon$-integrable and for
any given $\epsilon>0$ we can find $\de>0$ such that for
$\upsilon$-a.e. $x\in I$
\[
  \lim_{n\to\infty}\frac1n\sum_{j=0}^{n-1}
  -\log d_\delta(f^j(x),\Gamma_0)=
  \int -\log d_\de(x,\Gamma_0)\, d\upsilon(x) <\epsilon.  
\]
This means that $f$ is \emph{non-uniformly expanding} and
has \emph{slow recurrence to the singular set}. Hence we are
in the setting of the following result.

\begin{theorem}[Existence of a positive frequency of
  hyperbolic times]
\label{le:tempos-hip-existem}
Let $f:I\to I$ be a $C^{1+\alpha}$ map, behaving like a
power of the distance to a singular set $\Gamma_0$, 
non-uniformly expanding and with slow recurrence to
$\Gamma_0$ with respect to an absolutely continuous
invariant probability measure $\upsilon$. Then for $b,c,\de>0$
small enough there exists $\theta=\theta(b,c,\de)>0$ such
that $\upsilon$-a.e. $x\in I$ has infinitely many
$(b,c,\de)$-hyperbolic times.  Moreover if we write
$0<n_1<n_2<n_2<\dots$ for the hyperbolic times of $x$ then
their asymptotic frequency satisfies
\[
\liminf_{N\to\infty}\frac{\#\{ k\ge1 : n_k\le
  N\}}{N}\ge\theta
\quad\mbox{for}\quad \upsilon\mbox{-a.e.  } x\in I.
\]
\end{theorem}

\begin{proof}
A complete proof can be found in~\cite[Section 5]{ABV00} with
weaker assumptions corresponding to Theorem C in that paper.
\end{proof}

From now on we fix  values of $(b,c,\de)$ so that the
conclusions of Theorem~\ref{le:tempos-hip-existem} are true.

We now outline the properties of these special times. For detailed
proofs see~\cite[Proposition 2.8]{ABV00} and~\cite[Proposition
2.6, Corollary 2.7, Proposition 5.2]{AA03}.

\begin{proposition}
  \label{pr:prophyptimes}
  There are constants $\beta_1,\beta_2>0$ depending on
  $(b,c,\de)$ and $f$ only such that, if $n$ is
  $(b,c,\delta)$-hyperbolic time for $x\in I$, then there
  are neighborhoods $W_k(x)\subset I$ of $f^{n-k}(x)$,
  $k=1,\dots, n$, such that
\begin{enumerate} 
\item $f^k\mid W_k(x)$ maps $W_k(x)$ diffeomorphically to
  the ball of radius $\beta_1$ around $f^n(x)$;
\item   for every $1\leq k\leq n$ and $y,z\in W_k(x)$
\[
d\big(f^{n-k}(y),f^{n-k}(z)\big)\le
  e^{-ck/2}\cdot d\big(f^{n}(y),f^{n}(z)\big);
\]
\item for $y,z\in W_n(x)$
\[
\frac1{\beta_2}\le
\frac{\big|(f^n)'(y)\big|}{\big|(f^n)'(z)\big|} \le
\beta_2.
\]
\end{enumerate}
\end{proposition}

The conjugacy $p\circ R = f\circ p$ between the actions
of the Poincar\'e map and the one-dimensional map on the
space of leaves, together with the bounds on the
derivative~~\eqref{eq:DpSigma}, enables us to extend
the properties given by Proposition~\ref{pr:prophyptimes} to
any $cu$-curve inside $B(\eta)$, as follows.

Let $\gamma:J\to\Xi$ be a $cu$-curve in $\Xi\setminus\Gamma$
such that $\gamma(s)\in B(\eta)$ for Lebesgue almost every
$s\in J$, $J$ a non-empty interval --- such a curve exists
since the basin $B(\eta)$ contains entire strips of some
section $\Sigma\in\Xi$ except for a subset of zero area.
Note that we have the following limit in the weak$^*$
topology
\[
\lim_{n\to+\infty} \lambda_\gamma^{\,n} = \eta\quad
\mbox{where}\quad
\lambda_\gamma^{\,n}=\frac1n\sum_{j=0}^{n-1}
R^j_*(\lambda_\gamma),
\]
by the choice of $\gamma$ and by an easy application of the
Dominated Convergence Theorem.

\begin{proposition}
  \label{pr:prop-cu-hyptimes}
  There are constants $\kappa_0,\kappa_1>0$ depending on
  $(b,c,\de)$ and $R_0,\beta_0,\beta_1,\beta_2$ only such
  that, if $x\in\gamma$ and $n$ is big enough and a
  $(b,c,\delta)$-hyperbolic time for $p(x)\in I$, then there
  are neighborhoods $V_k(x)$ of $R^{n-k}(x)$ on
  $R^{n-k}(x)(\gamma)$, $k=1,\dots, n$, such that
\begin{enumerate} 
\item $R^k\mid V_k(x)$ maps $V_k(x)$ diffeomorphically
  to the ball of radius $\kappa_0$ around
  $R^n(x)$ on $R^n(\gamma)$;
\item for every $1\leq k\leq n$ and $y,z\in V_k(x)$
\[
d_{R^{n-k}(\gamma)}\big(R^{n-k}(y),R^{n-k}(z)\big)\le
\beta_0\cdot e^{-ck/2}\cdot
d_{R^n(\gamma)}\big(R^{n}(y),R^{n}(z)\big);
\]
\item for $y,z\in V_n(x)$
\[
\frac1{\kappa_1}\le\frac{\big|D(R^n\mid
  \gamma)(y)\big|}{\big|D(R^n\mid \gamma)(z)\big|} \le
  \kappa_1;
\]
\item the inducing time of $R^k$ on $V_k(x)$ is constant,
  i.e.  $r^{n-k}\mid V_k(x)\equiv$ const..
\end{enumerate}
\end{proposition}

Here $d_\gamma$ denotes the distance along $\gamma$
given by the shortest smooth curve in $\gamma$ joining two
given points and $\lambda_\gamma$ denotes the normalized Lebesgue
length measure induced on $\gamma$ by the
area form $\lambda^2$ on $\Xi$.


\begin{proof}[Proof of Proposition~\ref{pr:prop-cu-hyptimes}]
  Let $x_0=p(x)$ and $W_k(x_0)$ be given by
  Proposition~\ref{pr:prophyptimes}, $k=1,\dots n$. We have
  that $p(\gamma)$ is an interval in $I$ and that
  $p\mid\gamma:\gamma\to p(\gamma)$ is a diffeomorphism ---
  we may take $\gamma$ with smaller length if needed.

 If $n$ is big enough, then $W_n(x_0)\subset
 p(\gamma)$. Moreover the conjugacy implies that the
 following maps are all diffeomorphisms
\[
\begin{array}[c]{rcl}
V_k(x) & \stackrel{R^k}{\longrightarrow} & R^k(V_k(x))
\\
p \downarrow \quad & & \quad \downarrow p
\\
W_k(x_0) & \stackrel{f^k}{\longrightarrow}
& B\big(f^k(x_0),\kappa_0\big)
\end{array},
\]
and the diagram commutes, where $V_k(x)=\big(p\mid
R^k(\gamma)\big)^{-1} \big(W_k(x_0)\big)$, $k=1,\dots,n$,
see Figure~\ref{fig:hyperb-times-project}.
Using the bounds~\eqref{eq:DpSigma} to compare derivatives
we get $\kappa_0=\beta_1/\beta_0$ and
$\kappa_1=\beta_0\cdot\beta_2$.

\begin{figure}[htpb]
  \centering
  \includegraphics[height=2cm,width=9cm]{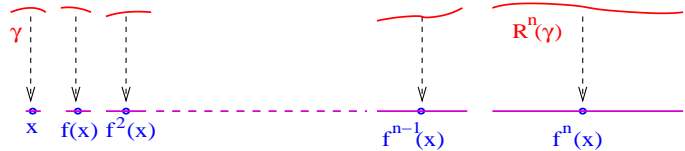}
  \caption{Hyperbolic times and projections.}
  \label{fig:hyperb-times-project}
\end{figure}

To get item (4) we just note that by definition of
$(b,c,\delta)$-hyperbolic time none of the sets $W_k(x_0)$
may intersect $\Gamma_0$. According to the definition of
$\Gamma_0$, this means that orbits through $x,y\in V_k(x)$
cannot cut different cross-sections in $\Xi$ before the next
return in time $\tau(x),\tau(y)$ respectively. Hence every
orbit through $W_k(x_0)$ cuts the same cross-sections in its
way to the next return cross-section.  In particular the
number of cross-section cuts is the same, i.e. $r\mid
V_k(x)$ is constant, $k=1,\dots,n$. Hence by definition of
$r^k$ we obtain the statement of item (4) since
$R(V_k(x))=V_{k-1}(x)$ by definition. This completes the
proof of the proposition.
\end{proof}

\subsubsection{Approximating $\eta$ by push forwards
  of Lebesgue measure at hyperbolic times}
\label{sec:approx-tildeeta-push}

We define for $n\ge1$
\[
H_n=\{ x\in\gamma : n \mbox{   is a  }
(b,c/2,\de)\mbox{-hyperbolic time for   } p(x)
\}.
\]
As a consequence of items (1-2) of
Proposition~\ref{pr:prop-cu-hyptimes}, we have that $H_n$ is
an open subset of $\gamma$ and for any $x\in\gamma\cap H_n$
we can find a connected component $\gamma^{\,n}$ of
$R^n(\gamma)\cap B(R^n(x),\kappa_0)$ containing $x$ such
that $R^n\mid V_n(x) : V_n(x)\to\gamma^{\,n}$ is a
diffeomorphism. In addition $\gamma^{\,n}$ is a $cu$-curve
according to Corollary~\ref{ccone}, and by item (3) of
Proposition~\ref{pr:prop-cu-hyptimes} we deduce that
\begin{equation}
  \label{eq:denslimitadahyp}
  \frac1{\kappa_1}\le
\frac{d\big(R^n_*(\lambda_\gamma)\mid B(R^n(x),\kappa_0)\big)}
{d\lambda_{\gamma^{\,n}}} 
\le \kappa_1,\quad
\lambda_{\gamma^{\,n}}-\mbox{a.e. on  } 
\gamma^{\,n},
\end{equation}
where $\lambda_{\gamma^{\,n}}$ is the Lebesgue induced
measure on $\gamma^{\,n}$ for any $n\ge1$, if we normalize
both measures so that $\big((R^n)_*(\lambda_\gamma)\mid
B(R^n(x),\kappa_0)\big)(\gamma^{\,n})=
\lambda_{\gamma^{\,n}}(\gamma^{\,n})$, i.e. their masses on
$\gamma^{\,n}$ are the same.

Moreover the set $R^n(\gamma\cap H_n)$ has an at most
countable number of connected components which are
diffeomorphic to open intervals. Each of these components is
a $cu$-curve with diameter bigger than $\kappa_0$ and hence
we can find a \emph{pairwise disjoint family
  $\gamma^{\,n}_i$ of $\kappa_0$-neighborhoods around
  $R^n(x_i)
  $ in $R^n(\gamma)$, for some
  $x_i\in H_n$, with maximum cardinality}, such that
\begin{equation}
  \label{eq:massaboa}
\Delta_n=\bigcup_i \gamma^{\,n}_i\subset R^n\big(\gamma\cap
H_n\big)
\quad\mbox{and}\quad
\big((R)^n_*(\lambda_\gamma)\mid\Delta_n\big)(\Delta_n)\ge
\frac1{2\kappa_1}\cdot\lambda_\gamma(H_n).
\end{equation}
Indeed since $R^n(\gamma\cap H_n)$ \emph{is one-dimensional}, for
each connected component the family $\Delta_n$ may miss a
set of points of length at most equal to the length of one
$\gamma^{\,n}_i$, for otherwise we would manage to include an
extra $\kappa_0$-neighborhood in $\Delta_n$. Hence we have
in the worst case (assuming that there is only one set
$\gamma^{\,n}_i$ for each connected component)
\[
\lambda_{\gamma^{\,n}}\big(R^n(\gamma\cap
H_n)\setminus\Delta_n\big) \le
\lambda_{\gamma^{\,n}}\big(\bigcup_i \gamma^{\,n}_i \big)
=\lambda_{\gamma^{\,n}}(\Delta_n) \quad\mbox{so that}\quad
\lambda_{\gamma^{\,n}}(\Delta_n)\ge\frac12\cdot
\lambda_{\gamma^{\,n}}\big(R^n(\gamma\cap H_n)\big)
\]
and the constant $\kappa_1$ comes from \eqref{eq:denslimitadahyp}.

For a fixed small $\rho>0$ we consider $\Delta_{n,\rho}$
given by the balls $\gamma^{\,n}_i$ with the same center
$x_{n,i}$ but a reduced radius of $\kappa_0-\rho$. Then the
same bound in~\eqref{eq:massaboa} still holds with
$2\kappa_1$ replaced by $3\kappa_1$. 

We write $D_n$ for the family of disks
from $\cup_{j\ge1}\Delta_j$ with the same expanding iterate
(the disks with the same centers as the ones from
$D_{n,\rho}$ but with their original size).  

We define the following sequences of measures
\[
\omega_{\rho}^n=\frac1n\sum_{j=0}^{n-1}
R^j_*(\lambda_\gamma)\mid \Delta_{j,\rho}
\quad\mbox{and}\quad
\overline\lambda_\gamma^{\,n}=\lambda_\gamma^{\,n}-
\omega_{\rho}^n,\quad n\ge1.
\]
Then 
any weak$^*$ limit point $\tilde\eta=\lim_k
\omega_{\rho}^{n_k}$ for some subsequence $n_1<n_2<\dots$
and $\overline\eta=\lim_k \overline\lambda_\gamma^{\,n'_k}$
(where $n_k'$ may be taken as a subsequence of $n_k$), are
$R$-invariant measures which satisfy
$\eta=\tilde\eta+\overline\eta$.

We claim that $\tilde\eta\not\equiv0$, thus
$\eta=\tilde\eta$ as a consequence of the ergodicity of
$\eta$.  In fact, we can bound the mass of $\omega_{\rho}^n$
from below using the density of hyperbolic times from
Theorem~\ref{le:tempos-hip-existem} and the bound
from~\eqref{eq:massaboa} 
through the following Fubini-Toneli-type argument.
Write $\#_n(J)=\# J/n$ for any $J\subset\{0,\dots,n-1\}$,
the uniform discrete measure on the first $n$ integers. Also
set $\chi_{i}(x)=1$ if $x\in H_i$ and zero otherwise,
$i=0,\dots,n-1$. Then
\begin{eqnarray*}
\omega_{\rho}^n(M)
&\ge&
\frac1{3\kappa_1\cdot n}
\sum_{j=0}^{n-1}\lambda_\gamma(H_j)
=
\frac{n}{3\kappa_1 n}
\int\!\!\int \chi_{i}(x)
\, d\lambda_\gamma(x) \, d\#_{n}(i)
\\
&=&
\frac{1}{3\kappa_1}
\int\!\!\int \chi_{i}(x)  \, d\#_{n}(i) \,
d\lambda_\gamma(x) 
\ge \frac{\theta}{6\kappa_1}>0,
\end{eqnarray*}
for every $n$ big enough by the choice of $\gamma$. 


\subsubsection{Approximating unstable curves by images of curves
  at hyperbolic times}
\label{sec:appr-limit-curv}

We now observe that since $\eta(F^u(x_0,\de_2))>0$ and $x_0$
is a density point of $\eta\mid F^u(x_0,\de_2)$, then
$\omega_{\rho}^n(F^u(x_0,\de_2))\ge c$ for some constant
$c>0$ for all big enough $n$.
If we assume that $\de_2<\rho$, which poses no restriction,
then we see that the $cu$-curves from $D_{j,\rho}$
intersecting $F^u(x_0,\de_2)$ will cross this horizontal
strip when we restore their original size. Thus the leaves
$\cup_{j=0}^{n-1} D_j$ in the support of
$\omega^n_0$ which intersect $F^u(x_0,\de_2)$
cross this strip.  Given any sequence
$\gamma^{\,n_k}$ of leaves in $D_{n_k}$
crossing $F^u(x_0,\de_2)$ with $n_1<n_2<n_3<\dots$, then
there exists a $C^1$-limit leaf $\gamma^\infty$ also
crossing $F^u(x_0,\de_2)$, by the Ascoli-Arzela Theorem. We
claim that this leaf coincides with the unstable manifold of
its points, i.e. $\gamma^\infty=W^u(x,\Sigma)$ for all
$x\in\gamma^\infty$.  This shows that the accumulation
curves $\gamma^\infty$ are defined independently of the
chosen sequence $\gamma^{\,n_k}$ of curves in
$\Sigma$.

To prove the claim let us fix $l>0$ and take a big $k$ so
that $n_k\gg l$.  We note that for any distinct
$x,y\in\gamma^\infty$ there are $x_k,y_k\in\gamma^{\,n_k}$
such that $(x_k,y_k)\to(x,y)$ when $k\to\infty$.  Then for
$x_k,y_k$ there exists a neighborhood $V_{n_k}$ of a point
$\gamma$ such that $\gamma^{\,n_k}=R^{\,n_k}(V_{n_k})$.

We take $j=n_k-l$. 
We can now write for some
$w_k,z_k\in V_{n_k}$
\begin{eqnarray*}
  \label{eq:expandingtimes}
d(x_k,y_k) 
&=&
d\Big( R^{n_k-j}\big( R^{j}(w_k)\big),
R^{n_k-j}\big( R^{j}(z_k)\Big) \ge
\frac{e^{lc/4}}{\beta_0}
\cdot d\big( R^{n_k-l}(w_k), R^{n_k-l}(z_k)\big).
\end{eqnarray*}
Note that each pair $R^{n_k-l}(w_k), R^{n_k-l}(z_k)$ belongs
to a section $\Sigma_k\in\Xi$ and that
$R^l\big(R^{n_k-l}(w_k)\big)=x_k$ and $R^l\big(R^{n_k-l}(z_k)\big)=y_k$.
Letting $k\to\infty$ we obtain limit points $\big(
R^{n_k-l}(w_k), R^{n_k-l}(z_k)\big) \to (w_l,z_l)$ in some
section $\Sigma\in\Xi$ (recall that $\Xi$ is a finite family
of compact adapted cross-sections) satisfying
\[
R^l(w_l)=x,\quad
R^l(z_l)=y\quad\text{and}\quad
d(w_l,z_l )
\le 
\beta_0 e^{-lc/4}\cdot d(x,y).
\]
Since this is true for any $l>0$ we conclude that $y$ is in
the unstable manifold of $x$ with respect to $R$, i.e. $y\in
W^u_R(x)$, thus $y\in W^u(x,\Sigma)$ by the following lemma.
This proves the claim.

\begin{lemma}
  \label{le:WuWsigma}
  In the same setting as above, we have $W^u_R(x)\subseteq
  W^u(x,\Sigma)$.
\end{lemma}

Notice that since both sets $W^u_R(x)$ and $W^u(x,\Sigma)$
are one-dimensional manifolds embedded in a neighborhood of
$x$ in $\Sigma$, then they coincide in a (perhaps smaller)
neighborhood of $x$.

\begin{proof}
  Let $y_0\in W^u(x,\Sigma)$. Then there exists $\epsilon$
  so that $z_0=X_\epsilon(y_0)\in W^{uu}(x)$, with
  $|\epsilon|$ small by Remark~\ref{r.trapaca} and tending
  to $0$ when we take $y_0\to x$. Let $t_l>0$ be such that
  $X_{-t_l}(x)=w_l\in\Sigma$ for $l\ge1$.  Then we have
  \begin{align}\label{eq:t_l0}
    \dist\big(X_{-t_l}(z_0),X_{-t_l}(x)\big)\xrightarrow[l\to\infty]{}0
  \end{align}
  and so there exists $\epsilon_l$ such that
  $X_{\epsilon_l-t_l}(z_0)=z_l=X_{\epsilon_l+\epsilon-t_l}(y_0)\in\Sigma$
  with $|\epsilon_l|$ small. Notice that \eqref{eq:t_l0}
  ensures that $|\epsilon_l|\to0$ also.

  Hence there exists $\delta=\delta(\epsilon,\epsilon_l)$
  satisfying $\delta\to0$ when $(\epsilon+\epsilon_l)\to0$
  and also $d(z_l,w_l)<\delta$ for all $l\ge1$.  Since
  $R^l(z_l)=y_0$ we conclude that $y_0\in W^u_R(x)$,
  finishing the proof.
\end{proof}

\subsubsection{Upper and lower bounds for  densities through
approximation}
\label{sec:bounds-dens-from}

We define $\F^u_\infty$ to be the family of all leaves
$\gamma^\infty$ obtained as $C^1$ accumulation points of
leaves in
\[
\F^u_n=\{ \xi\in\cup_{j=0}^{n-1} D_j :
\xi\quad\mbox{crosses}\quad
F^s(x_0,\de_2)\}.
\]
We note that $\F^u_\infty\subset\F^u(x_0,\de_2)$.  Since for
all $n$ we have $\omega^n_{0}\ge \omega^n_{\rho}$ and so
$\omega_0^n(\cup\F^u_n)>c$, we get that
$\eta\big(\cup\F^u_\infty\big)\ge c$.  By definition of
$\F^u_n$ and by~\eqref{eq:denslimitadahyp} we see that
$\omega^n_{0}\mid F^u_n$ disintegrates along the partition
$\F^u_n$ of $F^u_n=\cup\F_n^u$ into measures $\omega^n_\xi$
having density with respect to $\lambda_\xi$ uniformly
bounded from above and below, for almost every
$\xi\in\F^u_n$.

To take advantage of this in order to prove
Proposition~\ref{pr:densidades} we consider a sequence of
increasing partitions $(\V_k)_{k\ge1}$ of $W^s(x_0,\Sigma)$
whose diameter tends to zero. This defines a sequence
$\cP_k$ of partitions of $\tilde\F=\cup_{0\le n \le
  \infty}\F^u_n$ as follows: we fix $k\ge1$ and say that two
elements $\xi\in\F^u_i,\xi'\in\F^u_j, 0\le i,j\le\infty$ are
in the same atom of $\cP_k$ when both intersect
$W^s(x,\Sigma)$ in the same atom of $\V_k$ and either
$i,j\ge k$ or $i=j<k$.

If $q$ is the projection $q: \tilde\F\to W^s(x_0,\Sigma)$
given by the transversal intersection $\xi\cap
W^s(x_0,\Sigma)$ for all $\xi\in\tilde\F$, then $\tilde\F$
can be identified with a subset of the real line.  Thus we
may assume without loss that the union $\partial\cP_k$ of
the boundaries of $\cP_k$ satisfies
$\eta(\partial\cP_k)=\hat\eta(\partial\cP_k)=0$ for all $k\ge1$, by
suitably choosing the sequence $\V_k$.

\subsubsection*{Upper and lower bounds for  densities}

\begin{figure}[htpb]
  \centering
  \psfrag{S}{$\Sigma$}
  \psfrag{F}{$F^s(x_0,\delta_2)$}
  \psfrag{B}{$B$}
  \psfrag{z}{$\zeta$}
  \psfrag{P}{$\cP_k(\zeta)$}
  \psfrag{p}{$p^{-1}(B)$}
  \includegraphics[width=6cm,height=3cm]{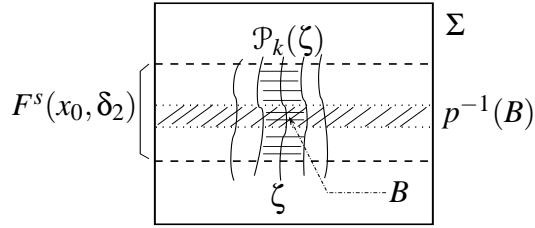}
  \caption{Leaves crossing $F^s(x_0,\delta_2)$ and the projection $p$.}
  \label{fig:leaves-crossing-fsx}
\end{figure}

Given $\zeta\in\tilde\F$ we write $p:F^u(x_0,\de_2)\to
\zeta$ the projection along stable leaves and $\omega$ for
$\omega_0$.  Writing $\cP_k(\zeta)$ for the atom of $\cP_k$
which contains $\zeta$, then since $\cP_k(\zeta)$ is a union
of leaves, for any given Borel set $B\subset\zeta$ and
$n\ge1$
\begin{equation}
  \label{eq:desintegra}
\omega^n\big(\cP_k(\zeta)\cap p^{-1}(B)\big)=
\int \omega^n_\xi\big(\cP_k(\zeta)\cap
p^{-1}(B)\big)\,d\hat\omega^n(\xi)  
\end{equation}
through disintegration, where $\hat\omega^n$ is the measure
on $\tilde\F$ induced by $\omega^n$.  Moreover
by~\eqref{eq:denslimitadahyp} and because each curve in
$\tilde\F$ crosses $F^u(x_0,\de_2)$
\begin{equation}
  \label{eq:desintegra2}
\frac1{\kappa_1\kappa_2}\cdot
\lambda_{\zeta}(B)
\le
\frac1{\kappa_1}\cdot\lambda_\xi\big(p^{-1}(B)\big)
\le
\omega_\xi^n\big(\cP_k(\zeta)\cap p^{-1}(B)\big)
\le
\kappa_1\cdot\lambda_\xi\big(p^{-1}(B)\big)
\le
\kappa_1\kappa_2\cdot
\lambda_{\zeta}(B)  
\end{equation}
for all $n,k\ge1$ and $\hat\omega^n$-a.e.
$\xi\in\tilde\F$, where $\kappa_2>0$ is a constant such that
\[
\frac1{\kappa_2}\cdot\lambda_\zeta\le
\lambda_\xi
\le\kappa_2\cdot\lambda_\zeta
\quad\mbox{for all}\quad
\xi\in\tilde\F,
\]
which exists since the angle between the stable leaves in
any $\Sigma\in\Xi$ and any $cu$-curve is bounded from below,
see Figure~\ref{fig:leaves-crossing-fsx}.

Finally letting $\zeta\in\F^u_\infty$ and choosing $B$ such
that $\eta\big(\partial p^{-1}(B)\big)=0$ (which poses no
restriction), assuming that
$\eta\Big(\partial\big(\cP_k(\zeta)\cap
p^{-1}(B)\big)\Big)=0$ we get from~\eqref{eq:desintegra}
and~\eqref{eq:desintegra2} for all $k\ge1$
\begin{equation}
  \label{eq:RadonNikodym}
\frac1{\kappa_1\kappa_2}\cdot\lambda_\zeta(B)
\cdot\hat\eta\big(\cP_k(\zeta)\big)
\le
\eta\big(\cP_k(\zeta)\cap p^{-1}(B)\big)
\le
\kappa_1\kappa_2\cdot\lambda_\zeta(B)
\cdot\hat\eta\big(\cP_k(\zeta)\big)  
\end{equation}
by the weak$^*$ convergence of $\omega^n$ to $\eta$. Thus to
conclude the proof we are left to check that
$\eta\big(\partial\big(\cP_k(\zeta)\cap
p^{-1}(B)\big)\big)=0$. For this we observe that $\cP_k(\zeta)\cap
p^{-1}(B)$ can be written as the product
$q(\cP_k(\zeta))\times B$. Hence the boundary is equal to
\[
\big(\partial q(\cP_k(\zeta))\times B\big)
 \,\cup\,
\big( q(\cP_k(\zeta))\times \partial B\big)
 \,\subset\,
q^{-1}\big(\partial q(\cP_k(\zeta))\big)\cup
p^{-1}(B)
\]
and the right hand side has $\eta$-zero measure by construction.

This completes the proof of Proposition~\ref{pr:densidades}
since we have $\{\zeta\}=\cap_{k\ge1}\cP_k(\zeta)$ for all
$\zeta\in\tilde\F$ and, by the Theorem of Radon-Nikodym, the
bounds in~\eqref{eq:RadonNikodym} imply that the
disintegration of $\eta\mid\cup\F^u_\infty$ along the curves
$\zeta\in\F^u_\infty$ is absolutely continuous with respect
to Lebesgue measure along these curves and with uniformly
bounded densities from above and from below.

\subsection{The support covers the whole attractor}
\label{sec:support-covers-whole}

\begin{figure}[htpb]
  \centering
\psfrag{a}{$X_t(z)$}
\psfrag{b}{$w$}
\psfrag{c}{$\gamma$}
\psfrag{e}{$z$}
\psfrag{d}{\large piece of stable manifold}
  \includegraphics[width=6.5cm]{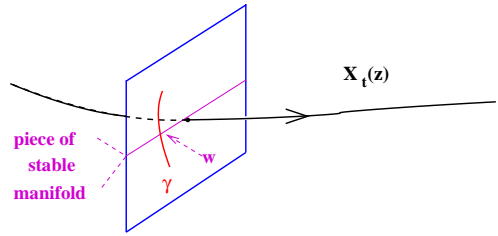}
  \caption{\label{fig:transit-gibbs-proper}Transitiveness
    and support of the physical measure.}
\end{figure}

Finally to conclude that $\suporte(\mu)=\Lambda$ it is
enough to show that $\suporte(\mu)$ contains some $cu$-curve
$\gamma:(a,b)\to\Sigma$ in some subsection $\Sigma\in\Xi$.
Indeed, see Figure~\ref{fig:transit-gibbs-proper}, letting
$x_0\in\Lambda\cap\Sigma$ be a point of a forward dense
regular $X$-orbit and fixing $c\in(a,b)$ and $\epsilon>0$
such that $a<c-\epsilon<c+\epsilon<b$, then for any $\rho>0$
there exists $t>0$ satisfying
$\dist\big(\gamma(c),X_t(x_0)\big)<\rho$.  Since
$W^s\big(X_t(x_0),\Sigma\big)\pitchfork
\big(\gamma\mid(c-\epsilon,c+\epsilon)\big)=\{ z\}$ (because
$\gamma$ is a $cu$-curve in $\Sigma$ and $\rho>0$ can be
made arbitrarily small, where $\pitchfork$ means transversal
intersection), then, by the construction of the adapted
cross-section $\Sigma$ (see
Section~\ref{sec:cross-sect-poinc}), this means that $z\in
W^s\big(X_t(x_0)\big)$. Hence the $\omega$-limit sets of $z$
and $x_0$ are equal to $\Lambda$.  Thus
$\suporte(\mu)\supseteq\Lambda$ because $\suporte(\mu)$ is
$X$-invariant and closed, and
$\Lambda\supseteq\suporte(\mu)$ because $\Lambda$ is an
attracting set.

We now use~\eqref{eq:RadonNikodym} to show that
$\hat\eta$-almost every $\gamma\in\tilde\F$ is contained in
$\suporte(\eta)$, which is contained in $\suporte(\mu)$ by
the construction of $\mu$ from $\eta$ in
Section~\ref{sec:proof-theorem-B}. In fact,
$\hat\eta$-almost every $\zeta\in\tilde\F$ is a density
point of $\hat\eta\mid\tilde\F$ and so for any one $\zeta$
of these curves we have $\hat\eta\big(\cP_k(\zeta)\big)>0$
for all $k\ge1$. Fixing $z\in\zeta$ and choosing
$\epsilon>0$ we may find $k\ge1$ big enough and a small
enough open neighborhood $B$ of $z$ in $\zeta$ such that
\[
\cP_k(\zeta)\cap p^{-1}(B)
\subset B(z,\epsilon)\cap\Sigma
\quad\mbox{and}\quad
\eta\big( \cP_k(\zeta)\cap p^{-1}(B) \big) >0,
\]
by the left hand side inequality in \eqref{eq:RadonNikodym}.
Since $\epsilon>0$ and $z\in\zeta$ where arbitrarily chosen,
this shows that $\zeta\in\suporte(\eta)\subset\suporte(\mu)$
and completes the proof of Theorem~\ref{thm:srbmesmo}.




\def\cprime{$'$}

\end{document}